\documentclass[12pt, reqno]{amsart}

\usepackage{amsmath, amsthm, amscd, amsfonts, amssymb, graphicx, xcolor}

\newtheorem{theorem}{\indent\sc Theorem}[section]
\newtheorem{lemma}[theorem]{\indent\sc Lemma}
\newtheorem{proposition}[theorem]{\indent\sc Proposition}
\newtheorem{corollary}[theorem]{\indent\sc Corollary}
\newtheorem{maintheorem}{\indent\sc Theorem}

\newtheorem{conjecture}[theorem]{\indent\sc Conjecture}

\theoremstyle{definition}
\newtheorem{definition}[theorem]{\indent\sc Definition}
\newtheorem{example}[theorem]{\indent\sc Example}

\theoremstyle{remark}
\newtheorem{remark}[theorem]{\indent\sc Remark}
\numberwithin{equation}{section}

\newcommand{\Broer}[1]{}

\usepackage{latexsym,amssymb}
\usepackage{graphicx}
\usepackage{tikz-cd}
\usepackage{amsfonts}
\usepackage[T1]{fontenc}
\usepackage[utf8]{inputenc}
\usepackage{amsmath}
\usepackage{IEEEtrantools}
\usepackage{amsthm}
\usepackage{euscript}
\usepackage{tikz}
\usepackage{bm}
\usepackage{mathtools}
\usepackage[citestyle=alphabetic,bibstyle=alphabetic]{biblatex}
\usepackage{indentfirst}
\usepackage[british]{babel}
\usepackage{csquotes}
\usepackage{dsfont}
\usepackage{enumitem}
\usepackage{soul}
\usepackage{setspace}
\usepackage{perpage}

\usetikzlibrary{arrows,chains,matrix,positioning,scopes,shapes}
\tikzset{join/.code=\tikzset{after node path={\ifx\tikzchainprevious\pgfutil@empty\else(\tikzchainprevious)edge[every join]#1(\tikzchaincurrent)\fi}}}
\tikzset{>=stealth',every on chain/.append style={join},
         every join/.style={->}}
\newcommand{\Thmstop}{\hglue-6pt.\kern6pt}
\graphicspath{ {./images/} }

\usepackage{caption}
\usepackage{subcaption}

\usepackage[bookmarksnumbered, colorlinks, plainpages=false]{hyperref}
\hypersetup{pdfborder=0 0 0}
\hypersetup{citecolor = blue}
\hypersetup{urlcolor = black}

\usepackage[toc]{appendix}

\usepackage{array,booktabs}
\usepackage{amsmath}
\usepackage{nccmath}
\usepackage{tabularray}
\UseTblrLibrary{amsmath}
\usepackage{nicematrix}

\usepackage{setspace}
\usepackage{comment}

\usepackage[myheadings]{fullpage}

\newenvironment{nohyphens}{
  \par
  \hyphenpenalty=10000
  \exhyphenpenalty=10000
  \sloppy
}{\par}

\DeclareMathSymbol{\mhyphen}{\mathord}{AMSa}{"39}
\newcommand{\itbf}[1]{\textit{\textbf{#1}}}
\newcommand{\slant}[2]{{\raisebox{.08em}{$#1$}\big/\raisebox{-.08em}{$#2$}}}
\newcommand{\restr}{\mathord\downarrow}

\makeatletter
\let\save@mathaccent\mathaccent
\newcommand*\if@single[3]{
  \setbox0\hbox{${\mathaccent"0362{#1}}^H$}
  \setbox2\hbox{${\mathaccent"0362{\kern0pt#1}}^H$}
  \ifdim\ht0=\ht2 #3\else #2\fi
  }
\newcommand*\rel@kern[1]{\kern#1\dimexpr\macc@kerna}
\newcommand*\widebar[1]{\@ifnextchar^{{\wide@bar{#1}{0}}}{\wide@bar{#1}{1}}}
\newcommand*\wide@bar[2]{\if@single{#1}{\wide@bar@{#1}{#2}{1}}{\wide@bar@{#1}{#2}{2}}}
\newcommand*\wide@bar@[3]{
  \begingroup
  \def\mathaccent##1##2{
    \let\mathaccent\save@mathaccent
    \if#32 \let\macc@nucleus\first@char \fi
    \setbox\z@\hbox{$\macc@style{\macc@nucleus}_{}$}
    \setbox\tw@\hbox{$\macc@style{\macc@nucleus}{}_{}$}
    \dimen@\wd\tw@
    \advance\dimen@-\wd\z@
    \divide\dimen@ 3
    \@tempdima\wd\tw@
    \advance\@tempdima-\scriptspace
    \divide\@tempdima 10
    \advance\dimen@-\@tempdima
    \ifdim\dimen@>\z@ \dimen@0pt\fi
    \rel@kern{0.6}\kern-\dimen@
    \if#31
      \overline{\rel@kern{-0.6}\kern\dimen@\macc@nucleus\rel@kern{0.4}\kern\dimen@}%
      \advance\dimen@0.4\dimexpr\macc@kerna
      \let\final@kern#2%
      \ifdim\dimen@<\z@ \let\final@kern1\fi
      \if\final@kern1 \kern-\dimen@\fi
    \else
      \overline{\rel@kern{-0.6}\kern\dimen@#1}
    \fi
  }
  \macc@depth\@ne
  \let\math@bgroup\@empty \let\math@egroup\macc@set@skewchar
  \mathsurround\z@ \frozen@everymath{\mathgroup\macc@group\relax}
  \macc@set@skewchar\relax
  \let\mathaccentV\macc@nested@a
  \if#31
    \macc@nested@a\relax111{#1}
  \else
    \def\gobble@till@marker##1\endmarker{}
    \futurelet\first@char\gobble@till@marker#1\endmarker
    \ifcat\noexpand\first@char A\else
      \def\first@char{}
    \fi
    \macc@nested@a\relax111{\first@char}
  \fi
  \endgroup
}
\makeatother

\newcommand{\s}[0]{_{\mathrm{s}}}
\newcommand{\n}[0]{_{\mathrm{n}}}
\newcommand{\g}[0]{\mathfrak{g}}
\newcommand{\fl}[0]{\mathfrak{l}}
\newcommand{\h}[0]{\mathfrak{h}}
\newcommand{\fu}[0]{\mathfrak{u}}
\newcommand{\p}[0]{\mathfrak{p}}

\newcommand{\z}[0]{\mathfrak{z}}

\newcommand{\Ru}[0]{\mathrm{R_u}}
\newcommand{\cN}[0]{\mathcal{N}}
\newcommand{\Ng}[0]{\mathcal{N}_{\mathfrak{g}}}
\newcommand{\N}[0]{\mathbb{N}}
\newcommand{\Np}[0]{\mathbb{N}^+}
\newcommand{\reg}[0]{\mathrm{reg}}
\newcommand{\Ad}[0]{\mathrm{Ad}}
\newcommand{\ad}[0]{\mathrm{ad}}
\newcommand{\md}[0]{\mathrm{d}}
\newcommand{\GL}[0]{\mathrm{GL}}
\newcommand{\gl}[0]{\mathfrak{gl}}

\newcommand{\midd}[0]{~\middle|~}
\newcommand{\bbK}[0]{\mathbb{K}}
\newcommand{\greg}[0]{{\g\text{-}\reg}}
\newcommand{\U}[0]{\mathrm{U}}
\newcommand{\Gm}[0]{\mathbb{G}_{\mathrm{m}}}

\newcommand{\fD}[0]{\mathfrak{D}}
\newcommand{\bbKx}[0]{\mathbb{K}^{\times}}
\newcommand{\ft}[0]{\mathfrak{t}}
\newcommand{\Le}[0]{\left(L;e_0\right)}
\newcommand{\fJ}[0]{\mathfrak{J}}
\newcommand{\am}[0]{\langle m \rangle}
\newcommand{\ls}[0]{\g_{\am}}
\newcommand{\dvJ}[0]{\overline{\mathfrak{J}}}

\DeclareMathOperator{\fc}{\mathfrak{c}}
\DeclareMathOperator{\C}{\mathrm{C}}
\DeclareMathOperator{\Ind}{\mathrm{Ind}}
\newcommand{\Co}[1]{\operatorname{\mathrm{C}}^{\,\circ}_{#1}}
\DeclareMathOperator{\Lie}{\mathrm{Lie}}

\newcommand{\J}[1]{\operatorname{\mathfrak{J}}_{#1}}
\DeclareMathOperator{\fd}{\mathfrak{d}}
\newcommand{\dreg}[1]{\operatorname{\mathfrak{d}}^{\mathrm{reg}}_{#1}}
\newcommand{\dv}[2]{\overline{\J{#1} #2}}
\newcommand{\llangle}[0]{\left\langle}
\newcommand{\rrangle}[0]{\right\rangle}

\begingroup
\makeatletter
\@ifundefined{ver@biblatex.sty}
  {\@latex@error
     {Missing 'biblatex' package}
     {The bibliography requires the 'biblatex' package.}
      \aftergroup\endinput}
  {}
\endgroup

\makeatletter
\DeclareRobustCommand\subsetneqdot{%
  \mathrel{% <-- or \mathbin, or ...
    \vphantom{\subsetneq}%
    \mathpalette\aw@qc\relax
  }%
}
\newcommand{\aw@qc}[2]{%
  \sbox\z@{$#1\subsetneq$}%
  \sbox\tw@{$#1\bm{\cdot}$}%
  \dimen@=.5\dimexpr\ht\z@-\ht\tw@\relax
  \ooalign{%
    $\m@th#1\subsetneq$\cr
    \hidewidth\raise\dimen@\box\tw@\hidewidth\cr
  }%
}
\makeatother

\newcommand{\boldpreceq}[0]{\mathrel{\bm{\preceq}}}
\newcommand{\boldprec}[0]{\mathrel{\bm{\prec}}}
\newcommand{\boldprecdot}[0]{\mathrel{\bm{\prec}\mathrel{\mkern-5mu}\mathrel{\bm{\cdot}}}}

\newcommand{\smallslant}[2]{{\raisebox{.05em}{$#1$}\big/\raisebox{-.05em}{$#2$}}}

\newcommand{\fm}[0]{\mathfrak{m}}

\newcommand{\bbZ}[0]{\mathbb{Z}}

\usepackage{scalerel,stackengine}

\newcommand\pig[1]{\scalerel*[5.5pt]{\Big#1}{%
  \ensurestackMath{\addstackgap[1.5pt]{\big#1}}}}
\newcommand\pigl[1]{\mathopen{\pig{#1}}}
\newcommand\pigr[1]{\mathclose{\pig{#1}}}

\let\originalcdot\cdot
\renewcommand{\cdot}{\originalcdot\nobreak}

\usepackage{xr} 
\newcommand{\scalemath}[2]{\scalebox{#1}{\mbox{\ensuremath{\displaystyle #2}}}}

\newcommand{\vbigm}[0]{\vphantom{\bigm|}}

\newcommand{\vpigm}[0]{\vphantom{\pig|}}
\usepackage{multirow}
\usepackage[normalem]{ulem}
\useunder{\uline}{\ul}{}
\usepackage{float}
\allowdisplaybreaks

\usepackage{tikz}
\usetikzlibrary{positioning, quotes, backgrounds}
\usetikzlibrary{backgrounds}
\usetikzlibrary{arrows}
\usetikzlibrary{shapes,shapes.geometric,shapes.misc}

\tikzstyle{tikzfig}=[baseline=-0.25em,scale=0.5]

\pgfkeys{/tikz/tikzit fill/.initial=0}
\pgfkeys{/tikz/tikzit draw/.initial=0}
\pgfkeys{/tikz/tikzit shape/.initial=0}
\pgfkeys{/tikz/tikzit category/.initial=0}

\pgfdeclarelayer{edgelayer}
\pgfdeclarelayer{nodelayer}
\pgfsetlayers{background,edgelayer,nodelayer,main}

\tikzstyle{none}=[inner sep=0mm]

\tikzstyle{every loop}=[]

\tikzstyle{red dot}=[fill=red, draw=black, shape=circle]
\tikzstyle{green dot}=[fill={rgb,255: red,0; green,185; blue,0}, draw=black, shape=circle]
\tikzstyle{cyan rectangle}=[fill=cyan, draw=cyan, shape=rectangle]
\tikzstyle{white ellipse}=[fill=white, draw=white, shape=ellipse]
\tikzstyle{white rectangle}=[fill=white, draw=white, shape=rectangle]
\tikzstyle{bold white ellipse}=[fill=white, draw=white, shape=ellipse, font={\boldmath}]
\tikzstyle{large red text white label}=[fill=white, draw=white, shape=ellipse, text=red, font={\boldmath}, font={\large}]
\tikzstyle{blue text white label}=[fill=white, draw=white, shape=ellipse, text=blue, font={\boldmath}]
\tikzstyle{large blue text white label}=[fill=white, draw=white, shape=ellipse, text=blue, font={\boldmath}, font={\large}]
\tikzstyle{bold white ellipse green text}=[fill=white, draw=white, shape=ellipse, font={\boldmath}, text=black!60!green]
\tikzstyle{bold white ellipse blue text}=[fill=white, draw=white, shape=ellipse, font={\boldmath}, text=blue]
\tikzstyle{bold white ellipse red text}=[fill=white, draw=white, shape=ellipse, font={\boldmath}, text=red]
\tikzstyle{bold blue rectangle}=[fill=white, draw=blue, shape=rectangle, line width=2, font={\boldmath}]
\tikzstyle{bold blue rectangle green text}=[fill=white, draw=blue, shape=rectangle, line width=2, font={\boldmath}, text=black!60!green]

\tikzstyle{new edge style 0}=[-, fill=black]
\tikzstyle{red dashed}=[-, dashed, line width=2, draw=red]
\tikzstyle{blue dashed}=[-, dashed, line width=1, draw=blue]
\tikzstyle{Thick Cyan Arrow}=[->, line width=2, draw=cyan]

\begin{document}

\setcounter{page}{1}

\centerline{}

\centerline{}

\title[Lie Algebra Decomposition Classes in Arbitrary Characteristic]{Lie Algebra Decomposition Classes for Reductive Algebraic Groups in Arbitrary Characteristic}

\author[Joel Summerfield]{Joel Summerfield}

\address{School of Mathematics, University of Birmingham, UK.}
\email{\textcolor[rgb]{0.00,0.00,0.84}{jns228@student.bham.ac.uk}}

\begin{abstract}
\begin{nohyphens}
In this paper, we investigate the decomposition classes of the Lie algebras of connected reductive algebraic groups, over algebraically closed fields of arbitrary characteristic.
We extend some results proved previously under restrictions on the characteristic, including a formula for the dimension of a decomposition class, and introduce Levi-type decomposition classes to account for some of the difficulties encountered in bad characteristic.
We also establish properties of Lusztig--Spaltenstein induction of non-nilpotent orbits, such as parabolic independence, extending the known results for nilpotent orbits.
Finally, we determine the covering relation for the closure order on decomposition classes, in the case of good characteristic.
\end{nohyphens}
\end{abstract} \maketitle

\section{Introduction}

In \cite[\S 5.2]{BK79}, Borho and Kraft introduced \emph{Zerlegungsklassen} (decomposition classes) as a tool for studying \emph{Schichten der Lie-Algebra} (sheets of a Lie algebra).
They considered a (connected) semisimple algebraic group $G$ of adjoint type over an algebraically closed field of characteristic $0$, acting via the adjoint action on its Lie algebra $\g = \Lie G$.
For an arbitrary element $x \in \g$, with Jordan decomposition $x = x\s + x\n$, define $\C_Gx\s \coloneqq \left\{ g \in G \midd g \cdot x\s = x\s \right\}$.
Another element $y \in \g$ then has \emph{\"{a}hnliche Jordanzerlegung} (similar Jordan decomposition) if there exists $g \in G$ such that the Jordan decomposition of $g \cdot y = y'\s + y'\n$ satisfies $\C_G y'\s = \C_G x\s$ and $\left(\C_Gx\s\right) \cdot y'\n = \left(\C_Gx\s\right) \cdot x\n$.
This yields an equivalence relation on $\g$, whose corresponding equivalence classes are decomposition classes.
Useful properties of these decomposition classes were then established in later parts of \cite{BK79} and \cite{B81}.

Spaltenstein demonstrated in \cite{S82} that some of the properties from \cite{BK79} generalised immediately to good characteristic; since, in this case, connected stabilisers of semisimple elements of $\g$ are Levi subgroups of $G$ (for a reductive algebraic group $G$).
Moreover, in \cite{S82}, certain properties related to nilpotent orbits were also shown to hold in the classical cases in bad characteristic (see \S\ref{Subsection Levi-Type Sheets}).

In \cite[\S 3]{B98}, Broer considered decomposition classes for the adjoint group $G$ of a semisimple Lie algebra, over an algebraically closed field.
Broer primarily used the additional assumption that the characteristic is very good, and generalised further results from \cite{BK79} relating to the closures of decomposition classes, as well as establishing that decomposition classes are smooth.
Further results on decomposition classes can be found in \cite{B98dv} and \cite{TY05} (both assuming characteristic $0$), \cite{PS18} (assuming the Standard Hypotheses), and \cite{A25} (partly assuming good characteristic). \\

\pagebreak

After establishing some preliminaries, we adapt the definition of decomposition classes from \cite{BK79} to an algebraic group $G$ over an arbitrary algebraically closed field $\bbK$, by replacing the stabiliser $\C_G x\s$ with its identity component $\Co{G}x\s$.
The decomposition class containing $x \in \g$ is then denoted $\J{G}x$, and we denote the set of decomposition classes by $\fD[G]$.

From \S\ref{Subsection Initial Properties} onwards, we assume that $G$ is a connected reductive algebraic group, and establish that our definition of decomposition classes coincides with that of packets used in \cite[\S 1.2]{S82}.
We prove initial properties of decomposition classes in \textsc{Theorem}~\ref{Theorem Initial Decomposition Class Properties}, including that they are $G$-stable, $\bbKx$-stable, irreducible, and constructible sets which form a finite partition of $\g$.
Moreover, we confirm the following formula for the dimension of a decomposition class, where $\fd_{\g}x\s = \z(\fc_{\g}x\s)$ is the centre of the centraliser $\fc_{\g}x\s = \left\{ y \in \g \midd [y,x\s] = 0 \right\}$.

\begin{maintheorem}\label{Maintheorem Dimension of Decomposition Class}
    Suppose $G$ is a connected reductive algebraic group, over an algebraically closed field.
    If $\fJ \in \fD[G]$ is a decomposition class, then $\dim \fJ = \dim G \cdot x + \dim \fd_{\g}x\s$, for any $x \in \fJ$.
\end{maintheorem}

This formula was previously proved by Premet--Stewart in \cite[Proposition 2.5]{PS18}, under certain conditions on $G$ known as the Standard Hypotheses (see \cite[\S 2.9]{J04} for an explanation thereof).
We then equip $\fD[G]$ with the closure order, defined so that $\J{G}x \boldpreceq \J{G}y$ if and only if $\J{G}x \subseteq \dv{G}{y}$.
In \S\ref{Section Preservation of Decomposition Classes}, we explore how the structure of decomposition classes is affected by central surjections.
In particular, we prove that separable central surjections preserve the closure order and interact uniformly with the dimensions of decomposition classes.

We define \emph{level sets} of $\g$ as fibres of either the stabiliser dimension map $\dim \C_G \colon \g \rightarrow\nobreak \N$ or the centraliser dimension map $\dim \fc_{\g} \colon \g \rightarrow \N$.
In \S\ref{Section Sheets}, we generalise known results about (stabiliser) sheets and investigate properties of the closure order when restricted to level sets.
We then describe the unique minimal and maximal decomposition classes, with respect to the closure order.

In \S\ref{Section Lusztig--Spaltenstein Induction}, we look at the generalisation of Lusztig--Spaltenstein induction to arbitrary orbits, building upon the work in \cite{S82}.
For any Levi subgroup $L \subseteq G$ and adjoint $L$-orbit $\mathcal{O}$, we confirm the existence of an \emph{induced orbit} $\Ind_{\fl}^{\g} \mathcal{O}$.
As with the nilpotent case, the construction of $\Ind_{\fl}^{\g} \mathcal{O}$ appears to depend on a choice of parabolic subgroup $P \subseteq G$ for which $L \subseteq P$ is a Levi factor.
However, using an outline provided by \cite[\S 2.1]{S82}, we provide a complete proof that this is not the case.
We also extend the following known results about nilpotent Lusztig--Spaltenstein induction to this more general setting.

\begin{maintheorem}\label{Theorem LS Induction Main Result}
    Suppose $G$ is a connected reductive algebraic group over an algebraically closed field, $L \subseteq G$ is a Levi subgroup, and $\mathcal{O}$ is an $L$-orbit.
\begin{itemize}
    \item[$\mathrm{(i)}$] The induced orbit $\Ind_{\fl}^{\g} \mathcal{O}$ is independent of the choice of parabolic used in its construction. 
    \item[$\mathrm{(ii)}$] Induction is transitive$:$ $\Ind_{\fl}^{\g} \mathcal{O} = \Ind_{\mathfrak{m}}^{\g} \Ind_{\fl}^{\mathfrak{m}} \mathcal{O}$, for Levi subgroups $L \subseteq\nobreak M \subseteq\nobreak G$.
    \item[$\mathrm{(iii)}$] Induction preserves codimension$:$ $\dim \Ind_{\fl}^{\g} \mathcal{O} = \dim \mathcal{O} + (\dim G - \dim L)$.
    \item[$\mathrm{(iv)}$] $\left(\Ind_{\fl}^{\g} \mathcal{O}\right) \cap \left(\mathcal{O} + \fu_{\p}\right)$ is a single $P$-orbit, where $P \subseteq G$ is any parabolic subgroup for which $L$ is a Levi factor, and $\fu_{\p} = \Lie \bigl(\Ru(P)\bigr)$.
\end{itemize}
\end{maintheorem}

Having worked in full generality up to this point, we narrow our scope in \S\ref{Section Levi-Type Decomposition Classes} to \mbox{\emph{Levi-type}} decomposition classes, which are defined as the decomposition classes of elements $x \in \g$ such that $\Co{G}x\s \subseteq G$ is a Levi subgroup.
They are all of the form $\J{G}(L;e_0) \coloneqq G \cdot \left(\mathfrak{z}(\mathfrak{l})^{\mathrm{reg}} + e_0\right)$, where $L \subseteq G$ is a Levi subgroup, $\mathfrak{z}(\mathfrak{l})^{\mathrm{reg}}$ denotes the set of elements $y$ of the Lie algebra centre $\mathfrak{z}(\mathfrak{l})$ for which $\dim \C_G y$ is minimal, and $e_0 \in \mathcal{N}_{\mathfrak{l}}$ is a nilpotent element of $\mathfrak{l}$.
These are introduced as a tool to avoid certain complications that arise in bad characteristic, and will allow us to prove the following (extending prior results in \cite{B81}, \cite{B98}, and \cite{A25}).

\begin{maintheorem}\label{Theorem Levi-Type Decomposition Varieties Main Result}
    Suppose $\J{G}\Le$ is a Levi-type decomposition class for a connected reductive algebraic group $G$, over an algebraically closed field.
    Let $P \subseteq G$ be a parabolic subgroup with Levi factor $L \subseteq P$, and unipotent radical $U_P = \mathrm{R_u}(P)$.
\begin{itemize}
    \item[$\mathrm{(i)}$] $\dv{G}{\Le} = G \cdot \left( \mathfrak{z}(\mathfrak{l}) + \overline{L \cdot e_0} + \mathfrak{u}_{\p}\right) = \bigcup\limits_{z \in \mathfrak{z}(\mathfrak{l})} \overline{\mathrm{Ind}_{\mathfrak{l}}^{\g}L \cdot \left(z + e_0\right)}$.
    \item[$\mathrm{(ii)}$] $\dv{G}{\Le}$ is a union of decomposition classes.
    \item[$\mathrm{(iii)}$] $\dv{G}{\Le}^{\,\reg} = \bigcup\limits_{z \in \mathfrak{z}(\mathfrak{l})} \mathrm{Ind}_{\mathfrak{l}}^{\g}L \cdot \left(z + e_0\right)$.
\end{itemize}
\end{maintheorem}

We then consider a conjecture of Spaltenstein (see \textsc{Conjecture}~\ref{Conjecture Spaltenstein Sheets Contain Unique Nilpotent}) regarding stabiliser sheets and nilpotent orbits.
In particular, we show that (regardless of characteristic) every Levi-type stabiliser sheet contains a unique nilpotent orbit.

Finally, we describe the Hasse diagram of $\fD[G]$ with respect to the closure order, which is defined as the graph $\bm{\Gamma}\fD[G]$ with vertex set $\fD[G]$ and a directed edge for each covering pair.
We equip the set $\mathfrak{C}_G \coloneqq \left\{ \Co{G}x\s \midd x \in \g \right\}$ with the covering relation $\subsetneqdot$ defined by inclusion of sets
For each $M \in \mathfrak{C}_G$, we equip its set of nilpotent orbits $\smallslant{\cN_{\fm}}{M}$ with the covering relation $\boldprecdot$ defined by its closure order.
Under the restriction to good characteristic, we then prove that every edge of $\bm{\Gamma}\fD[G]$ is either the result of a minimal Lusztig--Spaltenstein induction step (case $\mathrm{(I)}$ below), or the result of a minimal covering of nilpotent orbits (case $\mathrm{(II)}$ below).

\begin{maintheorem}\label{Maintheorem Hasse Diagram Edges}
    Suppose $G$ is a connected reductive algebraic group, over an algebraically closed field whose characteristic $p \geq 0$ is good for $G$.
    Then every edge of the Hasse diagram of $\fD[G]$ is exactly one of the following types\emph{:}
    \begin{itemize}
        \item[$\mathrm{(I)}$] An edge from $\J{G}(L';\Ind_{\fl}^{\fl'}\mathcal{O})$ to $\J{G}(L;\mathcal{O})$, for $L,L' \in \mathfrak{C}_G$ with $L \subsetneqdot L'$ and a nilpotent $L$-orbit $\mathcal{O} \in \smallslant{\cN_{\fl}}{L}$.
        \item[$\mathrm{(II)}$] An edge from $\J{G}(L;\mathcal{O}')$ to $\J{G}(L;\mathcal{O})$, for $L \in \mathfrak{C}_G$ and nilpotent $L$-orbits $\mathcal{O},\mathcal{O}' \in \smallslant{\cN_{\fl}}{L}$ with $\mathcal{O}' \boldprecdot \mathcal{O}$.
    \end{itemize}
\end{maintheorem}

This result characterises the covering relation for the closure order on decomposition classes, at least for the good characteristic case.
Decomposition varieties (closures of decomposition classes) were first described as unions of orbit closures in \cite[Proposition 3.1]{B81}, for the characteristic $0$ case, allowing one to determine the closure order on $\fD[G]$.
This was later extended to good characteristic \cite{A25}, and one can use \textsc{Theorem}~\ref{Theorem Levi-Type Decomposition Varieties Main Result} for Levi-type decomposition varieties in arbitrary characteristic.
However, determining the covering relation is a more intricate problem, hence the careful consideration required to deduce \textsc{Theorem}~\ref{Maintheorem Hasse Diagram Edges}.

\subsection*{Acknowledgements}

The author thanks their PhD supervisor, Simon Goodwin, for their continued guidance and support, as well as Matthew Westaway for originally introducing them to this topic.
The author additionally thanks Alexander Fr\"{u}h and Lauren Hazel for assisting in the translation of \cite{BK79}, \cite{B81}, and \cite{B81Uni}.
The author also thanks Simon Goodwin, Lauren Hazel, Lewis Groves, and Lewis Topley for all of their feedback and proofreading, as well as the reviewers for their comments.
The author's research is supported by EPSRC grant EP/V520275/1.

\section{Preliminaries}

We begin by establishing some notation, and recalling general facts from algebraic geometry and topology that will be used later.

\subsection{Notation}

Throughout, $\bbK$ is an algebraically closed field of characteristic $p \geq 0$, with non-zero elements $\bbK^{\times}$.
All varieties and vector spaces are over $\bbK$, and all spaces are equipped with the Zariski topology.
All algebraic groups are affine (hence also linear), and we identify an algebraic group with its set of $\bbK$-points.
The Lie algebras of algebraic groups will be denoted by the corresponding lowercase fraktur letter.

We denote the set of all non-negative integers by $\N$, and the set of strictly positive integers by $\N^+$.
For each $n \in \N^+$, we let $\GL_n$ denote the group of $n \times n$ invertible matrices, and $\gl_n$ the Lie algebra of all $n \times n$ matrices, both with entries in $\bbK$.
For each $r \in \mathbb{N}^+$, we define $\mathbb{K}^r_{\neq} \coloneqq \left\{ (a_1,\ldots,a_r) \in \mathbb{K}^r \midd \text{for all } 1 \leq i \neq j \leq r \colon a_i \neq a_j \right\}$.

Fix an algebraic group $G$, and let $[-,-] \colon \g \times \g \rightarrow \g$ denote the Lie bracket on its Lie algebra $\g$.
Suppose $H \subseteq G$ is a closed subgroup, and $X \subseteq \g$ is an arbitrary subset.
The \emph{identity component} of $H$, denoted $H^{\circ}$, is defined as the connected component of $H$ containing the identity element.
The centre of $H$ is denoted $\mathcal{Z}(H)$, and its identity component is also denoted $\mathcal{Z}^{\circ}(H) = \left(\mathcal{Z}(H)\right)^{\circ}$.

The adjoint action of $H$ on $\g$ is denoted $h \cdot x \coloneqq \Ad(h)(x)$, for any $h \in H$ and $x \in \g$, and the corresponding $H$-\emph{orbit} is $H \cdot x \coloneqq \left\{ h \cdot x \midd h \in H \right\}$.
The set of all adjoint $H$-orbits in $\g$ is denoted $\slant{\g}{H}$.
More generally, $H \cdot X \coloneqq \bigcup_{x \in X} H \cdot x$ denotes the $H$-\emph{saturation} of $X$, and we say that $X$ is $H$-\emph{stable} if $H \cdot X = X$.

For any $x \in \g$, let $x = x\s +x\n$ be the (\emph{additive}) \emph{Jordan decomposition}, as explained in \cite[\S 4.4.19]{S98}.
Then $x \in \g$ is \emph{semisimple} if and only if $x = x\s$, and \emph{nilpotent} if and only if $x = x\n$.
The set of all nilpotent elements of $\g$ (the \emph{nilpotent cone}) is denoted $\cN_{\g}$, which is a closed subset of $\g$.

\subsection{Stabilisers and Centralisers}

For each $x \in \g$, we define its $H$\mbox{-\emph{stabiliser}} as the closed subgroup $\C_H x \coloneqq \left\{ g \in H \midd g \cdot x = x \right\}$, and its $\h$-\emph{centraliser} as the subspace $\fc_{\h} x \coloneqq \left\{ y \in \h \midd [x,y] = 0 \right\}$.
More generally, $\C_H X \coloneqq \bigcap_{x \in X} \C_H x$ and $\fc_{\h} X \coloneqq \bigcap_{x \in X} \fc_{\h} x$.
We also let $\z(\h) \coloneqq \fc_{\h} \h = \left\{ x \in \h \midd [x,y] = 0, \text{~for~all~} y \in \h \right\}$ denote the \emph{centre} of $\h$.

When there is no ambiguity, we refer to the $G$-stabiliser and $\g$-centraliser as simply the \emph{stabiliser} and \emph{centraliser}, respectively.
Let $\Co{G} x \coloneqq \left(\C_G x\right)^{\circ}$ denote the \emph{connected stabiliser} of $x$.
The \emph{double centraliser} of $x \in \g$ is defined as $\fd_{\g} x \coloneqq \fc_{\g}\left(\fc_{\g} x\right)$, from which it follows readily that $\fd_{\g} x = \left\{ y \in \g \midd \fc_{\g} x \subseteq \fc_{\g} y \right\} = \z\left(\fc_{\g} x\right)$ (that is, the double centraliser coincides with the centre of the centraliser).
We observe that $g \cdot \fc_{\g} x = \fc_{\g}\left(g \cdot x\right)$, for any $g \in G$.

The following result records the well-known interactions between stabilisers/centralisers and the Jordan decomposition.
The proof of part (i) is a direct consequence of the above definitions, and that of (ii) can be shown using immersive representations alongside \cite[Proposition 4.2(2)]{B91}.

\begin{lemma}\label{Lemma Stabiliser/Centraliser of Jordan Decomposition}
    Suppose $x \in \g$.
\begin{itemize}
    \item[$\mathrm{(i)}$] $\C_G x = \C_G x\s \cap \C_G x\n = \C_{\mathrm{C}_G^{\,\phantom{\circ}} x^{\vphantom{\circ}}\s} x\n$, and $\Co{G} x = \left(\Co{G} x\s \cap \C_G x\n\right)^{\circ} = \Co{\mathrm{C}^{\,\circ}_{G}\,x^{\vphantom{\circ}}\s} x\n$.
    \item[$\mathrm{(ii)}$] $\fc_{\g} x = \fc_{\g} x\s \cap \fc_{\g} x\n = \fc_{\mathfrak{c}^{\vphantom{\circ}}_{\g}x^{\vphantom{\circ}}\s} x\n$.  
\end{itemize}
\end{lemma}

We define the \emph{stabiliser dimension map} $\dim \C_G \colon \g \rightarrow \N$ via $x \mapsto \dim\left(\C_G x \right) \in \N$; likewise for the \emph{centraliser dimension map} $\dim \fc_{\g} \colon \g \rightarrow \N$.
Both of these are constant on each $G$-orbit.
Using the version of Chevalley's Semi-Continuity Theorem from \cite[Corollary AG10.3]{B91}, we can establish the following lemma, in which a map $f \colon \g \rightarrow \N$ is \emph{upper semi-continuous} if $\left\{ x \in \g \midd f(x) \geq n \right\}$ is closed for all $n \in \N$.

\begin{lemma}\label{Lemma Stabiliser/Centraliser Dimension Maps are USC}
    Both the stabiliser and centraliser dimension maps $\dim \C_G \colon \g \rightarrow \N$ and $\dim \fc_{\g} \colon \g \rightarrow \N$ are upper semi-continuous.
\end{lemma}

We then define the \emph{stabiliser level sets} of $\g$ as the fibres of the stabiliser dimension map, and denote them $\g_{(m)} \coloneqq \left\{ x \in \g \midd \dim \C_G x = m \right\}$, for each $m \in \N$.
Analogously, we define the \emph{centraliser level sets} of $\g$ as the fibres of the centraliser dimension map, and denote them $\g_{[m]} \coloneqq \left\{ x \in \g \midd \dim \fc_{\g} x = m \right\}$, for each $m \in \N$.

We use the term \emph{level set} of $\g$ to refer collectively to any subset of $\g$ which is (at least one of) a stabiliser level set or a centraliser level set, and use $\ls$ to denote either $\g_{(m)}$ or $\g_{[m]}$.
It is clear that level sets are $\bbKx$-stable, and it follows from \textsc{Lemma}~\ref{Lemma Stabiliser/Centraliser Dimension Maps are USC} that each level set is locally closed in $\g$.
More generally, for any subspace $V \subseteq \g$, we define $V_{\am} \coloneqq V \cap \g_{\am}$, and observe that it is also locally closed in $\g$. \\

Let $X \subseteq \g$ be an arbitrary subset.
We define the set of $G$-\emph{regular elements} of $X$ as $X^{G\mhyphen\reg} \coloneqq \left\{ x \in X \midd \dim \C_G x \leq \dim \C_G y, \text{for~all~} y \in X \right\}$, the set of elements of $X$ with minimal stabiliser dimension.
Whenever the underlying group is unambiguous, we shall denote this set $X^\reg$ instead.
Moreover, for any connected reductive subgroup $H \subseteq G$, we use $X^\reg$ to refer to $X^{G\mhyphen\reg}$ instead of $X^{H\mhyphen\reg}$.

We define the set of $\g$-\emph{regular elements} of $X$ as the subset of elements with minimal centraliser dimension, and we denote it $X^\greg \coloneqq \left\{ x \in X \midd \dim \fc_{\g} x \leq \dim \fc_{\g} y, \text{for~all~} y \in X \right\}$.
Since $X^\reg = X \cap \g_{(m)}$, where $m \in \N$ is minimal such that this intersection is non-empty, it follows that $X^\reg$ is open in $X$; a similar argument holds for $X^\greg$.

If $V \subseteq \g$ is a subspace, then it follows that both $V^\reg$ and $V^\greg$ are open dense irreducible subsets of $V$, and are thus both irreducible and locally closed in $\g$.
In particular, since $\fd_{\g} x \subseteq \g$ is a subspace (for any $x \in \g$), we know that $\left(\fd_{\g} x\right)^\greg = \left\{ y \in \g \midd \fc_{\g} y = \fc_{\g} x \right\}$ is irreducible and locally closed in $\g$.

\begin{lemma}\label{Lemma Closure inside a Level Set}
    ~
\begin{itemize}
    \item[$\mathrm{(a)}$] If $Y \subseteq \g_{(m)}$, then $Y \subseteq \overline{Y}^{\,\reg} = \overline{Y} \cap \g_{(m)}$. 
    \item[$\mathrm{(b)}$] If $Y \subseteq \g_{[m]}$, then $Y \subseteq \overline{Y}^{\,\greg} = \overline{Y} \cap \g_{[m]}$. 
\end{itemize}
\begin{proof}
    If $Y \subseteq \g_{(m)}$, then \textsc{Lemma}~\ref{Lemma Stabiliser/Centraliser Dimension Maps are USC} implies that $\overline{Y} \subseteq \overline{\g_{(m)}} \subseteq \bigsqcup_{n \geq m} \g_{(n)}$.
    Therefore, the minimal $k \in \N$ such that $\overline{Y} \cap \g_{(k)} \neq \emptyset$ must be $m$, and hence $\overline{Y}^{\,\reg} = \overline{Y} \cap \g_{(m)}$.
    The inclusion $Y \subseteq \overline{Y}^{\,\reg}$ is then immediate since $Y \subseteq \overline{Y} \cap \g_{(m)}$.
    This proves (a), and (b) follows by an analogous argument.
\end{proof}
\end{lemma}

By considering the differential of the $G$-orbit map $g \mapsto g \cdot x$, for some fixed $x \in \g$, we have the inclusion $\Lie\left(\C_G x\right) \subseteq \fc_{\g} x$.
Using \cite[\S 9.1]{B91}, this is an equality if and only if the orbit map is a separable morphism of affine varieties $G \rightarrow \g$.
For example, we have equality under the Standard Hypotheses (see \cite[\S 2.9]{J04}).
On the other hand, \cite[Proposition 9.1(2)]{B91} shows that we have equality whenever $x \in \g$ is semisimple.
Consequently, if $X \subseteq \g$ consists of semisimple elements, then $X^\reg = X^\greg$.

\subsection{Closure Results}

Suppose that $V = \bigoplus V_i$ is a vector space, decomposed as a direct sum of finitely many subspaces, and let $X_i \subseteq V_i$ be a collection of arbitrary non-empty subsets.
Each $V_i \subseteq V$ is a closed irreducible subset, and so $\overline{X_i} \subseteq V_i$.
A simple induction argument then shows that $\overline{\sum X_i} = \sum \overline{X_i}$.

\begin{lemma}\label{Lemma Closure of Saturated Set}
    Suppose $\eta \colon V \rightarrow W$ is a linear map between vector spaces, and $X \subseteq V$ is such that $X = X + \ker \eta$ $($that is $X$ is invariant under addition by elements of $\ker \eta)$.
    Then $\eta\bigl(\:\!\overline{X}\:\!\bigr) = \overline{\eta(X)}$.
\end{lemma}

This result can be readily deduced from properties of (topological) quotient maps.
In particular, it demonstrates that linear injections between vector spaces preserve closures.

We now recall a crucial property of parabolic subgroups, as found in \cite[Proposition 0.15]{H95conj}.
Suppose $G$ is connected, and $P \subseteq G$ is a parabolic subgroup.
If $X$ is a $G$-variety, and $Y \subseteq X$ is a closed $P$-stable subset, then $G \cdot Y \subseteq X$ is also closed.
The following result encapsulates how we will make use of this property henceforth.

\begin{lemma}\label{Lemma Parabolic Property Result}
    Suppose $P \subseteq G$ is a parabolic subgroup of a connected algebraic group, and $X \subseteq \g$ is any subset.
    Then $\overline{G \cdot X} = G \cdot \bigl(\:\! \overline{P \cdot X} \:\!\bigr)$.
\end{lemma}

\subsection{Connected Reductive Algebraic Groups}

Suppose, for now, that $G$ is a connected reductive algebraic group.
We say that a subgroup $H \subseteq G$ is a \emph{connected reductive regular subgroup} if it is a connected reductive closed subgroup of $G$ which contains a maximal torus of $G$.

We use the term \emph{Levi subgroup} to mean a Levi factor of a parabolic subgroup, and observe that all Levi subgroups are connected reductive regular subgroups.
If $L \subseteq G$ is a Levi subgroup, then we let $\mathfrak{P}(G,L)$ denote the (finite) set of all parabolic subgroups of $G$ for which $L$ is a Levi factor.
Given any parabolic subgroup $P \subseteq G$, we let $U_P \coloneqq \Ru(P)$ denote its unipotent radical, with corresponding Lie algebra $\fu_{\p} \subseteq \p$.

For a fixed choice of maximal torus $T \subseteq G$, let $\Phi = \Phi(G,T)$ denote the corresponding \emph{root system}.
Moreover, let $\mathrm{X}(T)$ denote its \emph{character group} and let $\mathrm{Y}(T)$ denote its \emph{cocharacter group}.
Observe that $\mathfrak{P}(G,T)$ is precisely the set of Borel subgroups of $G$ which contain $T$.
For each $\alpha \in \Phi$, we let $\g_{\alpha}$ and $\U_{\alpha}$ denote the corresponding \emph{root subspace} and \emph{root subgroup}, respectively.
For any subset of roots $\Psi \subseteq \Phi$, we let $\g(\Psi)$ be shorthand for the subspace $\bigoplus_{\alpha \in \Psi} \g_{\alpha} \subseteq \g$.

If $H \subseteq G$ is a connected reductive regular subgroup, and $T \subseteq H$, then we identify its root system $\Phi(H,T)$ with the corresponding subset of $\Phi$.
For any $y \in \ft$, let $\Phi_y \coloneqq \left\{ \alpha \in \Phi \midd \md \alpha (y) = 0 \right\}$ denote the set of roots $\alpha \colon T \rightarrow \Gm$ whose differential $\md \alpha \colon \ft \rightarrow \bbK$ has kernel containing $y$.
Then \cite[Lemma 3.7]{S75} shows that $\Co{G} y = \llangle T, \U_{\alpha} \midd \alpha \in \Phi_y \rrangle$ is a connected reductive regular subgroup with root system $\Phi\left(\Co{G} y,T\right) = \Phi_y$.

We note that if $\alpha \in \Phi$ satisfies $\mathrm{d}\alpha = 0$, then $\alpha \in \Phi_y$ for all $y \in \mathfrak{t}$.
Conversely, we shall see in the proof of \textsc{Theorem}~\ref{Theorem Maximal Decomposition Class} that the elements $y \in \mathfrak{t}$ with $\Phi_y = \left\{ \alpha \in \Phi \midd \mathrm{d}\alpha = 0 \right\}$ are precisely the regular semisimple elements of $\mathfrak{g}$, up to $G$-conjugation.

\begin{corollary}\label{Corollary Semisimple Centraliser Description}
    For each $y \in \ft$, we have that $\fc_{\g} y = \ft \oplus \g(\Phi_y)$, which moreover contains only finitely many nilpotent $\Co{G} y$-orbits.
\begin{proof}
    Since $\Co{G} y \subseteq G$ is a closed subgroup containing $T$, we can apply \cite[Proposition 13.20]{B91} to \cite[Lemma 3.7]{S75} to deduce that $\fc_{\g} y = \ft \oplus \g(\Phi_y)$.
    The rest of the statement is a consequence of the fact that a connected reductive group has only finitely many nilpotent orbits in its Lie algebra (see \cite[\S 2.8, Theorem 1]{J04}, for example).
\end{proof}
\end{corollary}

The centralisers of semisimple Lie algebra elements have been further explored in some specific cases---see the proof of \cite[Theorem 2.8]{S82}, for example.
However, we were unable to find any general descriptions in the literature (even for $G$ simple), hence our decision to deduce \textsc{Corollary}~\ref{Corollary Semisimple Centraliser Description} here.
See \textsc{Appendix}~\ref{Appendix: Example} for a more detailed description when $G = \mathrm{Sp}_4$.
We also have the following property of connected reductive regular subgroups which will be used later.

\begin{proposition}\label{Proposition Property of Connected Reductive Regular Subgroups}
    Suppose $H \subseteq G$ is a connected reductive regular subgroup, and let $y \in \g$.
    Then $y \in \z(\h)_{[\dim \h]}$ if and only if $\h = \fc_{\g} y$, if and only if $H = \Co{G}y$.
\begin{proof}
    Since $\z(\h)_{[\dim \h]} = \left\{ z \in \z(\h) \midd \dim \fc_{\g} z = \dim \h \right\}$, the first equivalence readily follows.
    Fix a maximal torus $T \subseteq H$.
    Using \cite[\mbox{Proposition} 7.14]{MT11} and \cite[Corollary 8.13(b)]{MT11}, we can determine that $\left\{ x \in \g \midd T \subseteq \C_G x \right\} = \ft$.
    If $H = \Co{G}y$, then $T \subseteq \C_G y$ implies that $y \in \ft$, and thus $y$ is semisimple.
    On the other hand, if $\h = \fc_{\g}y$, then $y \in \z(\h) \subseteq \ft$, and thus $y$ is again semisimple.
    The second equivalence then follows from \cite[Lemma 3.7]{S75} and \textsc{Corollary}~\ref{Corollary Semisimple Centraliser Description}.
\end{proof}
\end{proposition}

It follows from \textsc{Corollary}~\ref{Corollary Semisimple Centraliser Description} and \textsc{Proposition}~\ref{Proposition Property of Connected Reductive Regular Subgroups} that $\Co{G}x \subseteq G$ is a connected reductive regular subgroup if and only if $x = x\s \in \g$ is semisimple.
We also have the following result that utilises \textsc{Lemma}~\ref{Lemma Parabolic Property Result}.

\begin{proposition}\label{Proposition Closure of Semisimple + Nilpotent}
    If $Z \subseteq \z(\g)$, and $Y \subseteq \Ng$ is $G$-stable, then $\overline{Z+Y} = \overline{Z} + \overline{Y}$.
\begin{proof}
    Suppose $X_1,X_2 \subseteq \g$ are arbitrary, and for each $x_1 \in X_1$, consider the homeomorphism $\g \rightarrow \g$ defined by $x \mapsto x_1 + x$.
    Since homeomorphisms preserve closures, we have that $x_1 + \overline{X_2} = \overline{x_1 + X_2}$.
    This implies that $x_1 + \overline{X_2} \subseteq \overline{X_1 + X_2}$, and thus $X_1 + \overline{X_2} \subseteq \overline{X_1 + X_2}$.
    Using the same argument, for each $x_2 \in \overline{X_2}$, we have that $\overline{X_1} + x_2 = \overline{X_1 + x_2} \subseteq \overline{X_1 + \overline{X_2}} \subseteq \overline{X_1 + X_2}$.
    Therefore, $\overline{X_1} + \overline{X_2} \subseteq \overline{X_1 + X_2}$, for any subsets $X_1,X_2 \subseteq \g$.

    For the reverse inclusion, fix a Borel subgroup $B \subseteq G$ with unipotent radical $U \coloneqq \Ru(B)$, and define $Y_0 \coloneqq Y \cap \fu$.
    Using the proof of \cite[Proposition 2.7(a)]{J04}, we know that $\Ng = G \cdot \fu$.
    Then $B \cdot (Z + Y_0) = Z + Y_0$ and $G \cdot (Z + Y_0) = Z + G \cdot Y_0 = Z + Y$.
    Moreover, $Z \subseteq \z(\g)$ and $Y_0 \subseteq \fu$, thus $\overline{Z + Y_0} = \overline{Z} + \overline{Y_0} \subseteq \z(\g) \oplus \fu$.
    Applying \textsc{Lemma}~\ref{Lemma Parabolic Property Result} with $P = B$ and $X = Z + Y_0$ shows that $\overline{Z + Y} = G \cdot \bigl(\;\!\overline{Z} + \overline{Y_0} \;\!\bigr) = \overline{Z} + G \cdot \overline{Y_0} \subseteq \overline{Z} + G \cdot \bigl(\;\!\overline{Y} \cap \fu\bigr) = \overline{Z} + \overline{Y}$.
\end{proof}
\end{proposition}

We note that \textsc{Proposition}~\ref{Proposition Closure of Semisimple + Nilpotent} is much simpler to prove if $\g = \z(\g) \oplus \Lie (G,G)$, where $(G,G)$ denotes the derived subgroup of $G$; see \cite[Corollary 2.3.9]{L05} for a sufficient condition on $p \geq 0$ for this to hold.

\subsection{Partial Orders and Hasse Diagrams}\label{Subsection Partial Orders and Hasse Diagrams}

Recall that a \emph{partial order} on a set $\Omega$ is a reflexive antisymmetric transitive binary relation $\boldpreceq$.
The pair $(\Omega,\boldpreceq)$ is then also referred to as a \emph{poset}.
We then form the corresponding strict partial order $\boldprec$, defined so that $x \boldprec y$ if and only if $x \boldpreceq y$ and $x \neq y$; likewise the corresponding covering relation $\boldprecdot$, defined so that $x \boldprecdot y$ if and only if $x \boldprec y$ and there does not exist $z \in X$ with $x \boldprec z \boldprec y$.

Given a set $\Omega$ equipped with a partial order $\boldpreceq$, we define its \emph{Hasse diagram} $\bm{\Gamma}\Omega$ as the directed graph with vertex set $\Omega$, and an edge from $x$ to $y$ whenever $x \boldprecdot y$.

An element $x \in \Omega$ is \emph{maximal} (\emph{in} $\Omega$) if there does not exist $y \in \Omega$ with $x \boldprec y$.
Similarly, an element $x \in \Omega$ is \emph{minimal} (\emph{in} $\Omega$) if there does not exist $y \in \Omega$ with $y \boldprec x$.
An element $x \in \Omega$ is \emph{isolated} (\emph{in} $\Omega$) if it is both maximal and minimal in $\Omega$.

\subsection{Closure Orders}\label{Subsection Closure Orders}

Fix a variety $Y$, and recall that a subset $A \subseteq Y$ is \emph{constructible} if it is a finite union of locally closed subsets of $Y$.
Using \cite[Corollary AG10.2]{B91}, each constructible set $A \subseteq Y$ contains a subset $U \subseteq A$ which is open and dense in $\overline{A}$.
Moreover, $A$ is locally closed if and only if $U$ can be taken to be all of $A$.

Suppose that $\mathfrak{X} \coloneqq \left\{ X_1, \ldots, X_r \right\}$ is a finite collection of irreducible constructible subsets of $Y$, which are all pairwise disjoint, and let $X \coloneqq \bigsqcup_i X_i$.
We equip $\mathfrak{X}$ with the \emph{closure order}, defined so that $X_i \boldpreceq X_j$ if and only if $X_i \subseteq \overline{X_j}$.

\begin{lemma}\label{Lemma Closure Order is a Partial Order}
    The closure order $\boldpreceq$ is a partial order on $\mathfrak{X}$.
\begin{proof}
    Since reflexivity and transitivity are immediate from its definition, it remains to be proven that $\boldpreceq$ is antisymmetric, so suppose that $X_i \boldpreceq X_j$ and $X_j \boldpreceq X_i$.
    Let $U_i \subseteq X_i$ be open and dense in $\overline{X_i}$, and likewise for $U_j \subseteq X_j$.
    Since $\overline{X_i} = \overline{X_j}$, we have $U_i \cap U_j \neq \emptyset$, and thus $X_i \cap X_j \neq \emptyset$.
    Therefore, $X_i = X_j$, as required.
\end{proof}
\end{lemma}

For example, this situation arises when $\mathfrak{X} = \smallslant{\cN_{\g}}{G}$ is the set of nilpotent orbits (of a connected reductive algebraic group $G$).
We shall see in \S\ref{Subsection The Closure Order on Decomposition Classes} that this also occurs when $\mathfrak{X}$ is the set of $G$-decomposition classes.

The following corollary of \textsc{Lemma}~\ref{Lemma Closure Order is a Partial Order} follows immediately from \cite[Proposition 1.22]{MT11}.

\begin{corollary}\label{Corollary Closure Order and Dimension}
    If $X_i \boldprec X_j$, then $\dim X_i < \dim X_j$.
\end{corollary}

This implies that it is possible to draw the Hasse diagram $\bm{\Gamma}\mathfrak{X}$ in the plane such that two vertices lie on the same horizontal line if and only if they have the same dimension, and all edges are directed upwards without touching any vertices other than their endpoints.

Let $\mathrm{Irr}(X)$ denote the set of irreducible components of $X$.
Since each $X_i \in \mathfrak{X}$ is irreducible, it is contained in some (not necessarily unique) element of $\mathrm{Irr}(X)$. 

\begin{proposition}\label{Proposition Unique Dense in Irreducible Component}
    Each irreducible component of $X$ contains a unique dense $X_i \in \mathfrak{X}$.
\begin{proof}
    Suppose $S \in \mathrm{Irr}(X)$, and note that $S = \overline{S} \cap X$.
    Moreover, $S = \bigcup_i \left(\:\! \overline{S \cap X_i} \cap X\right)$ is a finite union of closed subsets of $X$.
    It follows from the fact that $S$ is an irreducible component of $X$ that there exists $X_i \in \mathfrak{X}$ such that $S = \overline{X_i} \cap X$.
    Therefore, $X_i \subseteq S$ with $\overline{X_i} = \overline{S}$.
    For uniqueness, suppose $X_j \in \mathfrak{X}$ is also dense in $S$.
    Then $\overline{X_i} = \overline{X_j}$, and so \textsc{Lemma}~\ref{Lemma Closure Order is a Partial Order} implies that $X_i = X_j$.
\end{proof}
\end{proposition}

Given $S \in \mathrm{Irr}(X)$, we let $X_S$ denote the unique $X_i \in \mathfrak{X}$ which is dense in $S$.
Then $X_S$ is also the unique $X_i \in \mathfrak{X}$ such that $S = \overline{X_i} \cap X = \overline{S} \cap X$.

\begin{proposition}\label{Proposition Irreducible Components and Maximal Elements}
    Suppose $X_i \in \mathfrak{X}$.
\begin{itemize}
    \item[$\mathrm{(i)}$] $X_i$ is maximal in $\mathfrak{X}$ if and only if $X_i = X_S$ for some $S \in \mathrm{Irr}(X)$.
    \item[$\mathrm{(ii)}$] The irreducible components of $X$ are in bijection with the maximal elements in $\mathfrak{X}$.
    \item[$\mathrm{(iii)}$] If $X_i$ coincides with an irreducible component of $X$, then it is isolated in $\mathfrak{X}$.
\end{itemize}
\begin{proof}
    Suppose $S \in \mathrm{Irr}(X)$, and $X_j \in \mathfrak{X}$ is such that $X_S \boldpreceq X_j$.
    Let $S' \in \mathrm{Irr}(X)$ such that $X_j \subseteq S'$, and observe that $S = \overline{X_S} \cap X \subseteq \overline{X_j} \cap X \subseteq \overline{S^{'}} \cap X = S'$.
    Thus $S = S'$, from which the inclusions imply that $S = \overline{X_j} \cap X$.
    Therefore, $X_j = X_S$, and so $X_S$ is maximal in $\mathfrak{X}$.

    Conversely, suppose $X_i$ is maximal in $\mathfrak{X}$, and let $S \in \mathrm{Irr}(X)$ be such that $X_i \subseteq S$.
    Then $X_i \subseteq \overline{X_S} \cap X \subseteq \overline{X_S}$, and so the maximality of $X_i$ implies that $X_i = X_S$, which proves (i).

    The bijection for (ii) is given by the map $S \mapsto X_S$; it is well-defined and surjective by (i), and injective since $S = \overline{X_S} \cap X$ for each $S \in \mathrm{Irr}(X)$.

    If $X_i = S \in \mathrm{Irr}(X)$, then $X_i = X_S$, and so $X_i$ is maximal in $\mathfrak{X}$ by (i).
    Moreover, if $X_j \in \mathfrak{X}$ satisfies $X_j \boldpreceq X_i$, then $X_j \subseteq \overline{X_i} \cap X = X_i$, and so $X_j = X_i$.
    Therefore, $X_i$ is also minimal in $\mathfrak{X}$, which proves (iii).
\end{proof}
\end{proposition}

Regarding the converse to \textsc{Proposition}~\ref{Proposition Irreducible Components and Maximal Elements}(iii): if $X_i$ is isolated in $\mathfrak{X}$, then \textsc{Proposition}~\ref{Proposition Irreducible Components and Maximal Elements}(i) implies that $X_i = X_S$ for some $S \in \mathrm{Irr}(X)$.
Yet even though $X_i$ is minimal in $\mathfrak{X}$, it is possible to have $X_i \subsetneq S$ without further assumptions.

\begin{proposition}\label{Proposition Converse to Isolated Result}
    Assume that every irreducible component of $X$ is a union of some subcollection of $\mathfrak{X}$.
    Then $X_i \in \mathfrak{X}$ coincides with an irreducible component of $X$ if and only if it is isolated in $\mathfrak{X}$.
\begin{proof}
    Using \textsc{Proposition}~\ref{Proposition Irreducible Components and Maximal Elements}(iii), it suffices to prove that if $X_i$ is isolated in $\mathfrak{X}$, then it coincides with an irreducible component of $X$.
    So, suppose that $X_i$ is isolated in $\mathfrak{X}$, and (as above) let $S \in \mathrm{Irr}(X)$ be such that $X_i = X_S$.
    From our assumption, it follows that $S = \bigsqcup_{j \in I} X_j$, for some $I \subseteq \left\{ 1, \ldots, r \right\}$.
    For each $j \in I$, we have $X_j \subseteq S \subseteq \overline{X_i}$, and so the minimality of $X_i$ in $\mathfrak{X}$ implies that $X_j = X_i$.
    Therefore, $S = X_i$, as required.
\end{proof}
\end{proposition}

Finally, we provide a sufficient condition to prove that an element $X_i \in \mathfrak{X}$ is locally closed.

\begin{theorem}\label{Theorem Locally Closed}
    Suppose that $X_i \in \mathfrak{X}$ and, for each $X_j \boldpreceq X_i$, we have that $\overline{X_j}$ is a union of some subcollection of $\mathfrak{X}$.
    Then $X_i$ is locally closed.
\begin{proof}
    In particular, $X_i \boldpreceq X_i$, and so there exists (by assumption) $I \subseteq \left\{ 1, \ldots, r \right\}$ such that $\overline{X_i} = \bigsqcup_{j \in I} X_j$.
    Observe that $i \in I$, and define $J \coloneqq I \setminus \{i\}$.
    It follows that $\overline{X_i} = X_i \cup \bigcup_{j \in J} \overline{X_j}$.

    Suppose, for a contradiction, that $X_i \cap \overline{X_j} \neq \emptyset$, for some $j \in J$.
    By assumption, $\overline{X_j}$ is a union of some subcollection of $\mathfrak{X}$, which implies that $X_i \subseteq \overline{X_j}$.
    However, $X_j \subseteq \overline{X_i}$, so the antisymmetry of the closure order implies that $X_i = X_j$, which is a contradiction.
    Therefore, if we let $Z \coloneqq \bigcup_{j \in J} \overline{X_j}$, then $X_i = \overline{X_i} \cap ( Y \setminus Z )$ is locally closed.
\end{proof}
\end{theorem}

\section{Decomposition Classes}

We shall now define the Lie algebra decomposition classes of an arbitrary algebraic group $G$, over the algebraically closed field $\bbK$.

\subsection{Definition}

Two elements $x, y \in \g$ are \emph{Jordan equivalent}, written $x \sim y$, if there exists $g \in G$ such that $\Co{G}\left(g \cdot x\s\right) = \Co{G}y\s$ and $g \cdot x\n = y\n$.
Then $\sim$ is an equivalence relation on $\g$, and thus we may consider its equivalence classes.

\begin{definition}\label{Definition Decomposition Class}
    The $G$\itbf{-decomposition class} of $x \in \g$ is defined as its equivalence class with respect to $\sim$, and is denoted $\J{G} x = \left\{ y \in \g \midd x \sim y \right\}$. 
\end{definition}

We let $\fD[G]$ denote the set of $G$-decomposition classes, and note that this definition is slightly different from that of packets found in \cite[\S 1.2]{S82}.
We shall prove in \textsc{Corollary}~\ref{Corollary Equivalent Jordan Equivalence Definition} that (assuming $G$ is connected reductive) the two definitions coincide.
When the algebraic group $G$ is self-apparent, we shall simply refer to these as \emph{decomposition classes}.

Using the definition of Jordan equivalence, alongside \textsc{Lemma}~\ref{Lemma Stabiliser/Centraliser of Jordan Decomposition}, we can readily deduce the following result.

\begin{lemma}\label{Lemma Constant Dimension}
    Suppose that $x,y \in \g$ satisfy $x \sim y$.
\begin{itemize}
    \item[$\mathrm{(a)}$] $\dim \C_G x = \dim\C_G y$.
    \item[$\mathrm{(b)}$] $\dim \fc_{\g}x = \dim\fc_{\g}y$.
\end{itemize}
\end{lemma}

It follows from \textsc{Lemma}~\ref{Lemma Constant Dimension} that each decomposition class lies in both a unique stabiliser level set and a unique centraliser level set.
Given a level set $\ls$ of $\g$, if $\fD_{\am}[G] \coloneqq \left\{ \fJ \in \fD[G] \midd \fJ \subseteq \ls \right\}$ denotes the set of decomposition classes contained in $\ls$, then $\ls$ is the (disjoint) union of all $\fJ \in \fD_{\am}[G]$.

\subsection{Initial Properties}\label{Subsection Initial Properties}

From here onward, we always assume that $G$ is a connected \mbox{reductive} algebraic group.
For any semisimple $y \in \g$, let $\dreg{\g} y \coloneqq \left(\fd_{\g} y\right)^\reg$, and suppose that $T \subseteq G$ is a maximal torus.

For each $\alpha \in \Phi$, we consider the subgroup $G_{\alpha} \coloneqq \llangle \U_{\alpha}, \U_{-\alpha} \rrangle$.
Using the results in \cite[\S 8.3]{MT11}, there exists an isomorphism of algebraic groups $\varphi \colon H \rightarrow G_{\alpha}$, where $H$ is either $\mathrm{SL}_2$ or $\mathrm{PGL}_2$, such that the image of the standard maximal torus coincides with $T \cap G_{\alpha}$, and the differential maps the two root spaces of $\h$ to $\g_{\alpha}$ and $\g_{-\alpha}$.

\begin{proposition}\label{Proposition Double Centraliser Properties} 
    Suppose $y \in \ft$.
\begin{itemize}
    \item[$\mathrm{(i)}$] $\fd_{\g} y = \left\{ z \in \ft \midd \Phi_y \subseteq \Phi_z \right\}$.
    \item[$\mathrm{(ii)}$] $\dreg{\g} y = \left\{ z \in \ft \midd \Phi_y = \Phi_z \right\} = \left\{ z \in \g \midd \fc_{\g} y = \fc_{\g} z \right\} = \left\{ z \in \g \midd \Co{G}y = \Co{G}z \right\}$.
\end{itemize}
\begin{proof}
    Using \textsc{Corollary}~\ref{Corollary Semisimple Centraliser Description}, we have that $\fd_{\g} y \subseteq \fc_{\g} y = \ft \oplus \g(\Phi_y)$.
    Since $\fc_{\g}y$ is $T$-stable, so is $\fd_{\g}y = \fc_{\g}\left(\fc_{\g}y\right)$.
    Fix $\alpha \in \Phi_y$, and let $\varphi \colon H \rightarrow G_{\alpha}$ be an isomorphism of algebraic groups as described above, where $H$ is either $\mathrm{SL}_2$ or $\mathrm{PGL}_2$.
    If $H = \mathrm{SL}_2$, take $x = \md \varphi \begin{psmallmatrix} 0 & 1 \\ 0 & 0 \end{psmallmatrix} \in \g_{\alpha}$ and $x' = \md \varphi \begin{psmallmatrix} 0 & 0 \\ 1 & 0 \end{psmallmatrix} \in \g_{-\alpha}$.
    Otherwise, for $H = \mathrm{PGL}_2$, let $\pi \colon \GL_2 \rightarrow \mathrm{PGL}_2$ denote the canonical quotient homomorphism, and take $x = \md \varphi \big( \md \pi \begin{psmallmatrix} 0 & 1 \\ 0 & 0 \end{psmallmatrix} \big) \in \g_{\alpha}$ and $x' = \md \varphi \big( \md \pi \begin{psmallmatrix} 1 & 0 \\ 0 & 0 \end{psmallmatrix} \big) \in \ft$.
    Simple calculations in either case show that $x \in \g_{\alpha} \subseteq \fc_{\g} y$ and $x' \in \ft \oplus \g_{-\alpha} \subseteq \fc_{\g} y$ satisfy $\left[x,x'\right] \neq\nobreak 0$.
    It follows that $x \notin \fc_{\g} x'$, and thus $x \notin \fc_{\g}\left(\fc_{\g} y\right) = \fd_{\g} y$.
    Consequently, $\fd_{\g} y \cap \g_{\alpha}$ is a proper subspace of $\g_{\alpha}$.
    Since $\dim\g_{\alpha} = 1$, we have that $\fd_{\g} y \cap \g_{\alpha} = 0$, thus $\fd_{\g} y \subseteq \ft$.
    Therefore, \textsc{Corollary}~\ref{Corollary Semisimple Centraliser Description} implies that $\fd_{\g} y = \left\{ z \in \g \midd \fc_{\g}y \subseteq \fc_{\g}z \right\} = \big\{ z \in \ft \bigm| \g(\Phi_y) \subseteq \g(\Phi_z) \big\}$, from which (i) is clear.

    Since $\fd_{\g} y \subseteq \ft$ consists of semisimple elements, $\dreg{\g} y = \left(\fd_{\g} y\right)^\greg = \left\{ z \in \g \midd \fc_{\g} y = \fc_{\g} z \right\}$.
    The remaining equalities in (ii) are then a consequence of \cite[Lemma 3.7]{S75} and \textsc{Corollary}~\ref{Corollary Semisimple Centraliser Description}.
\end{proof}
\end{proposition}

We note that, unless $p = 2$ and $H = \mathrm{PGL}_2$, the element $x'$ in the proof of \textsc{Proposition}~\ref{Proposition Double Centraliser Properties}(i) can be chosen to lie in $\g_{-\alpha}$; explicitly, if we let $x' = \md \varphi \big( \md \pi \begin{psmallmatrix} 0 & 0 \\ 1 & 0 \end{psmallmatrix} \big) \in \g_{-\alpha}$, then $\left[x,x'\right] = 0$ if and only if $p = 2$.

Using \textsc{Proposition}~\ref{Proposition Double Centraliser Properties}(ii), we can prove the following corollary, which demonstrates that (for connected reductive algebraic groups) our definition of decomposition classes given in \textsc{Definition}~\ref{Definition Decomposition Class} coincides with that of packets given in \cite[\S 1.2]{S82}.

\begin{corollary}\label{Corollary Equivalent Jordan Equivalence Definition}
    Suppose $x, y \in \g$.
    Then $x \sim y$ if and only if there exists $g \in G$ such that $\fc_{\g}\left( g \cdot x\s\right) = \fc_{\g}y\s$ and $g \cdot x\n = y\n$.  
\end{corollary}

We then have the following initial properties of decomposition classes---some of which are found in \cite[\S 3.3]{B98}, but not in the level of generality presented here.

\begin{theorem}\label{Theorem Initial Decomposition Class Properties}
    Suppose $x \in \g$.
\begin{itemize}
    \item[$\mathrm{(i)}$] $\J{G} x = G \cdot \left( \dreg{\g} x\s + x\n \right)$.  
    \item[$\mathrm{(ii)}$] $\J{G} x$ is $G$-stable, and $\bbKx$-stable. 
    \item[$\mathrm{(iii)}$] $\J{G} x$ is irreducible and constructible.  
    \item[$\mathrm{(iv)}$] $\J{G} x + \z(\g) = \J{G} x$.  
    \item[$\mathrm{(v)}$] $\J{G}x\n = \z(\g) + G \cdot x\n$.  
    \item[$\mathrm{(vi)}$] There are only finitely many $G$-decomposition classes.
\end{itemize}
\begin{proof}
    If $y \in \J{G} x$, then there exists $g \in G$ such that $\fc_{\g}\left( g \cdot x\s\right) = \fc_{\g}y\s$ and $g \cdot x\n = y\n$.
    Using \textsc{Proposition}~\ref{Proposition Double Centraliser Properties}(ii), $g^{-1} \cdot y\s \in \dreg{\g} x\s$, and thus $y = g \cdot \left( g^{-1} \cdot y\s + x\n\right) \in G \cdot \left( \dreg{\g} x\s + x\n \right)$.
    Conversely, if $y \in G \cdot \left( \dreg{\g} x\s + x\n \right)$, then there exists $h \in G$ and $z \in \dreg{\g} x\s$ such that $y = g \cdot \left(z + x\n\right)$.
    Since $y\s = g \cdot z$, it follows from \textsc{Proposition}~\ref{Proposition Double Centraliser Properties}(ii) that $\fc_{\g} \left( g^{-1} \cdot y\s\right) = \fc_{\g} x\s$.
    Therefore, $x \sim y$, and so $y \in \J{G} x$, as required for (i).

    The first part of (ii) is immediate from (i), so suppose $\lambda \in \mathbb{K}^{\times}$.
    Using \textsc{Corollary}~\ref{Corollary Semisimple Centraliser Description} and \cite[Lemma 2.10]{J04}, we know that there exists $g \in \Co{G}x\s$ such that $g \cdot x\n = \lambda x\n$.
    Then $\fc_{\g} \left( g \cdot x\s\right) = \fc_{\g}x\s = \fc_{\g} \left(\lambda x\s\right)$ and $g \cdot x\n = \lambda x\n$, hence $\lambda x \in \J{G}x$, which proves the second part of (ii).

    Since $\dreg{\g} x\s$ is irreducible and locally closed, so is $\dreg{\g} x\s + x\n$.
    Since $G$ is connected, and $\J{G}x$ is the image of $G \times \left( \dreg{\g} x\s + x\n \right)$ under the adjoint action, we know $\J{G}x$ is also irreducible.
    Moreover, \cite[Corollary AG10.2]{B91} shows that $\J{G}x$ is constructible, completing the proof of (iii).

    Suppose $z \in \dreg{\g} x\s$ and $g \in G$.
    If $y \in \z(\g)$, then the Jordan decomposition of $g \cdot \left( z + x\n\right) + y$ is $g \cdot \left(z + y\right) + g \cdot x\n$, and $\fc_{\g}\left(z+y\right) = \fc_{\g} z$.
    Therefore, (i) implies that $x \sim g \cdot \left( z + x\n\right) \sim g \cdot \left(z + x\n\right) + y$, from which (iv) is immediate.
    Observing that $\dreg{\g} 0 = \z(\g)^\reg = \z(\g)$ is $G$-stable shows that (v) follows by applying (i) to $x\n$.

    Suppose $T \subseteq G$ is a maximal torus, with root system $\Phi$.
    Using \textsc{Corollary}~\ref{Corollary Semisimple Centraliser Description}, each $y \in \ft$ determines a subset of roots $\Phi_y \subseteq \Phi$ from which $\fc_{\g} y$ is determined.
    Since each semisimple element of $\g$ is $G$-conjugate to an element of $\ft$, there are (up to $G$-conjugacy) only finitely many centralisers of semisimple elements.
    For each semisimple element $z \in \g$, \textsc{Corollary}~\ref{Corollary Semisimple Centraliser Description} shows there are only finitely many nilpotent $\Co{G} z$-orbits in $\fc_{\g} z$, hence (vi) follows.
\end{proof}
\end{theorem}

A ($G$-)\emph{decomposition datum} corresponding to a decomposition class $\fJ \in \fD[G]$ is defined as a pair $\left(\Co{G}x\s;x\n\right)$, for some $x \in \fJ$.
It is clear that decomposition data are unique up to the action of $G$, where $G$ acts via conjugation on the first argument and via the adjoint action on the second argument.

For any $M \in \mathfrak{C}_G \coloneqq \left\{ \Co{G}x\s \midd x \in \g \right\}$ and $e_0 \in \cN_{\fm}$, we let $\J{G}\left(M;e_0\right) \coloneqq \J{G}\left(y + e_0\right)$ denote the corresponding decomposition class, where $y = y\s \in \g$ is such that $M = \Co{G}y$.
A \emph{saturated} ($G$-)\emph{decomposition datum} is any pair $\left(M;\mathcal{O}\right)$ consisting of some $M \in \mathfrak{C}_G$ and a nilpotent $M$-orbit $\mathcal{O} \in \smallslant{\cN_{\fm}}{M}$.
We then define $\J{G}\left(M;\mathcal{O}\right) \coloneqq \J{G}(M;e_0)$, for any $e_0 \in \mathcal{O}$, and note that $\J{G}\left(M;\mathcal{O}\right) = G \cdot \big(\z(\fm)^\reg + \mathcal{O}\big)$ by \textsc{Theorem}~\ref{Theorem Initial Decomposition Class Properties}(i).

\subsection{The Dimension of a Decomposition Class}

We proceed to prove \textsc{Theorem}~\ref{Maintheorem Dimension of Decomposition Class}, previously established under the Standard Hypotheses in \cite[Proposition 2.5]{PS18}.

\begin{theorem}\label{Theorem Dimension of Decomposition Class}
    Let $\fJ \in \fD[G]$.
    Then $\dim \fJ = \dim G \cdot x + \dim \fd_{\g}x\s$, for any $x \in \fJ$.  
\begin{proof}
    Fix  $x \in \fJ$, and let $\z \coloneqq \fd_{\g}x\s$.
    Let $T \subseteq G$ be a maximal torus such that $x\s \in \ft$, with corresponding Weyl group $W = \mathrm{W}(G,T)$, and recall from \textsc{Proposition}~\ref{Proposition Double Centraliser Properties}(i) that $\z \subseteq \ft$.
    Let $\alpha \colon G \times \left(\z^\reg + x\n\right) \rightarrow \fJ$ denote the adjoint action morphism, which is a surjective morphism of irreducible varieties by \textsc{Theorem}~\ref{Theorem Initial Decomposition Class Properties}.
    
    Let $u \in \fJ$, and fix $\left(g_0,z_0 + x\n\right) \in \alpha^{-1}(u)$.
    If we let $G$ act on $G \times (\mathfrak{z}^{\mathrm{reg}} + x_{\mathrm{n}})$ via left multiplication on $G$ (and the identity on $\mathfrak{z}^{\mathrm{reg}} + x_{\mathrm{n}}$), then $\alpha$ is $G$-equivariant, from which it follows that $\alpha^{-1}(u) = \big\{ \left(g_0 g, z + x\n\right) \in G \times \left(\z^\reg + x\n\right) \bigm| \left(g, z + x\n\right) \in \alpha^{-1}\left(z_0 + x\n\right) \big\}$, and thus $\dim \alpha^{-1}(u) = \dim \alpha^{-1}\left(z_0 + x\n\right)$.

    Using \cite[\S 7.12(1)]{J04}, we know that $G \cdot z_0 \cap \z^\reg \subseteq G \cdot z_0 \cap \ft = W \cdot z_0$ is finite.
    Let $Z \coloneqq \left\{ z \in G \cdot z_0 \cap \mathfrak{z}^{\mathrm{reg}} \midd \exists g \in \C_G x_{\mathrm{n}} \colon g \cdot z = z_0 \right\}$, which is non-empty since $z_0 \in Z$.
    For each $z \in Z$, fix $n_z \in \C_G x\n$ such that $n_z \cdot z = z_0$.
    It follows that $\alpha^{-1}\left(z_0 + x\n\right)$ is the finite disjoint union of the subsets $A_z \coloneqq \big\{(gn_z, z + x\n) \bigm| g \in G,\, g \cdot z = z,\, g \cdot x\n = x\n \big\}$, over $z \in Z$.
    Since $g \cdot z = z$ and $g \cdot x\n = x\n$ if and only if $g \in \C_G(z+x\n)$, we have that $\dim A_z = \dim \C_G(z + \nobreak x\n) = \dim \C_G(z_0 + x\n)$.
    Hence, $\dim \alpha^{-1}(z_0+x\n) = \dim \C_G(z_0+x\n) = \dim \C_G x$.

    Since we have shown that all fibres of $\alpha$ have dimension $\dim \C_G x$, \cite[Theorem AG10.1]{B91} implies that $\dim G + \dim \fd_{\g}x\s = \dim \fJ + \dim \C_G x$.
    Now the result follows from $\dim G \cdot x = \dim G - \dim \C_G x$.
\end{proof}
\end{theorem}

We can immediately deduce the following alternative formulations of \textsc{Theorem}~\ref{Theorem Dimension of Decomposition Class}.

\begin{corollary}\label{Corollary Alternative Dimension Formula}
    Suppose $x \in \g$, and let $M \coloneqq \Co{G}x\s \in \mathfrak{C}_G$.
    Then
    \begin{IEEEeqnarray*}{RCL}
        \dim \J{G}x & = & \dim G + \dim \z(\fm) - \dim \C_G x \\
        & = & \dim G + \dim \z(\fm) + \dim M \cdot x\n - \dim M.
    \end{IEEEeqnarray*}
\end{corollary}

\subsection{The Closure Order on Decomposition Classes}\label{Subsection The Closure Order on Decomposition Classes}

We refer to the closure of a \mbox{($G$-)decomposition class} as a ($G$-)\emph{decomposition variety}.
It follows from \textsc{Theorem}~\ref{Theorem Initial Decomposition Class Properties} that each decomposition variety is $G$-stable, irreducible, and $\bbKx$-stable.
In fact, they are stable under arbitrary scalar multiplication: if $z \in \dvJ$, then $\bbKx z \subseteq \dvJ$, and thus $\bbK z = \overline{\bbKx z} \subseteq \dvJ$.
Since $0 \in \dvJ$, we have that $\dvJ \cap \Ng \neq \emptyset$.
Moreover, it follows from $\J{G}0 = \z(\g)$ that $\fJ = \dvJ$ if and only if $\fJ = \z(\g)$, which is stated in \cite[\S 3.1]{A25}.

Using \textsc{Theorem}~\ref{Theorem Initial Decomposition Class Properties}, we know that $\fD[G]$ is a finite collection of irreducible constructible subsets of $\g$ such that $\g$ is their disjoint union.
Therefore, we can make use of the results in \S\ref{Subsection Closure Orders} regarding the \emph{closure order} $\boldpreceq$ on $\fD[G]$, defined so that $\fJ' \boldpreceq \fJ$ if and only if $\fJ' \subseteq \overline{\fJ}$.
In particular, using \textsc{Lemma}~\ref{Lemma Closure Order is a Partial Order} and \textsc{Corollary}~\ref{Corollary Closure Order and Dimension}, we can deduce the following.

\begin{corollary}\label{Corollary Closure Order for Decomposition Classes}
    The closure order $\boldpreceq$ is a well-defined partial order on $\fD[G]$.
    Moreover, if $\fJ' \boldprec \fJ$, then $\dim \fJ' < \dim \fJ$.
\end{corollary}

The following result shows that the stabiliser and centraliser dimensions weakly increase as we descend in the closure order.

\begin{lemma}\label{Lemma Stabiliser/Centraliser Dimension and Closure Order}
    Suppose $x,y \in \g$ with $\J{G}y \boldpreceq \J{G}x$.
\begin{itemize}
    \item[$\mathrm{(a)}$] $\dim \C_G y \geq \dim \C_G x$.
    \item[$\mathrm{(b)}$] $\dim \fc_{\g} y \geq \dim \fc_{\g} x$.
\end{itemize}
\begin{proof}
    For (a), let $m \coloneqq \dim \C_G x$, and recall from \textsc{Lemma}~\ref{Lemma Constant Dimension}(a) that $\J{G}x \subseteq \mathfrak{g}_{(m)}$.
    Using \textsc{Lemma}~\ref{Lemma Closure inside a Level Set}(a), we know $\overline{\J{G}x} \subseteq \bigsqcup_{\,n \geq m} \g_{(n)}$.
    Then $\J{G}y \subseteq \bigsqcup_{\,n \geq m} \g_{(n)}$ implies that $\dim \C_G y \geq \dim \C_G x$.
    The proof of (b) follows by an analogous argument.
\end{proof}
\end{lemma}

We also have the following result, which provides an analogue of \textsc{Theorem}~\ref{Theorem Initial Decomposition Class Properties}(iv) and (v) for decomposition varieties, and can be readily deduced from the corresponding parts of \textsc{Theorem}~\ref{Theorem Initial Decomposition Class Properties}, alongside \textsc{Lemma}~\ref{Lemma Closure of Saturated Set} and \textsc{Proposition}~\ref{Proposition Closure of Semisimple + Nilpotent}.

\begin{corollary}\label{Corollary Decomposition Varieties Properties}
    Suppose $x \in \g$.
\begin{itemize}
    \item[$\mathrm{(i)}$] $\dv{G}{x} + \z(\g) = \dv{G}{x}$. 
    \item[$\mathrm{(ii)}$] $\dv{G}{x\n} = \z(\g) + \overline{G \cdot x\n}$. 
\end{itemize}
\end{corollary}

Recall that \textsc{Theorem}~\ref{Theorem Initial Decomposition Class Properties}(iii) established that every decomposition class is constructible and thus a finite union of locally closed sets.
Therefore, the following consequence of \textsc{Theorem}~\ref{Theorem Locally Closed} strengthens this result for certain decomposition classes.

\begin{corollary}\label{Corollary Locally Closed Decomposition Classes}
    Suppose that $\fJ \in \fD[G]$ and, for each $x \in \overline{\fJ}$, we have that $\overline{\J{G}x}$ is a union of decomposition classes.
    Then $\fJ$ is locally closed.
\end{corollary}

We shall see in \textsc{Proposition}~\ref{Proposition Levi-type Decomposition Variety as Union of Decomposition Classes} that $\overline{\J{G}x}$ is a union of decomposition classes whenever $\J{G}x$ is a Levi-type decomposition class---see \textsc{Definition}~\ref{Definition Levi-Type}.
In particular, this always holds when $p \geq 0$ is good for $G$.
On the other hand, \textsc{Appendix}~\ref{Appendix: Example} demonstrates that $\overline{\J{\mathrm{Sp}^{\vphantom{\circ}}_4}x}$ is a union of decomposition classes for any $x \in \mathfrak{sp}_4$, even in bad characteristic.

\begin{remark}
    There is an analogue of Lie algebra decomposition classes inside the connected reductive algebraic group $G$ itself, which are referred to as \emph{Jordan classes}.
    These were originally introduced by Lusztig in \cite[\S3.1]{L84}, and have also seen renewed attention in recent years.
    Following \cite[\S2.1]{C20}, if an element $g \in G$ has multiplicative Jordan decomposition $g=su$, then its Jordan class is given by $J(su) \coloneqq G \cdot \left( \bigl(\mathcal{Z}^{\circ}(\Co{G}s)s\bigr)^{\mathrm{reg}}u\right)$.
    For the benefit of the reader, we have highlighted some of the differences between the theory of Jordan classes and that of Lie algebra decomposition classes:

        \begin{itemize}
        \item As observed in \cite[\S3]{ACE21}, not all closures of Jordan classes contain a unipotent conjugacy class, even in characteristic $0$. In contrast, we have already seen that the closure of every decomposition class contains $0$.
        
        \item The stabiliser of any subtorus is a Levi subgroup, hence $\C_G\bigl(\mathcal{Z}^{\circ}(\Co{G}s)\bigr)$ is always a Levi subgroup; as shown in \cite[\S3.1]{L84}, it is actually the smallest Levi subgroup of $G$ containing $\Co{G}s$. This allows for certain properties of Levi subgroups to be utilised in the study of Jordan classes, even though $\Co{G}s$ is not itself a Levi subgroup.
        There is no way to emulate this construction for Lie algebra decomposition classes---\textsc{Proposition}~\ref{Proposition Double Centraliser Properties} implies that $\fc_{\mathfrak{g}}\bigl(\mathfrak{z}(\fc_{\mathfrak{g}}y)\bigr) = \fc_{\mathfrak{g}}y$, for any semisimple $y \in \mathfrak{g}$.

        \item As a consequence, uniform descriptions of the closures of Jordan classes have been determined.
        On the other hand, these are only known for certain decomposition classes---see \textsc{Theorem}~\ref{Theorem Decomposition Variety as Union of Orbit Closures}.

        \item Jordan classes are known to be locally closed in all cases, as proved in \cite[\S3.1]{L84}. However, the same is not yet known for decomposition classes---see \textsc{Corollary}~\ref{Corollary Strongly-Levi-Type Decomposition Classes are Locally Closed} and \textsc{Corollary}~\ref{Corollary Good Characteristic Decomposition Classes are Locally Closed} for results regarding this.
        The principal obstruction in the Lie algebra case is the issue regarding Levi subgroups raised above.
    \end{itemize}
\end{remark}

\section{Preservation of Decomposition Classes}\label{Section Preservation of Decomposition Classes}

In this section, we shall explore how decomposition classes interact with direct products, central surjections, and separable central surjections.

\subsection{Preservation by Direct Products}

Suppose that $H$ is a connected reductive algebraic group, just as $G$ is.
The direct product $G \times H$ is likewise a connected reductive algebraic group, with Lie algebra $\g \oplus \h$.
If $x \in \g$ and $y \in \h$, it follows readily from the definitions that $\J{G \times H}(x + y) = \J{G}x + \J{H}y$.
Consequently, $\fD[G \times H] = \fD[G] \times \fD[H]$ as sets, where the closure order on $\fD[G \times H]$ coincides with the product order induced from $\fD[G]$ and $\fD[H]$.
By induction, this extends to arbitrary finite direct products of connected reductive algebraic groups.

\subsection{Preservation by Central Surjections}\label{Subsection Preservation by Central Surjections}

A surjective homomorphism of algebraic groups $\varphi \colon G \rightarrow H$ is said to be \emph{central} if $\ker \varphi \subseteq \mathcal{Z}(G)$ and $\ker \md\varphi \subseteq \z(\g)$.
Suppose that $\check{G}$ is also a connected reductive algebraic group, such that there exists a central surjection $\varphi \colon G \rightarrow \check{G}$.
Given any $x \in \g$, let $\check{x} \coloneqq \md \varphi (x) \in \check{\g}$ denote its image under the differential, and note that $\md\varphi \colon \g \rightarrow \check{\g}$ is not necessarily surjective (see \S\ref{Subsection Preservation by Separable Central Surjections}).

In \cite[Proposition 2.7(a)]{J04}, Jantzen proves that $\md\varphi \colon \g \rightarrow \check{\g}$ restricts to a bijection $\cN_{\g} \rightarrow \cN_{\check{\g}}$, which in turn induces a bijection $\smallslant{\cN_{\g}}{G} \rightarrow \smallslant{\cN_{\check{\g}}}{\check{G}}$.
Moreover, for each $x \in \cN_{\g}$, the differetial $\md\varphi$ further restricts to a bijection $G \cdot x \rightarrow \check{G} \cdot \check{x}$, with $\dim G \cdot x = \dim \check{G} \cdot \check{x}$ and $\varphi(\C_G x) = \C_{\check{G}} \check{x}$.
We expand upon this with the following result.

\begin{lemma}\label{Lemma Central Surjections preserve Nilpotent Closures}
    If $Y \subseteq \cN_{\g}$ is $G$-stable, then $\overline{\md\varphi(Y)} = \md\varphi\big(\:\!\overline{Y}\:\!\big)$. 
\begin{proof}
    If $X \coloneqq Y + \ker \md\varphi$, then $\md\varphi(X) = \md\varphi(Y)$, and so \textsc{Lemma}~\ref{Lemma Closure of Saturated Set} implies that $\md\varphi\big(\:\!\overline{X}\:\!\big) = \overline{\md\varphi(X)} = \overline{\md\varphi(Y)}$.
    Moreover, \textsc{Proposition}~\ref{Proposition Closure of Semisimple + Nilpotent} implies that $\overline{X} = \overline{Y} + \ker \md\varphi$, and thus $\md\varphi\big(\:\!\overline{X}\:\!\big) = \md\varphi\big(\:\!\overline{Y}\:\!\big)$, from which the result follows.
\end{proof}
\end{lemma}

Suppose $T \subseteq G$ is a maximal torus, and $B \in \mathfrak{P}(G,T)$.
The induced comorphism $\varphi^* \colon \bbK[\check{G}] \rightarrow \bbK[G]$ restricts to a homomorphism of character groups $\left(\varphi\restr_T\right)^* \colon \mathrm{X}(\check{T}) \rightarrow \mathrm{X}(T)$.
Then \cite[Proposition 22.4]{B91} implies that this further restricts to a bijection of root systems $\check{\Phi} =\Phi(\check{G},\check{T}) \rightarrow \Phi = \Phi(G,T)$.
Moreover, $\check{\Phi} \rightarrow \Phi$ again restricts to a bijection $\check{\Phi}_{\check{z}} \rightarrow \Phi_z$, for any $z \in \ft$.

\pagebreak

\begin{lemma}\label{Lemma Central Surjections and Semisimple Centralisers}
    If $y \in \g$ is semisimple, then $\md\varphi\left(\dreg{\g}y\right) = \md\varphi(\g) \cap \dreg{\check{\g}} \check{y}$.  
\begin{proof}
    Fix a maximal torus $T \subseteq G$ such that $y \in \ft$, and let $\Phi = \Phi(G,T)$.
    If $x \in \dreg{\g}y$, then \textsc{Proposition}~\ref{Proposition Double Centraliser Properties}(ii) implies that $x \in \ft$ and $\Phi_y = \Phi_x$, and thus $\check{\Phi}_{\check{y}} = \check{\Phi}_{\check{x}}$.
    It follows that $\md\varphi\left(\dreg{\g}y\right) \subseteq \md\varphi(\g) \cap \dreg{\check{\g}} \check{y}$.
    
    Conversely, suppose $x \in \g$ is such that $\check{x} \in \dreg{\check{\g}} \check{y} \subseteq \Lie \check{T}$.
    Since $\ker \md\varphi \subseteq \z(\g) \subseteq \ft$ and $\md\varphi(\ft) = \md\varphi(\g) \cap \Lie \check{T}$, it follows that $x \in \ft$.
    Using \textsc{Proposition}~\ref{Proposition Double Centraliser Properties}(ii) once more, we have that $\check{\Phi}_{\check{y}} = \check{\Phi}_{\check{x}}$, and so $\Phi_y = \Phi_x$.
    Therefore, $x \in \dreg{\g}y$, from which it follows that $\md\varphi(\g) \cap \dreg{\check{\g}} \check{y} \subseteq \md\varphi\left(\dreg{\g}y\right)$.
\end{proof}
\end{lemma}

Observe that $\check{x}\s = \md\varphi\left(x\s\right)$ and $\check{x}\n = \md\varphi\left(x\n\right)$, for any $x \in \g$.

\begin{theorem}\label{Theorem Central Surjections and Decomposition Classes}
    Suppose $x,y \in \g$.
\begin{itemize}
    \item[$\mathrm{(i)}$] $\md\varphi\left(\J{G}x\right) = \md\varphi(\g) \cap \J{\check{G}}\check{x}$. 
    \item[$\mathrm{(ii)}$] $\md\varphi\left(\:\!\dv{G}{x}\:\!\right) = \overline{\md\varphi(\mathfrak{\g}) \cap \J{\check{G}} \check{x}}$. 
    \item[$\mathrm{(iii)}$] $\J{G}x \subseteq \dv{G}{y}$ if and only if $\md\varphi\left(\J{G}x\right) \subseteq \md\varphi\left(\:\!\dv{G}{y}\:\!\right)$.  
\end{itemize}
\begin{proof}
    Using \textsc{Theorem}~\ref{Theorem Initial Decomposition Class Properties}(i) in conjunction with \textsc{Lemma}~\ref{Lemma Central Surjections and Semisimple Centralisers}, we have that $\md\varphi\left(\J{G}x\right) = \varphi(G) \cdot \big(\md\varphi\left(\dreg{\g}x\s\right) + \md\varphi\left(x\n\right)\big) = \check{G} \cdot \left((\md\varphi(\g) \cap \dreg{\check{\g}}\check{x}\s) + \check{x}\n\right)$.
    Since $\md\varphi(\g)$ is $\check{G}$\mbox{-stable}, we have that $\check{G} \cdot \left((\md\varphi(\g) \cap \dreg{\check{\g}}\check{x}\s) + \check{x}\n\right) = \md\varphi(\g) \cap \check{G} \cdot \left(\dreg{\check{\g}}\check{x}\s + \check{x}\n\right)$, and thus (i) follows by using \textsc{Theorem}~\ref{Theorem Initial Decomposition Class Properties}(i) once more.

    Since $\ker\md\varphi \subseteq \z(\g)$, it follows from \textsc{Theorem}~\ref{Theorem Initial Decomposition Class Properties}(iv) and \textsc{Lemma}~\ref{Lemma Closure of Saturated Set} that $\md\varphi\left(\:\!\dv{G}{x}\:\!\right) = \overline{\md\varphi\left(\J{G}x\right)}$.
    Therefore, (ii) follows from (i).
    
    The forward direction of (iii) is immediate, so suppose that $\md\varphi\left(\J{G}x\right) \subseteq \md\varphi\left(\:\!\dv{G}{y}\:\!\right)$, and let $z \in \J{G}x$.
    Then $\md\varphi(z) \in \md\varphi\left(\:\!\dv{G}{y}\:\!\right)$, and so $z \in (\md\varphi)^{-1}\pigl(\md\varphi\left(\:\!\dv{G}{y}\:\!\right)\pigr) = \dv{G}{y} + \ker\md\varphi \subseteq \dv{G}{y} +\nobreak \z(\g)$.
    It then follows from \textsc{Corollary}~\ref{Corollary Decomposition Varieties Properties}(i) that $z \in \dv{G}{y}$, which proves the other direction of (iii).
\end{proof}
\end{theorem}

One might expect that \textsc{Theorem}~\ref{Theorem Central Surjections and Decomposition Classes}(i) could be strengthened to say that $\mathrm{d}\varphi(\J{G}x) = \J{\check{G}}\check{x}$.
However, as the following example demonstrates, this is not always the case.

\begin{example}
    Let $G \coloneqq \mathrm{SL}_3$ and $\check{G} \coloneqq \mathrm{PGL}_3$, both defined over an algebraically closed field $\mathbb{K}$ with characteristic $p=3$.
    Let $\pi \colon \mathrm{GL}_3 \rightarrow \mathrm{PGL}_3$ be the canonical projection, and observe that its restriction to $\mathrm{SL}_3$ is a central surjection $\varphi \colon G \rightarrow \check{G}$.
    If we identify the prime subfield $\mathbb{F}_3 = \{0,1,2\}$ with its image in $\mathbb{K}$, then $x = x_{\mathrm{s}} \coloneqq \begin{psmallmatrix}
        0 & 0 & 0 \\
        0 & 1 & 0 \\
        0 & 0 & 2
    \end{psmallmatrix} \in \mathfrak{g}$.
    We claim that $\J{\check{G}}\check{x} \not\subseteq \mathrm{d}\varphi(\mathfrak{g})$.

    Since $x \in \mathfrak{g}$ is semisimple, we have that $\J{\check{G}}\check{x} = \check{G} \cdot \fd_{\mathfrak{g}}^{\mathrm{reg}} \check{x}$.
    Using \textsc{Proposition}~\ref{Proposition Double Centraliser Properties}(ii), one can show that $\fd_{\mathfrak{g}}^{\mathrm{reg}} \check{x} = \left\{ \mathrm{d}\pi \scalemath{0.8}{\begin{psmallmatrix}
        a_1 & 0 & 0 \\
        0 & a_2 & 0 \\
        0 & 0 & a_3
    \end{psmallmatrix}} \midd (a_1,a_2,a_3) \in \mathbb{K}^3_{\neq} \right\}$.
    On the other hand, $\mathrm{d}\varphi(\mathfrak{g}) = \left\{ \mathrm{d}\pi(y) \midd y \in \mathfrak{gl}_3, \mathrm{tr}\,y=0\right\}$.
    The field $\mathbb{K}$ is an infinite dimensional vector space over $\mathbb{F}_3$, hence there exists $b_1,b_2,b_3 \in \mathbb{K}$ which are linearly independent over $\mathbb{F}_3$.
    In particular, $(b_1,b_2,b_3) \in \mathbb{K}^3_{\neq}$ and $b_1+b_2+b_3 \neq 0$.
    It follows that $\tilde{y} \coloneqq \mathrm{d}\pi \scalemath{0.8}{\begin{psmallmatrix}
        b_1 & 0 & 0 \\
        0 & b_2 & 0 \\
        0 & 0 & b_3
    \end{psmallmatrix}} \in \fd_{\mathfrak{g}}^{\mathrm{reg}} \check{x} \subseteq \J{\check{G}}\check{x}$, but $\tilde{y} \notin \mathrm{d}\varphi(\J{G}x)$.
\end{example}

\subsection{Preservation by Separable Central Surjections}\label{Subsection Preservation by Separable Central Surjections}

Suppose still that $\varphi \colon G \rightarrow \check{G}$ is a central surjection (of connected reductive algebraic groups), and retain the other notation from \S\ref{Subsection Preservation by Central Surjections}.
It then follows from \cite[Theorem 4.3.7(iii)]{S98} that $\varphi$ is separable if and only if $\md\varphi \colon \g \rightarrow \check{\g}$ is surjective.

\begin{theorem}\label{Theorem Separable Central Surjections and Decomposition Classes}
    Suppose that $\varphi \colon G \rightarrow \check{G}$ is a separable central surjection, and let $x,y \in \g$.
    Then $\J{G}x \mapsto \J{\check{G}}\check{x}$ is a bijection $\fD[G] \rightarrow \fD[\check{G}]$ with the following properties$:$
\begin{itemize}
    \item[$\mathrm{(i)}$] $\md\varphi \colon \g \rightarrow \check{\g}$ restricts to a surjection $\J{G}x \rightarrow \J{\check{G}}\check{x}$. 
    \item[$\mathrm{(ii)}$] $\md\varphi\left(\:\!\dv{G}{x}\:\!\right) = \overline{\J{\check{G}} \check{x}}$.  
    \item[$\mathrm{(iii)}$] $\J{G}x \boldpreceq \J{G}y$ if and only if $\J{\check{G}}\check{x} \boldpreceq \J{\check{G}}\check{y}$.  
    \item[$\mathrm{(iv)}$] $\dim \J{G}x = \dim\ker \md\varphi +  \dim \J{\check{G}}\check{x}$.  
\end{itemize}
\begin{proof}
    Using \textsc{Theorem}~\ref{Theorem Central Surjections and Decomposition Classes}(i), we know that $\md\varphi\left(\J{G}x\right) = \J{\check{G}}\check{x} \in \fD[\check{G}]$, from which we can conclude that $\J{G}x \mapsto \J{\check{G}}\check{x}$ is a well-defined map $\fD[G] \rightarrow \fD[\check{G}]$.
    Its surjectivity follows immediately from that of $\md\varphi \colon \g \rightarrow \check{\g}$.
    Since \textsc{Theorem}~\ref{Theorem Initial Decomposition Class Properties}(iv) implies that $\J{G}x = \J{G}x + \ker\md\varphi = (\md\varphi)^{-1}\left(\J{\check{G}}\check{x}\right)$, its injectivity also follows.

    Since $\md\varphi \colon \g \rightarrow \check{\g}$ is a surjection, so is its restriction to $\J{G}x$, which proves (i).
    Moreover, (ii) and (iii) immediately follow from \textsc{Theorem}~\ref{Theorem Central Surjections and Decomposition Classes}(ii) and (iii), respectively.
    
    Using (i) and (ii), we know that $\md\varphi \colon \g \rightarrow \check{\g}$ restricts to a surjective morphism $\dv{G}{x} \rightarrow \overline{\J{\check{G}} \check{x}}$ of irreducible varieties.
    Then (iv) follows from \cite[Theorem 5.1.6]{S98}.
\end{proof}
\end{theorem}

It follows that, whenever there is a separable central surjection $\varphi \colon G \rightarrow \check{G}$, the Hasse diagrams of $\fD[G]$ and $\fD[\check{G}]$ are the same, up to a relabelling of the vertices.
That is, there is an induced isomorphism of directed graphs $\bm{\Gamma}\fD[G] \rightarrow \bm{\Gamma}\fD[\check{G}]$.

Let $G_{\ad}$ denote the adjoint group of $G$, with corresponding central surjection $\pi \colon G \rightarrow G_{\ad}$.
Using \cite[Corollary 2.3.7]{L05}, we have that $\pi \colon G \rightarrow G_{\ad}$ is separable if and only if $p$ does not divide $\big\lvert \left(\smallslant{\mathrm{X}(T)}{\mathbb{Z}\Phi}\right)_{\mathrm{tor}}\big\rvert$.
Therefore, if $p \geq 0$ is very good for $G$, then we can apply \textsc{Theorem}~\ref{Theorem Separable Central Surjections and Decomposition Classes} to the separable central surjection $\pi \colon G \rightarrow G_{\ad}$.

On the other hand, if $G = \GL_n$, then $G_{\ad} = \mathrm{PGL}_n$ and $\left(\smallslant{\mathrm{X}(T)}{\mathbb{Z}\Phi}\right)_{\mathrm{tor}}$ is trivial.
Therefore, (for any characteristic) the canonical projection $\pi \colon \GL_n \rightarrow \mathrm{PGL}_n$ is a separable central surjection, and thus we can apply \textsc{Theorem}~\ref{Theorem Separable Central Surjections and Decomposition Classes} to conclude that $\bm{\Gamma}\fD\left[\mathrm{PGL}_n\right]$ and $\bm{\Gamma}\fD\left[\GL_n\right]$ are isomorphic as directed graphs.

\section{Sheets}\label{Section Sheets}

As explained in the introduction, decomposition classes were originally introduced in \cite{BK79} as a tool to study the sheets of $\g$.
In the existing literature, these are the maximal irreducible subsets of $\g$ consisting of equal-dimension orbits.
However, we will make a departure with the following definitions.

\begin{definition}\label{Definition Sheets}
~ 
    \begin{itemize}
        \item An irreducible component $S$ of a non-empty level set is called a \itbf{sheet} of $\g$.
        \item $S$ is a \itbf{stabiliser sheet} if it is an irreducible component of a stabiliser level set.
        \item $S$ is a \itbf{centraliser sheet} if it is an irreducible component of a centraliser level set.
    \end{itemize}
\end{definition}

This terminology distinguishes stabiliser sheets as the type that have been studied so far in the literature.
Such a change in nomenclature will allow us to uniformly state certain results about both types of sheet, while also highlighting differences (see \S\ref{Subsection Levi-Type Sheets}).

Each stabiliser sheet lies in a unique stabiliser level set, and similarly, each centraliser sheet lies in a unique centraliser level set.
Given a (non-empty) level set $\ls$, we say that $S$ is a \emph{sheet of} $\ls$ if $S$ is an irreducible component of $\ls$; that is, $S \in \mathrm{Irr}(\g_{\am})$.

\subsection{Properties of Sheets}

Many of the results in this section have already been established for stabiliser sheets, sometimes with additional assumptions on the characteristic (see \cite{BK79}, \cite{B81}, and \cite{S82}, for example).
However, we can now extend these to arbitrary sheets in all characteristics.

We note that, if $\Lie\left(\C_Gx\right) = \fc_{\g}x$ for all $x \in \g$ (such as when assuming the Standard Hypotheses), then $\g_{(m)} = \g_{[m]}$ for all $m \in \N$, implicating that stabiliser sheets and centraliser sheets coincide in this case.

\begin{lemma}\label{Lemma Centre is a Sheet}
    $\z(\g) = \g_{(\dim G)} = \g_{[\dim G]}$, and thus $\z(\g)$ is both a stabiliser sheet and a centraliser sheet.
\begin{proof}
    It follows immediately from the definitions that $\g_{[\dim G]} = \z(\g)$.
    Using \cite[Lemma 13.3.3(i)]{S98}, $\C_G x = G$ if and only if $x \in \z(\g)$, from which $\g_{(\dim G)} = \z(\g)$ is clear.
    Since $\z(\g) = \g_{(\dim G)} = \g_{[\dim G]}$ is irreducible, the result follows.
\end{proof}
\end{lemma}

For each level set $\ls$, \textsc{Lemma}~\ref{Lemma Constant Dimension} and \textsc{Theorem}~\ref{Theorem Initial Decomposition Class Properties} imply that $\fD_{\am}[G]$ is a finite collection of irreducible constructible subsets of $\g$ such that $\ls$ is equal to their disjoint union.
Therefore, we can apply the results in \S\ref{Subsection Closure Orders} to the closure order on $\fD_{\am}[G]$.

Observe that each decomposition class is a subset of at least one stabiliser sheet and at least one centraliser sheet.
Using \textsc{Proposition}~\ref{Proposition Unique Dense in Irreducible Component}, we know that each sheet $S$ of $\ls$ contains a unique dense decomposition class, which we denote $\fD_S$.

\begin{lemma}\label{Lemma Closure of Sheet in Level Set}
    Suppose $S$ is a sheet of $\ls$.
\begin{itemize}
    \item[$\mathrm{(i)}$] $S = \bbKx S = \overline{S} \cap \ls = \overline{\bbKx S} \cap \ls = \overline{\bbK S} \cap \ls$. 
    \item[$\mathrm{(ii)}$] $S = \overline{G \cdot S} \cap \ls = \overline{\fD_S} \cap \ls$ is locally closed and $G$-stable. 
    \item[$\mathrm{(iii)}$]
        \begin{itemize}[leftmargin=18pt]
            \item[$\mathrm{(a)}$] If $S$ is a stabiliser sheet, then $S = \overline{S}^{\,\reg} = \overline{\fD_S}^{\,\reg} = \overline{\bbKx S}^{\,\reg} = \overline{\bbK S}^{\,\reg}$.
            \item[$\mathrm{(b)}$] If $S$ is a centraliser sheet, then $S = \overline{S}^{\,\greg} = \overline{\fD_S}^{\,\greg} = \overline{\bbKx S}^{\,\greg} = \overline{\bbK S}^{\,\greg}$.
        \end{itemize}
\end{itemize}
\begin{proof}
    Both $\bbKx$ and $S$ are irreducible, and thus so is $\bbKx S$.
    Since $\bbKx \ls =\nobreak \ls$, we have that $\bbKx S \subseteq \ls$, and hence $S = \bbKx S$ by maximality.
    Since $S \subseteq \ls$ is closed, it follows that $S = \overline{S} \cap \ls = \overline{\bbKx S} \cap \ls$.
    For (i), the last equality remains to be proven.
    Observe that $\bbK S = \{0\} \cup \bbKx S$, and so $\overline{\bbK S} = \{0\} \cup \overline{\bbKx S}$.
    If $0 \notin \ls$, then $\overline{\bbK S} \cap \ls = \overline{\bbK^{\times} S} \cap \ls$.
    Otherwise, $0 \in \ls$, and so \textsc{Lemma}~\ref{Lemma Centre is a Sheet} implies that $\ls = S = \z(\g)$, and hence $\overline{\bbK S} \cap \ls = \z(\g) = \overline{\bbK^{\times} S} \cap \ls$.

    For (ii), since $S = \overline{S} \cap \ls$ and $\ls$ is locally closed, it follows that $S$ is likewise locally closed.
    Since $\overline{G \cdot S} \cap \ls$ is also irreducible in $\ls$, maximality implies that $S = \overline{G \cdot S} \cap \ls$, which is $G$-stable.
    Moreover, $\overline{S} = \overline{\fD_S}$, and so (i) implies that $S = \overline{\fD_S} \cap \ls$.
    
    If $S$ is a stabiliser sheet, then \textsc{Lemma}~\ref{Lemma Closure inside a Level Set}(a) implies that $\overline{S} \cap \ls = \overline{S}^{\,\reg}$, and so the first three equalities in (iii)(a) follow from (i) and knowing that $\overline{S} = \overline{\fD_S}$.
    A similar argument using \textsc{Lemma}~\ref{Lemma Closure inside a Level Set}(b) shows the same for (iii)(b).
    The final equality in both follows from $\overline{\bbK S} = \{0\} \cup \overline{\bbKx S}$.
\end{proof}
\end{lemma}

\subsection{The Closure Order within Level Sets}\label{Subsection The Closure Order within Level Sets}

For each level set $\ls$, we can consider the Hasse diagram $\bm{\Gamma}\fD_{\am}[G]$.
This is precisely the subgraph of $\bm{\Gamma}\fD[G]$ induced by the subset of vertices $\fD_{\am}[G] \subseteq \fD[G]$.
We say that $\fJ \in \fD_{\am}[G]$ is \emph{maximal in} $\ls$ if it is maximal in $\fD_{\am}[G]$ (as defined in \S\ref{Subsection Partial Orders and Hasse Diagrams}), and similarly for \emph{minimal} and \emph{isolated}.
The following result is then given by \textsc{Proposition}~\ref{Proposition Irreducible Components and Maximal Elements}.

\begin{corollary}\label{Corollary Irreducible Components of Level Sets}
    Suppose $\ls$ is a level set, and $\fJ \in \fD_{\am}[G]$.
\begin{itemize}
    \item[$\mathrm{(i)}$] $\fJ$ is maximal in $\ls$ if and only if $\fJ = \fD_S$, for some sheet $S$ of $\ls$.  
    \item[$\mathrm{(ii)}$] The sheets of $\ls$ are in bijection with the maximal decomposition classes in $\ls$. 
    \item[$\mathrm{(iii)}$] If $\fJ$ coincides with a sheet of $\ls$, then it is isolated in $\ls$. 
\end{itemize}
\end{corollary}

This leads to an open question: does the converse of \textsc{Corollary}~\ref{Corollary Irreducible Components of Level Sets}(iii) hold in all cases?
If we make an additional assumption about the sheets of $\ls$, then we can use \textsc{Proposition}~\ref{Proposition Converse to Isolated Result} to prove the following.

\begin{proposition}\label{Proposition Coincide with Irreducible Component iff Isolated}
    Suppose $\ls$ is a level set, and assume that every sheet of $\ls$ is a union of decomposition classes.
    Then a decomposition class $\fJ \in \fD_{\am}[G]$ coincides with a sheet of $\ls$ if and only if it is isolated in $\ls$.  
\end{proposition}

We come back to the assumption required for \textsc{Proposition}~\ref{Proposition Coincide with Irreducible Component iff Isolated} in \textsc{Proposition}~\ref{Proposition Coincide with Levi-Type Sheets iff Isolated}.
In particular, it always holds if the characteristic is good for $G$.
On the other hand, \textsc{Lemma}~\ref{Lemma Centre is a Sheet} shows that \textsc{Proposition}~\ref{Proposition Coincide with Irreducible Component iff Isolated} always holds for the level set $\mathfrak{g}_{(\dim G)} = \mathfrak{g}_{[\dim G]} = \mathfrak{z}(\mathfrak{g})$.
Moreover, we shall see in \textsc{Appendix}~\ref{Appendix: Example} that \textsc{Proposition}~\ref{Proposition Coincide with Irreducible Component iff Isolated} holds for $G = \mathrm{Sp}_4$, even in bad characteristic.

\subsection{Sheets of Regular Elements}

Using \textsc{Theorem}~\ref{Theorem Initial Decomposition Class Properties}(v) and \textsc{Corollary}~\ref{Corollary Decomposition Varieties Properties}(i), we know that $\J{G}0 = \z(\g) \subseteq \overline{\fJ}$ for each $\fJ \in \fD[G]$.
Therefore, $\z(\g)$ is the unique minimal decomposition class in the closure order on $\fD[G]$.

Since $\g$ is the finite union of its decomposition varieties, all of which are closed irreducible subsets, there exists $\fJ \in \fD[G]$ such that $\overline{\fJ} = \g$, and any such decomposition class is necessarily maximal in the closure order on $\fD[G]$.
Using the antisymmetry of the closure order, we know such a decomposition class is unique, and we denote it $\J{G}^{\mathrm{max}}$.
It is clear from \textsc{Corollary}~\ref{Corollary Closure Order for Decomposition Classes} that $\J{G}^{\mathrm{max}}$ is the unique decomposition class $\fJ \in \fD[G]$ with $\dim \fJ = \dim \g$.

\begin{proposition}\label{Proposition Regular Sheets}
    ~
    \begin{itemize}
        \item[$\mathrm{(a)}$] $\g^\reg$ is a stabiliser sheet.
        \item[$\mathrm{(b)}$] $\g^\greg$ is a centraliser sheet.
    \end{itemize}
\begin{proof}
    Let $m \coloneqq \dim \C_G y$, for any $y \in \J{G}^{\mathrm{max}}$.
    For any $x \in \g$, \textsc{Lemma}~\ref{Lemma Stabiliser/Centraliser Dimension and Closure Order}(a) implies that $\dim \C_G x \geq \dim \C_G y$, and thus $\J{G}^{\mathrm{max}} \subseteq \g_{(m)} = \g^\reg$.
    It follows that $\J{G}^{\mathrm{max}}$ is the unique maximal decomposition class in $\g_{(m)}$, and so \textsc{Corollary}~\ref{Corollary Irreducible Components of Level Sets}(ii) implies that $\g_{(m)} = \g^\reg$ is a single sheet.
    This proves (a), and an analogous argument proves (b).
\end{proof}
\end{proposition}

It follows that \textsc{Proposition}~\ref{Proposition Coincide with Irreducible Component iff Isolated} always holds for the level sets $\mathfrak{g}_{(m)}$ and $\mathfrak{g}_{[m']}$, where $m \coloneqq \dim \C_G y$ and $m' \coloneqq \dim \fc_{\mathfrak{g}} y$, for any (and all) $y \in \J{G}^{\mathrm{max}}$.

The \emph{regular semisimple elements} of $\g$ are the semisimple elements $y \in \g$ such that $\dim G \cdot y$ is maximal -- or equivalently, $\dim \C_G y$ is minimal -- across all semisimple elements of $\g$.
For any maximal torus $T \subseteq G$, the set of regular semisimple elements of $\g$ is precisely $G \cdot \ft^\reg$.
The \emph{regular nilpotent elements} are defined analogously, and we let $\cN_{\g}^{\raisebox{2pt}{$\scriptstyle\,\reg$}} \coloneqq (\cN_{\g})^\reg$ denote the set of all such elements, which necessarily forms a single $G$-orbit.

\begin{theorem}\label{Theorem Maximal Decomposition Class}
    Suppose $z$ is a regular semisimple element of $\g$, and $e_0 \in \cN_{\mathfrak{c}^{\vphantom{\circ}}_{\g}z}^{\raisebox{2pt}{$\scriptstyle\,\reg$}}$ is a regular nilpotent element of $\fc_{\g}z$.
    Then $\J{G}^{\mathrm{max}} = \J{G}(z+e_0)$.
\begin{proof}
    Suppose $x \in \J{G}^{\mathrm{max}}$, and let $M \coloneqq \Co{G}x\s \in \mathfrak{C}_G$.
    Since $\dim \J{G}^{\mathrm{max}} = \dim \g$, it follows from \textsc{Corollary}~\ref{Corollary Alternative Dimension Formula} that $\dim \C_M x\n = \dim \C_G x = \dim \fd_{\g}x\s = \dim \z(\fm)$.
    Therefore, $\dim M \cdot x\n = \dim \fm - \dim \z(\fm)$.
    Fix a maximal torus $T \subseteq M \subseteq G$, let $\Phi_M \coloneqq \Phi(M,T)$ and $\Phi_G \coloneqq \Phi(G,T)$, and assume without loss of generality that $x_{\mathrm{s}} \in \mathfrak{t}$.
    Using \cite[Theorem 6.4]{J04}, we know that $\dim \cN_{\fm} = \lvert \Phi_M \rvert$, thus $\dim M \cdot x\n \leq \lvert \Phi_M \rvert$.
    Since $\dim M \cdot x\n = \dim \fm - \dim \z(\fm) \geq \dim \fm - \dim \ft = \lvert \Phi_M \rvert$, it follows that $\dim M \cdot x\n = \dim \cN_{\fm}$, and thus $x\n \in \cN_{\fm}^{\raisebox{2pt}{$\scriptstyle\,\reg$}}$ is a regular nilpotent element of $\fm$.

    Moreover, $\ft = \z(\fm) = \bigcap_{\alpha \in \Phi_M} \ker(\md \alpha)$, where the second equality is given in \cite[\S 2.3.2]{L05}.
    The differential map $\mathrm{d} \colon \mathrm{X}(T) \rightarrow \mathfrak{t}^*$ is $\mathbb{Z}$-linear with kernel $p\mkern 1.5mu \mathrm{X}(T)$, hence $\Phi_G \cap p\mkern 1.5mu \mathrm{X}(T) = \left\{ \alpha \in \Phi_G \midd \mathrm{d}\alpha = 0 \right\}$.
    If $\alpha \in \Phi_M$, then $\mathfrak{t} \subseteq \ker(\mathrm{d}\alpha)$, so $\mathrm{d}\alpha = 0$, and thus $\Phi_M \subseteq \Phi_G \cap p\mkern 1.5mu \mathrm{X}(T)$.
    Conversely, if $\alpha \in \Phi_G \cap p\mkern 1.5mu \mathrm{X}(T)$, then in particular $\mathrm{d}\alpha(x_{\mathrm{s}}) = 0$, so $\alpha \in (\Phi_G)_{x_{\mathrm{s}}} = \Phi_M$ where the equality follows from \textsc{Proposition}~\ref{Proposition Double Centraliser Properties}(ii).
    Therefore, $\Phi_M = \Phi_G \cap p\mkern 1.5mu \mathrm{X}(T)$.
    For any $y \in \ft$, \cite[Lemma 3.7]{S75} implies that $M = \Co{G}x\s \subseteq \Co{G} y$, and thus $\dim G \cdot y \leq \dim G \cdot x\s$.
    Therefore, $x\s \in \ft^{\reg}$ is a regular semisimple element, from which the result is clear.
\end{proof}
\end{theorem}

It follows from the proof of \textsc{Theorem}~\ref{Theorem Maximal Decomposition Class} that $\g^\reg = \g_{(\dim T)}$ (for any maximal torus $T \subseteq G$), and thus \textsc{Proposition}~\ref{Proposition Regular Sheets}(a) implies that $\g_{(\dim T)}$ is a stabiliser sheet.

\section{Lusztig--Spaltenstein Induction}\label{Section Lusztig--Spaltenstein Induction}

The Lusztig--Spaltenstein induction of nilpotent orbits is already well-studied for connected reductive algebraic groups; see, for example, \cite[\S 2]{B81} (over an algebraically closed field of characteristic $0$), \cite[\S 7]{CM93} (over $\mathbb{C}$), and \cite[\S 2.1]{S82}.
However, as demonstrated in \cite[\S 2.2]{S82}, we need not be limited to considering only nilpotent orbits.
In this section, we shall establish that many of the known properties of nilpotent Lusztig--Spaltenstein induction can be extended to arbitrary orbits in all characteristics.

\subsection{Existence of the Induced Orbit}

Before we prove the existence of the induced orbit, we require the following result.
Although this can be deduced from \cite[Proposition 2.7.1(ii)]{L05}, no proof is outlined, so we provide one here.

\begin{proposition}\label{Proposition Semisimple Part Parabolic Conjugacy}
    Suppose $P \subseteq G$ is a parabolic subgroup, with Levi factor $L \subseteq P$, and let $y \in \fl$.
    If $x \in L \cdot y + \fu_{\p}$, then $x\s \in P \cdot y\s$.
\begin{proof}
    Fix a maximal torus $T \subseteq L$, and a Borel subgroup $B \subseteq G$ such that $T \subseteq B \subseteq P$.
    Since $B' = B \cap L \in \mathfrak{P}(L,T)$, \cite[Proposition 14.25]{B91} implies that (up to $L$-conjugacy) we have $x \in y + \fu_{\p}$ and $y \in \mathfrak{b}'$.
    Then there exists a unique $y_0 \in \ft$ such that $y \in y_0 + \fu_{\mathfrak{b}'} + \fu_{\p} = y_0 + \fu_{\mathfrak{b}}$.
    It follows from \cite[Lemma 2.6.6]{L05} that $y_0 + \fu_{\mathfrak{b}}$ is $U_B$-stable, and hence $B$-stable.

    Suppose $z \in y_0 + \fu_{\mathfrak{b}}$.
    Since $z\s$ is semisimple, there exists $b \in B$ such that $b \cdot z\s \in \ft$.
    Moreover, $\cN_{\mathfrak{b}} = \fu_{\mathfrak{b}}$, so $b \cdot z\n \in \fu_{\mathfrak{b}}$.
    Then $b \cdot z\s + b \cdot z\n \in y_0 + \fu_{\mathfrak{b}}$ implies that $b \cdot z\s = y_0$, and thus $z\s \in B \cdot y_0$.
    Applying this with $x,y \in y_0 + \fu_{\mathfrak{b}}$ yields $x\s \in B \cdot y\s$, from which the result follows.
\end{proof}
\end{proposition}

Let $L \subseteq G$ be a Levi subgroup, and $\mathcal{O} \in \smallslant{\fl}{L}$ an $L$-orbit.
For each parabolic $P \in \mathfrak{P}(G,L)$, \textsc{Proposition}~\ref{Proposition Semisimple Part Parabolic Conjugacy} shows that every element of $G \cdot (\mathcal{O} + \fu_{\p})$ has the same semisimple part (up to $G$-conjugation).
It follows that there exists a semisimple element $z \in \g$ such that $G \cdot (\mathcal{O} + \fu_{\p}) \subseteq G \cdot \left(z + \cN_{\mathfrak{c}^{\vphantom{\circ}}_{\g}z}\right)$, and thus \textsc{Corollary}~\ref{Corollary Semisimple Centraliser Description} implies that $G \cdot (\mathcal{O} + \fu_{\p})$ consists of finitely many $G$-orbits.
Since $G \cdot (\mathcal{O} + \fu_{\p})$ is $G$-stable and irreducible, it follows that it contains a unique dense $G$-orbit, which we denote $\Ind_{\fl,\p}^{\g} \mathcal{O}$.
We refer to $\Ind_{\fl,\p}^{\g} \mathcal{O}$ as the \emph{induced} orbit (corresponding to the choices of $L$, $P$, and $\mathcal{O}$).

Since $\Ind_{\fl,\p}^{\g} \mathcal{O} \subseteq G \cdot (\mathcal{O} + \fu_{\p})$ is $G$-stable, it follows that $\bigl(\Ind_{\fl,\p}^{\g} \mathcal{O}\bigr) \cap (\mathcal{O} + \fu_{\p}) \neq \emptyset$.
We also have the following description of the closure of the induced orbit.

\begin{lemma}\label{Lemma Closure of Induced Orbit}
    Suppose $L \subseteq G$ is a Levi subgroup, $P \in \mathfrak{P}(G,L)$, and $\mathcal{O} \in \smallslant{\fl}{L}$.
    Then $\overline{\Ind_{\fl,\p}^{\g} \mathcal{O}} = G \cdot \big(\:\! \overline{\mathcal{O}} + \fu_{\p} \big)$.
\begin{proof}
    Observe that $\overline{\mathcal{O} + \fu_{\p}} = \overline{\mathcal{O}} + \fu_{\p}$, that $\overline{\Ind_{\fl,\p}^{\g} \mathcal{O}} = \overline{G \cdot (\mathcal{O} + \fu_{\p})}$, and $\overline{\mathcal{O}} + \fu_{\p}$ is $P$-stable.
    Therefore, applying \textsc{Lemma}~\ref{Lemma Parabolic Property Result} with $X = \mathcal{O} + \fu_{\p}$ yields $\overline{G \cdot ( \mathcal{O} + \fu_{\p} )} = G \cdot \big(\:\! \overline{\mathcal{O}} + \fu_{\p} \big)$, from which the result follows.
\end{proof}
\end{lemma}

\subsection{Nilpotent Induction}

For now, we focus on nilpotent Lusztig--Spaltenstein induction.
Our first result is an alternative description of the closure of an induced nilpotent orbit.

\begin{lemma}\label{Lemma Closure of Induced Nilpotent Orbit}
    Suppose $L \subseteq G$ is a Levi subgroup, $P \in \mathfrak{P}(G,L)$, and $\mathcal{O} \in \smallslant{\cN_{\fl}}{L}$.
    If $Z \subseteq \z(\fl)$ is any subset with $0 \in Z$, then $\overline{\Ind_{\fl,\p}^{\g} \mathcal{O}} = \cN_{\g} \cap G \cdot \big( Z + \overline{\mathcal{O}} + \fu_{\p} \big)$.
\begin{proof}
    Suppose $x \in \cN_{\g} \cap G \cdot \big( Z + \overline{\mathcal{O}} + \fu_{\p} \big)$, and thus $x\s = 0$.
    Up to $G$-conjugacy, there exists $z \in Z$, and $y \in \overline{\mathcal{O}}$ such that $x \in z + y + \fu_{\p}$.
    Since $(z+y)\s = z$, \textsc{Proposition}~\ref{Proposition Semisimple Part Parabolic Conjugacy} shows that $x\s \in P \cdot z$, which implies that $z = 0$ and thus $x \in G \cdot \big(\:\! \overline{\mathcal{O}} + \fu_{\p} \big)$.
    The reverse inclusion, $G \cdot \big(\;\! \overline{\mathcal{O}} + \fu_{\p} \big) \subseteq \cN_{\g} \cap G \cdot \big( Z + \overline{\mathcal{O}} + \fu_{\p} \big)$, is immediate, and so the result follows from \textsc{Lemma}~\ref{Lemma Closure of Induced Orbit}.
\end{proof}
\end{lemma}

Next, we need some results regarding the scalar action of $\bbK$ on non-nilpotent orbits.
We note that the proof of part (i) is outlined in \cite[\S 1.2(b)]{S82}.

\begin{lemma}\label{Lemma Dimension of Null Fibre}
    Suppose $y \in \g$ is not nilpotent.
\begin{itemize}
    \item[$\mathrm{(i)}$] $\dim \big(\:\! \overline{G \cdot \bbK y} \cap \cN_{\g} \big) = \dim G \cdot y$. 
    \item[$\mathrm{(ii)}$] If $H \subseteq G$ is a closed subgroup with $y \in \h$, then $\dim H \cdot \bbK y = \dim H \cdot y + 1$. 
\end{itemize}
\begin{proof}
    Define $\mathcal{O} \coloneqq G \cdot y$, and let $\chi \colon \overline{\bbK\mathcal{O}} \rightarrow \g/\!\!/G$ denote the restriction of the adjoint quotient morphism $\g \rightarrow \g/\!\!/G$ to $\overline{\bbK\mathcal{O}}$.
    Then ${\chi}^{-1}(0) = \overline{\bbK\mathcal{O}} \cap \cN_{\g} \neq \emptyset$.
    Since $\bbKx$ and $\mathcal{O}$ are both locally closed, \cite[Corollary AG10.2]{B91} implies that $\bbKx\mathcal{O}$ is constructible, and thus there exists $V \subseteq \bbKx\mathcal{O}$ such that $V \subseteq \overline{\bbKx\mathcal{O}} = \overline{\bbK\mathcal{O}}$ is open and dense.

    Let $x \in \overline{\bbK\mathcal{O}}$ and $a \in \bbKx$.
    Using the proof of \cite[Proposition 7.13]{J04}, we have that ${\chi}^{-1}\big({\chi}(x)\big) = \overline{\bbK\mathcal{O}} \cap G \cdot \big(x\s + \cN_{\mathfrak{c}^{\vphantom{\circ}}_{\g}x\s}\big)$, and it follows from \cite[Lemma 2.10]{J04} that ${\chi}^{-1}\big({\chi}(ax)\big) =\nobreak a{\chi}^{-1}\big({\chi}(x)\big)$.
    Therefore, $\dim {\chi}^{-1}\big({\chi}(x)\big)$ is constant across all $x \in \bbKx\mathcal{O}$, particularly on $V \subseteq \overline{\bbK\mathcal{O}}$.
    Moreover, $\mathcal{O} \subseteq {\chi}^{-1}\big({\chi}(x)\big)$ for each $x \in \mathcal{O}$, and thus $\dim {\chi}^{-1}\big({\chi}(x)\big) \geq \dim \mathcal{O}$ for each $x \in V$.

    Hence, for all $z \in \overline{\bbK\mathcal{O}}$ and $x \in V$, we have that $\dim {\chi}^{-1}\big({\chi}(z)\big) \geq \dim {\chi}^{-1}\big({\chi}(x)\big) \geq \dim \mathcal{O}$.
    In particular, $\dim \big(\:\!\:\! \overline{\bbK\mathcal{O}} \cap \cN_{\g} \big) = \dim {\chi}^{-1}(0) \geq \dim \mathcal{O}$.
    Given that $\mathcal{O}$ is not nilpotent, we have that ${\chi}^{-1}(0) \subsetneq \overline{\bbK\mathcal{O}}$, and so $\dim {\chi}^{-1}(0) < \dim \overline{\bbK\mathcal{O}}$.
    Since $\dim \bbKx \mathcal{O} \leq \dim \mathcal{O} + 1$, it follows that $\dim \overline{\bbK\mathcal{O}} = \dim \mathcal{O} + 1$ and $\dim \big(\:\! \overline{\bbK\mathcal{O}} \cap \cN_{\g} \big) = \dim \mathcal{O}$.

    For (ii), let $T \subseteq G$ be a maximal torus with $y\s \in \ft$, and let $W = \mathrm{W}(G,T)$ be the corresponding Weyl group.
    Using \cite[\S 7.12]{J04}, we have that $G \cdot y\s \cap \ft = W \cdot y\s$, and thus $H \cdot y\s \cap\nobreak \ft \subseteq W \cdot y\s$.
    It follows that $\mathrm{F} \coloneqq \left\{ a \in \bbKx \midd ay \in H \cdot y \right\}$ is finite.
    The algebraic group $H \times\nobreak \bbKx$ acts on $\h$, and the orbit of $y$ under this action is $H \cdot \bbKx y$, so $\dim (H \cdot \bbKx y) = \dim (H \times \bbKx) -\nobreak \dim \C_{H \times \bbKx}y$.
    Since $\mathrm{F}$ is finite, we can deduce that $\Co{H \times \bbKx}y = \Co{H}y \times \{1\}$, and thus $\dim \C_{H \times \bbKx}y = \dim H -\nobreak \dim H \cdot y$, from which the result follows.
\end{proof}
\end{lemma}

We now provide a proof that nilpotent Lusztig--Spaltenstein induction is independent of the choice of parabolic, in arbitrary characteristic.
The analogous result for induced unipotent classes was first proved in \cite[Theorem 2.2]{LS79}, and alternatively in \cite[\S 3]{B81Uni}.
The result for nilpotent orbits in characteristic zero can be found in \cite[Satz 2.6]{B81} and \cite[Theorem 7.1.3]{CM93}, and under the Standard Hypotheses in \cite[Remark 2.7]{PS18}.

The outline of our argument was first presented in \cite[\S 2.1]{S82} and takes inspiration from \cite[\S 3]{B81Uni}.
Here we expand on Spaltenstein's original proof for the reader's convenience.
Moreover, our argument proves that the known formula for the dimension of an induced nilpotent orbit also holds in arbitrary characteristic.

\begin{theorem}\label{Theorem Parabolic-Independence and Dimension Formula}
    Suppose $L \subseteq G$ is a Levi subgroup, and $\mathcal{O} \in \smallslant{\cN_{\fl}}{L}$ is a nilpotent $L$-orbit.
    \begin{itemize}
    \item[$\mathrm{(i)}$] For any $P,Q \in \mathfrak{P}(G,L)$, we have that $\Ind_{\fl,\p}^{\g} \mathcal{O} = \Ind_{\fl,\mathfrak{q}}^{\g} \mathcal{O}$. 
    \item[$\mathrm{(ii)}$] $\dim \Ind_{\fl,\p}^{\g} \mathcal{O} = \dim G - \dim L + \dim \mathcal{O}$. 
\end{itemize}
\begin{proof}
    Let $P,Q \in \mathfrak{P}(G,L)$ be arbitrary.
    We shall prove both statements simultaneously via strong induction on $\dim G$, observing that the base case (when $\dim G = 0$) is trivial.
    If $\dim L = \dim G$, the result is again trivial, so we can assume that $\dim L < \dim G$.
    Since the statements in (i) and (ii) hold for $G$ if and only if they hold for its derived subgroup, we can assume that $G$ is semisimple.

    There exists a central isogeny $\pi \colon G \rightarrow G_{\mathrm{ad}}$ to the adjoint group of $G$.
    Using \cite[\S 22]{B91} and \cite[Proposition 2.7(a)]{J04}, alongside \textsc{Lemma}~\ref{Lemma Central Surjections preserve Nilpotent Closures}, we can establish that the statements in (i) and (ii) hold for $G$ if and only if they hold for $G_{\mathrm{ad}}$.
    Therefore, we can assume that $G$ is semisimple of adjoint type.

    Fix a maximal torus $T \subseteq L$, a choice of simple roots $\Pi_G \subseteq \Phi_G \coloneqq \Phi(G,T)$, and a subset of simple roots $\Pi_L \subsetneq \Pi_G$ such that $\Phi_L \coloneqq \Phi(L,T) = \bbZ\Pi_L \cap \Phi_G$.
    Let $\{ \omega_{\alpha}^\vee \bigm| \alpha \in \Pi_G \}$ be the corresponding fundamental coweights, which (since $G$ is of adjoint type) form a $\mathbb{Z}$-basis of $\mathrm{Y}(T)$.
    It follows that $\omega_{\alpha}^\vee \notin p\mathrm{Y}(T)$, and consequently $\md \omega_{\alpha}^\vee \neq 0$, for each $\alpha \in \Pi_G$.

    Fix $\alpha \in \Pi_G \setminus \Pi_L$, and define $z \coloneqq \md\omega_{\alpha}^\vee(1) \in \ft \setminus \{0\}$.
    For each $\beta = \sum_{\gamma \in \Pi^{\vphantom{\circ}}_G} n^{\vphantom{\circ}}_{\gamma} \gamma \in \Phi$, we have that $\md\beta(z) = 0$ if and only if $p \mid n_{\alpha}$, from which it follows that $\Phi_L \subseteq \bbZ \big(\Pi_G \setminus \{\alpha\}\big) \cap \Phi_G \subseteq \left\{ \beta \in \Phi_G \midd \md\beta(z) = 0 \right\}$.
    Using \cite[Lemma 3.7]{S75}, we then know that $L \subseteq H \coloneqq \Co{G}z$.
    Since $z \in \ft \subseteq \fl \subseteq \fc_{\g}z$, it follows that $z \in \z(\fl)$.
    As $G$ is semisimple, \cite[Proposition 6.20]{MT11} and \cite[Proposition 2.3.4]{L05} show that $\z(\g) = \Lie\big(\mathcal{Z}(G)\big) = 0$.
    Therefore, $z \neq 0$ implies that $\fc_{\g}z \neq \g$, and thus $H = \Co{G}z \subsetneq G$.

    Observe that $\Phi_L$ is a $\mathbb{Q}$-closed root subsystem of $\Phi_G$, hence also of $\Phi_H$.
    Using the results in \cite[\S II.1.2]{S71}, there exists a base of $\Phi_L$ which can be extended to a base of $\Phi_H$, from which it follows that $L$ is also a Levi subgroup of $H$.
    Since $P \in \mathfrak{P}(G,L)$, there exists a Borel subgroup $B \subseteq G$ such that $T \subseteq B \subseteq P$.
    Then $B \cap H$ is a Borel subgroup of $H$ such that $T \subseteq B \cap H \subseteq P \cap H$, from which it follows that $P \cap H \in \mathfrak{P}(H,L)$.
    Using \cite[Proposition 13.20]{B91}, the corresponding Lie algebra is $\Lie(P \cap H) = \mathfrak{p} \cap \mathfrak{h}$, and likewise $Q \cap H \in \mathfrak{P}(H,L)$ with $\Lie(Q \cap H) = \mathfrak{q} \cap \mathfrak{h}$.
    Since $\dim H <\nobreak \dim G$, our inductive hypothesis implies that $\Ind_{\fl,\p \cap \h}^{\h} \mathcal{O} = \Ind_{\fl,\mathfrak{q} \cap \h}^{\h} \mathcal{O}$ and $\dim \Ind_{\fl,\p \cap \h}^{\h} \mathcal{O} = \dim H -\nobreak \dim L + \dim \mathcal{O}$.
    
    Let $P_1 \coloneqq P \cap H$, $\mathcal{O}_1 \coloneqq \Ind_{\fl,\p^{\vphantom{\circ}}_1}^{\h} \mathcal{O} \in \smallslant{\cN_{\h}}{H}$, and $\fu_1 \coloneqq \fu_{\p^{\vphantom{\circ}}_1}$.
    If $x_0 \in \mathcal{O}$ and $x \in \mathcal{O}_1$, then $\dim \C_H x = \dim \C_L x_0 = \dim L - \dim \mathcal{O}$.
    We will now prove that, for any choice of $x \in \mathcal{O}_1$, we have the following equality:
    \begin{equation}\label{Equation Claim within Parabolic-Independence Proof}
        \overline{\Ind_{\fl,\p}^{\g} \mathcal{O}} = \overline{G \cdot \bbK(z+x)} \cap \cN_{\g}\tag{$\star$}
    \end{equation}

    Since $G \cdot \bbK(z+h \cdot x) = G \cdot \bbK(z+x)$, for any $h \in H$, we can assume (without loss of generality) that $x \in \mathcal{O}_1 \cap (\mathcal{O} + \fu_1) \neq \emptyset$.
    Define $y \coloneqq z + x \in \h$, and suppose that $a \in \bbKx$.
    Using \cite[Lemma 2.10]{J04} and \cite[Lemma 2.6.6]{L05}, we can see that $P \cdot ay \subseteq az + \mathcal{O} + \fu_{\p}$, and thus $P \cdot \bbKx y \subseteq \bbKx z + \mathcal{O} + \fu_{\p}$.

    Given that $\mathcal{O} + \fu_1$ is $P_1$-stable, we have that $\dim P_1 \cdot x \leq \dim \mathcal{O} + \dim \fu_1$, and therefore $\dim \C_{P^{\vphantom{\circ}}_1}x \geq \dim P_1 - \dim \fu_1 - \dim \mathcal{O} = \dim L - \dim\mathcal{O} = \dim \C_H x $.
    Then $P_1 \subseteq H$ implies that $\dim \C_{P^{\vphantom{\circ}}_1}x = \dim \C_H x = \dim P - \dim \fu_{\p} - \dim \mathcal{O}$.
    Since $y_{\mathrm{s}} = z$ and $y_{\mathrm{n}} = x$, it follows from \textsc{Lemma}~\ref{Lemma Stabiliser/Centraliser of Jordan Decomposition}(i) and $H = \Co{G}z$ that $\dim \C_P y = \dim \C_{P^{\vphantom{\circ}}_1} x$, and thus $\dim P \cdot y = \dim \mathcal{O} + \dim \mathfrak{u}_{\mathfrak{p}}$.

    Since $z \neq 0$, we know that $y \neq y\n$ is not nilpotent, and so it follows from \textsc{Lemma}~\ref{Lemma Dimension of Null Fibre}(ii) that $\dim P \cdot \bbKx y = \dim P \cdot y + 1$.
    Therefore, we have that $P \cdot \bbKx y \subseteq \bbKx z + \mathcal{O} + \fu_{\p}$ with $\dim P \cdot \bbKx y = \dim(\mathcal{O}+ \fu_{\p}) + 1 = \dim \big(\bbKx z + \mathcal{O} + \fu_{\p}\big)$.
    Since $P \cdot \bbKx y $ is irreducible, it follows that $\overline{P \cdot \bbKx y} = \overline{\bbKx z + \mathcal{O} + \fu_{\p}} = \bbK z + \overline{\mathcal{O}} + \fu_{\p}$, where the last equality makes use of \textsc{Proposition}~\ref{Proposition Closure of Semisimple + Nilpotent}.
    Applying \textsc{Lemma}~\ref{Lemma Parabolic Property Result} with $X = \bbKx y$ shows that $\overline{G \cdot \bbKx y} = G \cdot \big(\:\!\overline{P \cdot \bbKx y} \:\!\big) = G \cdot \big( \bbK z + \overline{\mathcal{O}} + \fu_{\p} \big)$.
    Since $0 \in \bbK z \subseteq \z(\fl)$, it follows from \textsc{Lemma}~\ref{Lemma Closure of Induced Nilpotent Orbit} that $\overline{ G \cdot \bbK y} \cap \cN_{\g} = G \cdot \big( \bbK z + \overline{\mathcal{O}} + \fu_{\p} \big) \cap \Ng = \overline{\Ind_{\fl,\p}^{\g} \mathcal{O}}$.
    This completes our proof of (\ref{Equation Claim within Parabolic-Independence Proof}). \\

    Since $\mathcal{O}_1 = \Ind_{\fl,\mathfrak{q} \cap \h}^{\h} \mathcal{O}$, (\ref{Equation Claim within Parabolic-Independence Proof}) also shows that $\overline{\Ind_{\fl,\p}^{\g} \mathcal{O}} = \overline{ G \cdot \bbK(z+x)} \cap \cN_{\g} = \overline{\Ind_{\fl,\mathfrak{q}}^{\g} \mathcal{O}}$, for any $x \in \mathcal{O}_1$, from which (i) is clear.
    
    For (ii), let $x_0 \in \mathcal{O}$ as before, and $x \in \mathcal{O}_1 \cap (\mathcal{O} + \fu_1)$ as in the proof of (\ref{Equation Claim within Parabolic-Independence Proof}).
    Since $\C_G(z+x) = \C_H x$ and $H \cdot x = \mathcal{O}_1$, it follows that $\dim \C_G (z+x) = \dim \C_H x = \dim \C_L x_0$.
    On the other hand, \textsc{Lemma}~\ref{Lemma Dimension of Null Fibre}(i) implies that $\overline{\Ind_{\fl,\p}^{\g} \mathcal{O}} = \overline{\bbK \big(G \cdot(z+x)\big)} \cap \cN_{\g}$ has the same dimension as $G \cdot(z+x)$.
    Therefore, $\dim \Ind_{\fl,\p}^{\g} \mathcal{O} = \dim G \cdot (z+x) = \dim G - \dim L + \dim \mathcal{O}$, which proves (ii).
\end{proof}
\end{theorem}

Consequently, for any Levi subgroup $L \subseteq G$ and $\mathcal{O} \in \smallslant{\cN_{\fl}}{L}$, we let $\Ind_{\fl}^{\g} \mathcal{O}$ denote the induced orbit $\Ind_{\fl,\p}^{\g} \mathcal{O}$ for any choice of parabolic $P \in \mathfrak{P}(G,L)$.
Moreover, \textsc{Theorem}~\ref{Theorem Parabolic-Independence and Dimension Formula}(ii) proves that nilpotent induction preserves codimension: if $x_0 \in \mathcal{O} \in \smallslant{\cN_{\fl}}{L}$ and $\tilde{x} \in \Ind_{\fl}^{\g} \mathcal{O}$, then $\dim \C_L x_0 = \dim \C_G \tilde{x}$.

\pagebreak

\subsection{Arbitrary Induction}

Suppose $L \subseteq G$ is a Levi subgroup, and $\mathcal{O} \in \smallslant{\fl}{L}$ is any \mbox{$L$-orbit}.
Recall that, for any $P \in \mathfrak{P}(G,L)$, we defined $\Ind_{\fl,\p}^{\g} \mathcal{O}$ as the unique dense $G$-orbit in $G \cdot (\mathcal{O} + \fu_{\p})$.
The following alternative construction of $\Ind_{\fl,\p}^{\g} \mathcal{O}$ can be found in \cite[\S 2.2]{S82}.

\begin{lemma}\label{Lemma Arbitrary Induced Orbit}
    Suppose $L \subseteq G$ is a Levi subgroup, $\mathcal{O} \in \smallslant{\fl}{L}$ is an arbitrary $L$-orbit, and $P \in \mathfrak{P}(G,L)$.
    If $x \in \mathcal{O}$, then $\Ind_{\fl,\p}^{\g} \mathcal{O} = G \cdot \big(x\s + \Ind_{\mathfrak{c}^{\vphantom{\circ}}_{\fl}x^{\vphantom{\circ}}\s}^{\mathfrak{c}^{\vphantom{\circ}}_{\g}x^{\vphantom{\circ}}\s} (\Co{L}x\s \cdot x\n) \big)$.
\begin{proof}
    Fix a maximal torus $T \subseteq L$ such that $x\s \in \ft$, and recall from \cite[Lemma 3.7]{S75} that $G_1 \coloneqq \Co{G}x\s$ is reductive.
    The cocharacter approach to parabolics is extremely classical, and can be found in (for example) \cite[\S2]{R88} and \cite[Proposition 8.4.5]{S98}.
    Using this approach -- specifically the results listed in \cite[\S 2.3]{BMR05} -- we can see that $P_1 \coloneqq \Co{P}x_{\mathrm{s}}$ is a parabolic subgroup of $G_1$, with Levi factor $L_1 \coloneqq \Co{L}x_{\mathrm{s}}$ and unipotent radical $U_1 \coloneqq \mathrm{U}_{P^{\vphantom{\circ}}_1} = U_P \cap G_1$.

    Let $\mathcal{O}_1 \coloneqq L_1 \cdot x\n \in \slant{\cN_{\fl^{\vphantom{\circ}}_1}}{L_1}$, fix $x_1 \in \bigl(\Ind_{\fl^{\vphantom{\circ}}_1}^{\g^{\vphantom{\circ}}_1} \mathcal{O}_1\bigr) \cap (\mathcal{O}_1 + \fu_{\p^{\vphantom{\circ}}_1}) \neq \emptyset$, and let $y \coloneqq x\s + x_1 \in \g_1$.
    Since $x_1 \in \g_1 = \fc_{\g}x\s$, it is clear that $y\s = x\s$ and $y\n = x_1$.
    Moreover, $y \in x\s + \mathcal{O}_1 + \fu_{\p_1} = L_1 \cdot x + \fu_{\p^{\vphantom{\circ}}_1} \subseteq L \cdot x + \fu_{\p} = \mathcal{O} + \fu_{\p}$, which implies that $G \cdot y \subseteq G \cdot (\mathcal{O} + \fu_{\p})$.

    It follows from \cite[Lemma 2.6.6]{L05} that $\mathcal{O} + \fu_{\p}$ is $P$-stable, and thus $\dim G \cdot y \leq \dim G \cdot (\mathcal{O} + \fu_{\p}) \leq (\dim G - \dim P) + \dim (\mathcal{O} + \fu_{\p}) = 2 \dim \fu_{\p} + \dim \mathcal{O}$.
    On the other hand, \textsc{Theorem}~\ref{Theorem Parabolic-Independence and Dimension Formula}(ii) implies that $\dim \C_{G^{\vphantom{\circ}}_1} x_1 = \dim \C_{L^{\vphantom{\circ}}_1} x_{\mathrm{n}}$, from which it follows that $\dim \C_G y = \dim \C_L x$.
    Thus $\dim G \cdot y = \dim G - \dim L + \dim \mathcal{O} = 2 \dim \fu_{\p} + \dim \mathcal{O} = \dim G \cdot (\mathcal{O} + \fu_{\p})$.
    Since $\Ind_{\mathfrak{l},\mathfrak{p}}^{\mathfrak{g}} \mathcal{O}$ is the unique dense orbit in $G \cdot (\mathcal{O} + \mathfrak{u}_{\mathfrak{p}})$, we therefore have that $\Ind_{\fl,\p}^{\g} \mathcal{O} = G \cdot y = G \cdot \big(x\s + \bigl(\Ind_{\fl^{\vphantom{\circ}}_1}^{\g^{\vphantom{\circ}}_1} \mathcal{O}_1\bigr) \cap (\mathcal{O}_1 + \fu_{\p^{\vphantom{\circ}}_1})\big) = G \cdot \big(x\s +\nobreak \Ind_{\fl^{\vphantom{\circ}}_1}^{\g^{\vphantom{\circ}}_1} \mathcal{O}_1 \big)$.
\end{proof}
\end{lemma}

As with nilpotent Lusztig--Spaltenstein induction, we can establish parabolic independence and calculate the dimension of the induced orbit.
The following is immediate from \textsc{Theorem}~\ref{Theorem Parabolic-Independence and Dimension Formula}(i) and the proof of \textsc{Lemma}~\ref{Lemma Arbitrary Induced Orbit}.

\begin{corollary}\label{Corollary Arbitrary Parabolic-Independence and Dimension Formula}
    Suppose $L \subseteq G$ is a Levi subgroup, and $\mathcal{O} \in \smallslant{\fl}{L}$ is an arbitrary $L$-orbit.
    \begin{itemize}
    \item[$\mathrm{(i)}$] For any $P,Q \in \mathfrak{P}(G,L)$, we have that $\Ind_{\fl,\p}^{\g} \mathcal{O} = \Ind_{\fl,\mathfrak{q}}^{\g} \mathcal{O}$.
    \item[$\mathrm{(ii)}$] $\dim \Ind_{\fl,\p}^{\g} \mathcal{O} = \dim G - \dim L + \dim \mathcal{O}$. 
\end{itemize}
\end{corollary}

Therefore, for any Levi subgroup $L \subseteq G$ and $\mathcal{O} \in \smallslant{\fl}{L}$, we let $\Ind_{\fl}^{\g} \mathcal{O}$ denote the \emph{induced} orbit $\Ind_{\fl,\p}^{\g} \mathcal{O}$ for any choice of parabolic $P \in \mathfrak{P}(G,L)$.
As before, \textsc{Corollary}~\ref{Corollary Arbitrary Parabolic-Independence and Dimension Formula}(ii) implies that induction preserves codimension: that is, $\dim \C_G \tilde{x} = \dim \C_L x$, for any $x \in \mathcal{O}$ and $\tilde{x} \in \Ind_{\fl}^{\g} \mathcal{O}$.

We also have transitivity of induction for arbitrary orbits, as the following result demonstrates; since the proof presented here does not rely on the transitivity of nilpotent induction, this also serves as a proof of that result in arbitrary characteristic.
The approach here is very similar to that used in \cite[Proposition 7.1.4]{CM93}, for the nilpotent case in characteristic $0$.
The result is also proved via a different method in \cite[\S 2.5]{PS18}, for the nilpotent case assuming the Standard Hypotheses.

\begin{theorem}\label{Theorem Transitivity of Arbitrary Induction}
    Suppose $L \subseteq M \subseteq G$ are nested Levi subgroups of $G$.
    Then, for any $L$-orbit $\mathcal{O} \in \smallslant{\fl}{L}$, we have that $\Ind_{\fl}^{\g} \mathcal{O} = \Ind_{\mathfrak{m}}^{\g} \Ind_{\fl}^{\mathfrak{m}} \mathcal{O}$. 
\begin{proof}
    Let $Q_1 \in \mathfrak{P}(M,L)$ and $Q_2 \in \mathfrak{P}(G,M)$.
    Then $P \coloneqq U_{Q^{\vphantom{\circ}}_2} \rtimes Q_1 \in \mathfrak{P}(G,L)$, with $U_P = U_{Q^{\vphantom{\circ}}_2} \rtimes U_{Q^{\vphantom{\circ}}_1}$, and thus $\fu_{\p} = \fu_{\mathfrak{q}^{\vphantom{\circ}}_2} \oplus \fu_{\mathfrak{q}^{\vphantom{\circ}}_1}$.
    There exists $y_1 \in \bigl(\Ind_{\fl}^{\mathfrak{m}} \mathcal{O}\bigr) \cap \left(\mathcal{O} + \fu_{\mathfrak{q}^{\vphantom{\circ}}_1} \right) \neq \emptyset$ and $y_2 \in \bigl(\Ind_{\mathfrak{m}}^{\g} \Ind_{\fl}^{\mathfrak{m}} \mathcal{O}\bigr) \cap \bigl(\Ind_{\fl}^{\mathfrak{m}} \mathcal{O} + \fu_{\mathfrak{q}^{\vphantom{\circ}}_2} \bigr) \neq \emptyset$.
    Since $\Ind_{\fl}^{\mathfrak{m}} \mathcal{O} = M \cdot y_1$, and $\fu_{\mathfrak{q}^{\vphantom{\circ}}_2}$ is $M$-stable, there exists $h \in M$ such that $h \cdot y_2 \in y_1 + \fu_{\mathfrak{q}^{\vphantom{\circ}}_2} \subseteq \mathcal{O} + \fu_{\mathfrak{q}^{\vphantom{\circ}}_1} + \fu_{\mathfrak{q}^{\vphantom{\circ}}_2} = \mathcal{O} + \fu_{\p}$.
    Therefore, $\Ind_{\mathfrak{m}}^{\g} \Ind_{\fl}^{\mathfrak{m}} \mathcal{O} = G \cdot y_2 \subseteq G \cdot \left( \mathcal{O} + \fu_{\p} \right)$.
    It follows from \textsc{Corollary}~\ref{Corollary Arbitrary Parabolic-Independence and Dimension Formula}(ii) that $\dim \Ind_{\mathfrak{m}}^{\g} \Ind_{\fl}^{\mathfrak{m}} \mathcal{O} = \dim \Ind_{\fl}^{\g} \mathcal{O}$, and thus $\Ind_{\mathfrak{m}}^{\g} \Ind_{\fl}^{\mathfrak{m}} \mathcal{O} = \overline{G \cdot ( \mathcal{O} + \fu_{\p} )}^{\,\reg} = \Ind_{\fl}^{\g} \mathcal{O}$.
\end{proof}
\end{theorem}

Suppose that $L \subseteq G$ is a Levi subgroup, $P \in \mathfrak{P}(G,L)$, and $\mathcal{O} \in \smallslant{\fl}{L}$.
When $\mathcal{O}$ is nilpotent, and the Standard Hypotheses are assumed, \cite[Lemma 3.2(3)]{ST23} shows that the intersection $\Ind_{\fl}^{\g} \mathcal{O} \cap (\mathcal{O} + \fu_{\p})$ is a single $P$-orbit by adapting the analogous result for unipotent classes found in \cite[Theorem 1.3(c)]{LS79}.
For the convenience of the reader, we present a proof of this result for the general case, which likewise takes inspiration from \cite[Theorem 1.3(c)]{LS79}.

\begin{theorem}\label{Theorem Intersection is a Single P-orbit}
    Suppose $L \subseteq G$ is a Levi subgroup and $\mathcal{O} \in \smallslant{\fl}{L}$.
    If $P \in \mathfrak{P}(G,L)$, then $\bigl(\Ind_{\fl}^{\g} \mathcal{O}\bigr) \cap (\mathcal{O} + \fu_{\p})$ is a single $P$-orbit.
\begin{proof}
    Suppose $y \in \bigl(\Ind_{\fl}^{\g} \mathcal{O}\bigr) \cap (\mathcal{O} + \fu_{\p}) \neq \emptyset$.
    It follows from \cite[Lemma 2.6.6]{L05} that $P \cdot y \subseteq \mathcal{O} + \fu_{\p}$, and hence $\dim P \cdot y \leq \dim (\mathcal{O} + \fu_{\p})$.
    Therefore, $\dim \C_G y \geq \dim \C_P y \geq \dim P - \dim \fu_{\p} - \dim \mathcal{O} = \dim \C_G y$, where the equality comes from \textsc{Corollary}~\ref{Corollary Arbitrary Parabolic-Independence and Dimension Formula}(ii).
    Then $\dim P \cdot y = \dim (\mathcal{O} + \fu_{\p})$, and so the irreducibility of $\mathcal{O} + \fu_{\p}$ implies that $\overline{P \cdot y} = \overline{\mathcal{O} + \fu_{\p}}$.
    It follows that $P \cdot y$ and $P \cdot z$ are open and dense in $\mathcal{O} + \mathfrak{u}_{\mathfrak{p}}$, for any $y,z \in \bigl(\Ind_{\mathfrak{l}}^{\mathfrak{g}} \mathcal{O}\bigr) \cap (\mathcal{O} + \mathfrak{u}_{\mathfrak{p}})$.
    Therefore, $P \cdot y \cap P \cdot z \neq \emptyset$, hence $z \in P \cdot y$, from which the result follows.
\end{proof}
\end{theorem}

This concludes our proof of all the properties of Lusztig--Spaltenstein induction that were claimed in \textsc{Theorem}~\ref{Theorem LS Induction Main Result}.

\subsection{Closures of Orbits}

We shall end this section by deriving an alternative description of the closure of an induced orbit, inspired by \cite[\S 2.2]{S82}, but we first need the following result about closures of Levi-type orbits---see \textsc{Definition}~\ref{Definition Levi-Type} later.

\begin{proposition}\label{Proposition Closure of Levi-Type Orbits}
    Suppose $x \in \g$ is such that $\Co{G}x\s \subseteq G$ is a Levi subgroup.
    Then $\overline{G \cdot x} = G \cdot \big(x\s + \overline{\Co{G}x\s \cdot x\n}\,\big)$. 
\begin{proof}
    Let $L \coloneqq \Co{G}x\s$, $\mathcal{O} \coloneqq L \cdot x\n$, and fix $P \in \mathfrak{P}(G,L)$.
    The map $y \mapsto x\s + y$ is a homeomorphism, hence $x\s + \overline{L \cdot y} = \overline{L \cdot (x\s + y)}$, for any $y \in \overline{\mathcal{O}}$.
    It follows from \cite[Lemma 2.6.6]{L05} that $P \cdot \left(x\s+y\right) = x\s + L \cdot y + \fu_{\p}$, and thus $\overline{P \cdot (x\s + y)} = x\s + \overline{L \cdot y} + \overline{\fu_{\p}}$.
    Since $\overline{\mathcal{O}}$ is a finite union of nilpotent $L$-orbits, it follows that $P \cdot \big(\:\!\overline{L \cdot x}\:\!\big) = \bigcup_{y \in \overline{\mathcal{O}}} P \cdot(x\s + y) = \linebreak x\s + \overline{\mathcal{O}} + \fu_{\p} = \overline{L \cdot x} + \fu_{\p}$ which is $P$-stable and closed.
    Therefore, $\overline{P \cdot x} \subseteq \overline{L \cdot x} + \mathfrak{u}_{\mathfrak{p}} = P \cdot \big(\mkern 1.5mu \overline{L \cdot x} \mkern 1.5mu\big)$, and thus $\overline{L \cdot x} \subseteq \overline{P \cdot x}$ implies that $\overline{P \cdot x} = P \cdot \big(\;\!\overline{L \cdot x}\;\!\big)$.
    Applying \textsc{Lemma}~\ref{Lemma Parabolic Property Result} with $X = \{x\}$ then shows that $\overline{G \cdot x} = G \cdot \big(\:\!\overline{L \cdot x} \:\!\big) = G \cdot \big(x\s + \overline{\mathcal{O}}\:\!\big)$.
\end{proof}
\end{proposition}

This leads to what we believe to be an open question: does \textsc{Proposition}~\ref{Proposition Closure of Levi-Type Orbits} remain true without the assumption that $\Co{G}x\s \subseteq G$ is a Levi subgroup?

\begin{corollary}\label{Corollary Closure of Induced Orbit}
    Suppose $L \subseteq G$ is a Levi subgroup, and $x \in \mathcal{O} \in \smallslant{\fl}{L}$.
    If $\Co{G}x\s \subseteq G$ is also a Levi subgroup, then $\overline{\Ind_{\fl}^{\g} \mathcal{O}} = G \cdot \Big(x\s + \overline{\Ind_{\mathfrak{c}^{\vphantom{\circ}}_{\fl}x^{\vphantom{\circ}}\s}^{\mathfrak{c}^{\vphantom{\circ}}_{\g}x^{\vphantom{\circ}}\s} (\Co{L}x\s \cdot x\n)}\,\Big)$. 
\begin{proof}
    
    For any $y \in\nobreak \Ind_{\mathfrak{c}^{\vphantom{\circ}}_{\fl}x^{\vphantom{\circ}}\s}^{\mathfrak{c}^{\vphantom{\circ}}_{\g}x^{\vphantom{\circ}}\s} (\Co{L}x\s \cdot x\n)$, \textsc{Lemma}~\ref{Lemma Arbitrary Induced Orbit} shows that $\Ind_{\fl}^{\g} \mathcal{O} = G \cdot (x\s + \Co{G}x\s \cdot y) = G \cdot (x\s + y)$.
    Since $y \in \cN_{\mathfrak{c}_{\g}x\s}$, we have that $(x\s + y)\n = y$.
    Then \textsc{Proposition}~\ref{Proposition Closure of Levi-Type Orbits} implies that $\overline{G \cdot (x\s + y)} = G \cdot \big(x\s +\nobreak \overline{\Co{G}x\s \cdot y}\:\!\big)$, from which the result follows.
\end{proof}
\end{corollary}

\section{Levi-Type Decomposition Classes}\label{Section Levi-Type Decomposition Classes}

As already mentioned, much of the existing literature on decomposition classes has been developed under the assumption that $p \geq 0$ is (at least) good for the connected reductive group $G$.
It is well-known (see \cite[\S1.2, Remark 1]{S82}, for example) that this is equivalent to the assertion that $\Co{G}y \subseteq G$ is a Levi subgroup -- or equivalently, $\fc_{\g}y \subseteq \g$ is a Levi subalgebra -- for each semisimple $y \in \g$.

\subsection{Levi-Type Elements}

Outside good characteristic, there exist semisimple elements $x \in \g$ such that $\Co{G}x \subseteq G$ is a Levi subgroup, and we shall see that the decomposition classes of such elements have certain desirable properties.
However, we note that not all Levi subgroups are of the form $\Co{G}x$, as will be expounded upon in \S\ref{Subsection Stabiliser-Type Levi Subgroups}.

\begin{definition}\label{Definition Levi-Type}
    An element $x \in \g$ (as well as its orbit $G \cdot x$, and its decomposition class $\J{G}x$) are all called \itbf{Levi-type} if $\Co{G}x\s \subseteq G$ is a Levi subgroup.
    A decomposition variety is called \itbf{Levi-type} if it is the closure of a Levi-type decomposition class. 
\end{definition}

Observe that for any $x \sim y$, we have that $\Co{G}x\s$ is a Levi subgroup if and only if $\Co{G}y\s$ is a Levi subgroup.
Therefore, the property of being Levi-type (or not) is consistent across orbits and decomposition classes, meaning that \textsc{Definition}~\ref{Definition Levi-Type} is well-defined.

We observe that $p \geq 0$ is good for $G$ if and only if every element of $\g$ is Levi-type--- equivalently, if and only if every orbit / decomposition class / decomposition variety is Levi-type.
We also note that \textsc{Proposition}~\ref{Proposition Closure of Levi-Type Orbits} establishes a description for the closure of Levi-type orbits.
For an example of an element which is not Levi-type, see \textsc{Appendix}~\ref{Appendix: Example}.

\subsection{Stabiliser-Type Levi Subgroups}\label{Subsection Stabiliser-Type Levi Subgroups}

Clearly, a decomposition class is Levi-type if and only if any (or all) of its decomposition data are of the form $\Le$, where $L \subseteq G$ is a Levi subgroup and $e_0 \in \cN_{\fl}$.
However, in general, not every such pair $\Le$ is a valid decomposition datum.

\begin{definition}\label{Definition Stabiliser-Type}
    A Levi subgroup $L \subseteq G$ is \itbf{stabiliser-type} if there exists $y \in \g$ such that $L = \Co{G}y$.
\end{definition}

It follows from \textsc{Proposition}~\ref{Proposition Property of Connected Reductive Regular Subgroups} that any $y \in \g$ as in \textsc{Definition}~\ref{Definition Stabiliser-Type} has to be semisimple, and moreover, an element of $\z(\fl)_{[\dim \fl]}$.
The following lemma, providing equivalent definitions for stabiliser-type, can be easily deduced from \textsc{Proposition}~\ref{Proposition Property of Connected Reductive Regular Subgroups} alongside \textsc{Proposition}~\ref{Proposition Double Centraliser Properties}.

\begin{lemma}\label{Lemma Stabiliser-Type Levis in Decomposition Data}
    If $L \subseteq G$ is a Levi subgroup, then the following are equivalent$\:\!:$
    \begin{itemize}
        \item[$\mathrm{(1)}$] $L$ is stabiliser-type.
        \item[$\mathrm{(2)}$] $\z(\fl)_{[\dim \fl]} = \z(\fl)^{\reg}$.
        \item[$\mathrm{(3)}$] For any $z \in \z(\fl)^\reg$, we have that $L = \Co{G}z$.
        \item[$\mathrm{(4)}$] For any $z \in \z(\fl)^\reg$, we have that $\fl = \fc_{\g}z$.
        \item[$\mathrm{(5)}$] $\Le$ is a decomposition datum for some decomposition class, for each $e_0 \in \cN_{\fl}$.
    \end{itemize}
\end{lemma}

Using \cite[Lemma 2.6.13(i)]{L05}, we can see that whenever $p \geq 0$ is good for $G$ and does not divide $\big\lvert \left(\slant{\mathrm{X}(T)}{\mathbb{Z}\Phi}\right)_{\mathrm{tor}}\big\rvert$, any Levi subgroup of $G$ is stabiliser-type.
An example of a Levi subgroup which is not stabiliser-type is $L = \mathrm{SL}_{1,1} \coloneqq \pigl\{ \begin{psmallmatrix} a & 0 \\ 0 & {a^{-1}} \end{psmallmatrix} \bigm| a \in \bbKx \pigr\} \subseteq \mathrm{SL}_2 = G$, with $p = 2$; here we have that $\z(\fl) = \z(\fl)_{[3]}$, and so $\z(\fl)_{[\dim \fl]} = \z(\fl)_{[1]} = \emptyset$.

If $L \subseteq G$ is a stabiliser-type Levi subgroup, and $e_0 \in \cN_{\fl}$, then we refer to such a pair $\Le$ as a \emph{Levi-type} ($G$-)\emph{decomposition datum}.
Its corresponding Levi-type decomposition class is denoted $\J{G}\Le \coloneqq G \cdot \left( \z(\fl)^\reg + e_0 \right)$.

\subsection{Levi-Type Decomposition Varieties}

We now establish a description of \mbox{Levi-type} decomposition varieties that was previously only proved under stricter assumptions: that $G$ is semisimple, and either $p = 0$ \cite{BK79} or $p \geq 0$ is very good \cite{B98}.
The structure of our argument is adapted from \cite[Lemma 3.5.1]{B98}.

It follows from \S\ref{Subsection Stabiliser-Type Levi Subgroups} that each Levi-type decomposition variety can be expressed in the form $\dv{G}{\Le} = \overline{G \cdot \left( \z(\fl)^\reg + e_0 \right)}$, where $L \subseteq G$ is a stabiliser-type Levi subgroup and $e_0 \in \cN_{\fl}$.

\begin{theorem}\label{Theorem Levi-type Decomposition Variety description}
    Suppose $\Le$ is a Levi-type decomposition datum, and $P \in \mathfrak{P}(G,L)$.
    Then $\dv{G}{\Le} = G \cdot \left(\z(\fl) + \overline{L \cdot e_0} + \fu_{\p}\right)$.  
\begin{proof}
    If $x \in \z(\fl)^\reg + e_0 \subseteq \fl$, then \textsc{Lemma}~\ref{Lemma Stabiliser-Type Levis in Decomposition Data} implies that $\Co{G}x\s = L$.
    Using \cite[Lemma 2.6.6]{L05}, we have that $U_P \cdot x = x\s + e_0 + \fu_{\p}$, and thus $U_P \cdot \big(\z(\fl)^\reg + e_0\big) = \z(\fl)^\reg + e_0 + \fu_{\p}$.
    Therefore, $P \cdot \big( \z(\fl)^\reg + e_0 \big) = \z(\fl)^\reg + L \cdot e_0 + \fu_{\p}$, and so it follows from \textsc{Proposition}~\ref{Proposition Closure of Semisimple + Nilpotent} that $\overline{P \cdot \big( \z(\fl)^\reg + e_0 \big)} = \z(\fl) + \overline{L \cdot e_0} + \fu_{\p}$.
    Applying \textsc{Lemma}~\ref{Lemma Parabolic Property Result} with $X = \z(\fl)^\reg + e_0$ shows that $\overline{G \cdot \big( \z(\fl)^\reg + e_0 \big)} = G \cdot \big(\;\! \overline{P \cdot ( \z(\fl)^\reg + e_0 )} \;\!\big) = G \cdot \big( \z(\fl) + \overline{L \cdot e_0} + \fu_{\p} \big)$, from which the result is clear.
\end{proof}
\end{theorem}

This proves the first equality in \textsc{Theorem}~\ref{Theorem Levi-Type Decomposition Varieties Main Result}(i).
Now that we have the description of Levi-type decomposition varieties provided by \textsc{Theorem}~\ref{Theorem Levi-type Decomposition Variety description}, we can prove that they are unions of decomposition classes.

\begin{proposition}\label{Proposition Levi-type Decomposition Variety as Union of Decomposition Classes}
    Suppose $\J{G}\Le \in \fD[G]$ is a Levi-type decomposition class, and let $x \in \dv{G}{\Le}$.
    Then $\J{G}{x} \subseteq \dv{G}{\Le}$.  
\begin{proof}
    It follows from \textsc{Theorem}~\ref{Theorem Levi-type Decomposition Variety description} that we can assume $x \in \z(\fl) + \overline{L \cdot e_0} + \fu_{\p}$.
    Since $\z(\fl) + \overline{L \cdot e_0} + \fu_{\p}$ is $P$-stable, we can use \textsc{Proposition}~\ref{Proposition Semisimple Part Parabolic Conjugacy} to further assume that $x\s \in \z(\fl)$.
    This implies that $\fl \subseteq \fc_{\g}x\s$, and hence $\fd_{\g}x\s \subseteq \fc_{\g} \fl = \z(\fl)$.
    Since $x\n = x - x\s \in \z(\fl) + \overline{L \cdot e_0} + \fu_{\p}$, it follows that $\dreg{\g}x\s + x\n \subseteq \z(\fl) + \overline{L \cdot e_0} + \fu_{\p} \subseteq \dv{G}{\Le}$, and so \textsc{Theorem}~\ref{Theorem Initial Decomposition Class Properties}(i) implies that $\J{G}x \subseteq \dv{G}{\Le}$.
\end{proof}
\end{proposition}

Therefore, a Levi-type decomposition variety coincides with the union of the decomposition classes that it intersects, which proves \textsc{Theorem}~\ref{Theorem Levi-Type Decomposition Varieties Main Result}(ii).
Making use of \textsc{Lemma}~\ref{Lemma Closure inside a Level Set} alongside \textsc{Lemma}~\ref{Lemma Constant Dimension}, it follows immediately that both $\dv{G}{\Le}^{\,\reg}$ and $\dv{G}{\Le}^{\,\greg}$ are also unions of decomposition classes, for a Levi-type decomposition class $\J{G}\Le$.
This generalises statements found in \cite[\S 3.5]{B81} (for characteristic $0$) and \cite[\S 3.1]{A25} (for good characteristic).
On the other hand, \textsc{Appendix}~\ref{Appendix: Example} demonstrates that $\overline{\mathfrak{J}}$, $\overline{\mathfrak{J}}^{\,\mathrm{reg}}$, and $\overline{\mathfrak{J}}^{\,\mathfrak{g}\text{-}\mathrm{reg}}$ are each unions of decomposition classes, for any $\mathrm{Sp}_4$-decomposition class $\mathfrak{J} \in \mathfrak{D}[\mathrm{Sp}_4]$, even those which are not Levi-type (which exist when $p = 2$).

For the next result, we must strengthen \textsc{Definition}~\ref{Definition Levi-Type}.
We say that $\fJ \in \fD[G]$ is a \emph{strongly Levi-type} decomposition class if each $x \in \overline{\fJ}$ is Levi-type.
Then the following is a direct consequence of \textsc{Corollary}~\ref{Corollary Locally Closed Decomposition Classes} and \textsc{Proposition}~\ref{Proposition Levi-type Decomposition Variety as Union of Decomposition Classes}.

\begin{corollary}\label{Corollary Strongly-Levi-Type Decomposition Classes are Locally Closed}
    If $\fJ \in \fD[G]$ is a strongly Levi-type decomposition class, then $\fJ \subseteq \g$ is locally closed. 
\end{corollary}

As noted before, if $p \geq 0$ is good for $G$, then every element of $\g$ is Levi-type, and so the following result is immediate from \textsc{Corollary}~\ref{Corollary Strongly-Levi-Type Decomposition Classes are Locally Closed}.

\begin{corollary}\label{Corollary Good Characteristic Decomposition Classes are Locally Closed}
    If $p \geq 0$ is good for $G$, then every $G$-decomposition class is locally closed.
\end{corollary}

This result is already known, but the only proofs we could find in the literature are for the characteristic $0$ setting.
We note that the following is currently an open problem: are decomposition classes locally closed in general, without the assumption of being strongly Levi-type?
We shall see in \textsc{Appendix}~\ref{Appendix: Example} that (for example) $\mathrm{Sp}_4$-decomposition classes are always locally closed, even in bad characteristic---within which some decomposition classes are not (strongly) Levi-type.

\subsection{Decomposition Varieties and Lusztig--Spaltenstein Induction}\label{Subsection Decomposition Varieties and Lusztig--Spaltenstein Induction}

Next, we generalise some of the results found in \cite[\S 3.1]{A25} (under the assumption of good characteristic) which link decomposition varieties with Lusztig--Spaltenstein Induction.
We note that the characteristic $0$ case was first proved by Borho in \cite[\S 3]{B81}.

\begin{theorem}\label{Theorem Decomposition Variety as Union of Orbit Closures}
    If $\J{G}\Le \in \fD[G]$ is Levi-type, then $\dv{G}{\Le} = \bigcup\limits_{z \in \z(\fl)} \overline{\Ind_{\fl}^{\g}L \cdot \left(z+e_0\right)}$.  
\begin{proof}
    Fix $P \in \mathfrak{P}(G,L)$, suppose that $z \in \z(\fl)$, and let $\mathcal{O} = L \cdot \left(z + e_0\right) = z + L \cdot e_0 \in\nobreak \smallslant{\fl}{L}$.
    Then $\mathcal{O} + \fu_{\p} \subseteq \z(\fl) + \overline{L \cdot e_0} + \fu_{\p}$, hence \textsc{Theorem}~\ref{Theorem Levi-type Decomposition Variety description} implies that $G \cdot \left(\mathcal{O} + \fu_{\p}\right) \subseteq \dv{G}{\Le}$.
    Moreover, $\Ind_{\fl}^{\g} \mathcal{O} = \big( G \cdot \left( \mathcal{O} + \fu_{\p} \right) \big)^\reg$, and thus $\overline{\Ind_{\fl}^{\g} \mathcal{O}} = \overline{ G \cdot \left( \mathcal{O} + \fu_{\p} \right)\vphantom{\Ind_{\fl}^{\g} \mathcal{O}}} \subseteq \overline{\J{G}\Le\vphantom{\Ind_{\fl}^{\g}L \cdot \left(z+e_0\right)}}$.
    Therefore, $\bigcup_{z \in \z(\fl)} \overline{\Ind_{\fl}^{\g}L \cdot \left(z+e_0\right)} \subseteq \overline{\J{G}\Le\vphantom{\Ind_{\fl}^{\g}L \cdot \left(z+e_0\right)}}$.

    Observe that $z + \overline{L \cdot e_0} + \fu_{\p} = \overline{z + L \cdot e_0 + \fu_{\p}} \subseteq \overline{G \cdot \left( z + L \cdot e_0 + \fu_{\p} \right)} = \overline{\Ind_{\fl}^{\g}L \cdot \left(z+e_0\right)}$, where the last equality uses \textsc{Lemma}~\ref{Lemma Closure of Induced Orbit}.
    Then $\z(\fl) + \overline{L \cdot e_0} + \fu_{\p} \subseteq \bigcup_{z \in \z(\fl)} \overline{\Ind_{\fl}^{\g}L \cdot \left(z+e_0\right)}$, and hence $G \cdot \left( \z(\fl) + \overline{L \cdot e_0} + \fu_{\p} \right) \subseteq \bigcup_{z \in \z(\fl)} \overline{\Ind_{\fl}^{\g}L \cdot \left(z+e_0\right)}$.
    Consequently, the result follows from \textsc{Theorem}~\ref{Theorem Levi-type Decomposition Variety description}.
\end{proof}
\end{theorem}

This proves the second equality in \textsc{Theorem}~\ref{Theorem Levi-Type Decomposition Varieties Main Result}(i).
We also have the following consequences, the first of which was stated in \cite[Proposition 3.1(a)]{B81} (for characteristic $0$), and the second of which appeared in \cite[\S 3.1]{A25} (for good characteristic) and \cite[Corollary 3.2]{B81} (for characteristic $0$).

\begin{corollary}\label{Corollary Regular Closure of Levi-type Decomposition Class}
    Suppose $\J{G}\Le \in \fD[G]$ is Levi-type.
\begin{itemize}
    \item[$\mathrm{(i)}$] $\dv{G}{\Le}^{\,\reg} = \bigcup\limits_{z \in \z(\fl)} \Ind_{\fl}^{\g}L \cdot \left(z+e_0\right)$.
    \item[$\mathrm{(ii)}$] $\dv{G}{\Le}^{\,\reg} \cap \Ng = \Ind_{\fl}^{\g} L \cdot e_0$.
\end{itemize}
\begin{proof}
    Fix a parabolic $P \in \mathfrak{P}(G,L)$, and let $z \in \z(\fl)$.
    It follows from \textsc{Corollary}~\ref{Corollary Arbitrary Parabolic-Independence and Dimension Formula}(ii) that $\dim \Ind_{\fl}^{\g}L \cdot \left(z+e_0\right) = \dim L \cdot e_0 + 2 \dim \fu_{\p}$.
    Therefore, $\dim \Ind_{\fl}^{\g}L \cdot \left(z+e_0\right)$ is constant across all $z \in \z(\fl)$, and so (i) follows from \textsc{Theorem}~\ref{Theorem Decomposition Variety as Union of Orbit Closures}.

    For any $y \in \Ind_{\fl}^{\g} L \cdot \left(z + e_0\right)$, \textsc{Lemma}~\ref{Lemma Arbitrary Induced Orbit} implies that $y\s \in G \cdot z$, and thus $y \in \Ng$ if and only if $z = 0$.
    Consequently, (ii) follows from (i).
\end{proof}
\end{corollary}

This completes the proof of \textsc{Theorem}~\ref{Theorem Levi-Type Decomposition Varieties Main Result}.

\subsection{Levi-Type Sheets}\label{Subsection Levi-Type Sheets}

Recall from \textsc{Proposition}~\ref{Proposition Unique Dense in Irreducible Component} that each sheet $S \subseteq \g$ contains a unique dense decomposition class $\fD_S$.
We then say a sheet $S \subseteq \g$ is \emph{Levi-type} if $\fD_S$ is a Levi-type decomposition class, and observe that, if we assume good characteristic, then every sheet is Levi-type.

\begin{corollary}\label{Corollary Levi-type Sheet as a Union}
    If $S \subseteq \g$ is a Levi-type sheet, then $S$ is a union of decomposition classes.
\begin{proof}
    As observed previously, \textsc{Proposition}~\ref{Proposition Levi-type Decomposition Variety as Union of Decomposition Classes} implies that (for any Levi-type decomposition class $\fJ$) both $\overline{\fJ}^{\,\reg}$ and $\overline{\fJ}^{\,\greg}$ are unions of decomposition classes.
    Applying this to $\fD_S$ shows that the result follows from \textsc{Lemma}~\ref{Lemma Closure of Sheet in Level Set}(iii).
\end{proof}
\end{corollary}

Recall that, in order to prove \textsc{Proposition}~\ref{Proposition Coincide with Irreducible Component iff Isolated}, we required the assumption that every sheet of a level set $\ls$ was a union of decomposition classes.
Since \textsc{Corollary}~\ref{Corollary Levi-type Sheet as a Union} shows that this holds if all the sheets of $\ls$ are Levi-type sheets, we can restate a version of \textsc{Proposition}~\ref{Proposition Coincide with Irreducible Component iff Isolated} as follows.

\begin{proposition}\label{Proposition Coincide with Levi-Type Sheets iff Isolated}
    Suppose $\ls$ is a level set, and assume that every sheet of $\ls$ is Levi-type.
    Then a decomposition class $\fJ \in \fD_{\am}[G]$ coincides with a sheet of $\ls$ if and only if it is isolated in $\ls$. 
\end{proposition}

Moreover, if the characteristic is good for $G$, then the assumption used in \textsc{Proposition}~\ref{Proposition Coincide with Levi-Type Sheets iff Isolated} always holds.

Since every sheet $S$ is $G$-stable by \textsc{Lemma}~\ref{Lemma Closure of Sheet in Level Set}(ii), we have that $S \cap \Ng$ is a finite (possibly empty) union of nilpotent orbits.
The following property of Levi-type stabiliser sheets is an immediate consequence of \textsc{Lemma}~\ref{Lemma Closure of Sheet in Level Set}(iii)(a) alongside \textsc{Corollary}~\ref{Corollary Regular Closure of Levi-type Decomposition Class}(ii).

\begin{corollary}\label{Corollary Levi-type Sheet Nilpotent Orbit}
    If $S \subseteq \g$ is a Levi-type stabiliser sheet, then $S$ contains a unique nilpotent orbit.
    In particular, if $\fD_S = \J{G}\Le$, then $S \cap \Ng = \Ind_{\fl}^{\g} L \cdot e_0$. 
\end{corollary}

\subsection{Spaltenstein's Conjecture}

We shall conclude this section by drawing attention to connections between \textsc{Corollary}~\ref{Corollary Levi-type Sheet Nilpotent Orbit} and the following conjecture of Spaltenstein's, which we have reworded slightly in line with the new terminology introduced in \textsc{Definition}~\ref{Definition Sheets}.

\begin{conjecture}[{\cite[\S 1.2]{S82}}]\label{Conjecture Spaltenstein Sheets Contain Unique Nilpotent}
    For any connected reductive algebraic group $G$ \textup{(}over an algebraically closed field of arbitrary characteristic\textup{)}, every stabiliser sheet of $\g$ contains exactly one nilpotent orbit.
\end{conjecture}

In \cite[\S 1.2(c)]{S82}, Spaltenstein establishes that every stabiliser sheet contains at least one nilpotent orbit, and observes that Borho--Kraft \cite{BK79} (although working in characteristic $0$) essentially prove \textsc{Conjecture}~\ref{Conjecture Spaltenstein Sheets Contain Unique Nilpotent} when the characteristic $p \geq 0$ is good for $G$.
Moreover, \mbox{Spaltenstein} proves in \cite[Theorem 2.8]{S82} that \textsc{Conjecture}~\ref{Conjecture Spaltenstein Sheets Contain Unique Nilpotent} holds when $G$ has no simple components of exceptional type.
In later work, they prove that \textsc{Conjecture}~\ref{Conjecture Spaltenstein Sheets Contain Unique Nilpotent} is also true when $G$ is a simple algebraic group of either type $\mathrm{E}_6$ \cite[\S 7, Corollary 3]{S83}, or type $\mathrm{F}_4$ when $p = 2$ \cite[\S 5, Theorem]{S84}.

However, it is noted in \cite[\S 3.1]{PS18} that \textsc{Conjecture}~\ref{Conjecture Spaltenstein Sheets Contain Unique Nilpotent} remains open for certain bad characteristics.
It follows from \textsc{Corollary}~\ref{Corollary Levi-type Sheet Nilpotent Orbit} that \textsc{Conjecture}~\ref{Conjecture Spaltenstein Sheets Contain Unique Nilpotent} is at least true for Levi-type stabiliser sheets (regardless of characteristic). \\

We remark that \textsc{Conjecture}~\ref{Conjecture Spaltenstein Sheets Contain Unique Nilpotent} is false (in general) for centraliser sheets.
To realise this, it suffices to show that there exist non-empty centraliser level sets that contain no nilpotent orbits.

For example, consider $G = \mathrm{PGL}_2$ with $p = 2$, and let $\pi \colon \GL_2 \rightarrow \mathrm{PGL}_2$ be the canonical quotient homomorphism.
The nilpotent cone $\Ng$ consists of only two nilpotent $G$-orbits: the zero orbit, and $G \cdot x$, where $x = \md \pi \begin{psmallmatrix} 0 & 1 \\ 0 & 0 \end{psmallmatrix}$.
Simple computation reveals that we have $\fc_{\g} x = \big\{ \md \pi \begin{psmallmatrix} a & b \\ c & a \end{psmallmatrix} \bigm| a,b,c \in \bbK \big\}$, and thus $\dim \fc_{\g} x = 2$, while $\dim \fc_{\g} 0 = 3$.
Similar computation shows that $\md \pi \begin{psmallmatrix} 1 & 0 \\ 0 & 0 \end{psmallmatrix} \in \g_{[1]}$, hence $\g_{[1]} \neq \emptyset$ is a non-empty centraliser level set containing no nilpotent orbits.
Moreover, all the sheets of $\g$ are Levi-type, and so this also proves that \textsc{Corollary}~\ref{Corollary Levi-type Sheet Nilpotent Orbit} is false for Levi-type centraliser sheets.

\section{The Edges of the Hasse Diagram}\label{Section The Edges of the Hasse Diagram}

We conclude this paper by proving \textsc{Theorem}~\ref{Maintheorem Hasse Diagram Edges}, which describes the edges of the Hasse diagram $\bm{\Gamma}\fD[G]$ for the closure order $\boldpreceq$ on $\fD[G]$.

\subsection{Determining the Closure Order}

Suppose $M = \Co{G}x\s \in \mathfrak{C}_G$.
Then $\fm = \fc_{\g}x\s$ and $\z(\fm) = \fc_{\g}\fm = \fd_{\g}x\s$.
It follows from \textsc{Proposition}~\ref{Proposition Double Centraliser Properties} that $\fc_{\g}\fm = \z(\fm)$ and $\fc_{\g}\z(\fm) = \fm$.
Therefore, $\fm$ (and thus also $M$) is recoverable from the centre $\z(\fm)$.
In particular, if $M' \in \mathfrak{C}_G$ with $M \subsetneq M'$, then we must have that $\z(\fm') \subsetneq \z(\fm)$.

We equip the set of nilpotent $M$-orbits $\smallslant{\cN_{\fm}}{M}$ with the closure order, denoted $\boldpreceq$.
As a result of \textsc{Corollary}~\ref{Corollary Semisimple Centraliser Description} and \textsc{Lemma}~\ref{Lemma Closure Order is a Partial Order}, this is a finite poset, and we let $\boldprec$ and $\boldprecdot$ denote the corresponding strict partial order and covering relation, respectively.

\begin{lemma}\label{Lemma Nilpotent Order Inducing Decomposition Class Order}
    Suppose $M \in \mathfrak{C}_G$ and $\mathcal{O}, \mathcal{O}' \in \smallslant{\cN_{\fm}}{M}$.
\begin{itemize}
    \item[$\mathrm{(i)}$] For any $z \in \z(\fm)^\reg$ and $e_0 \in \mathcal{O}$, we have $\dim G \cdot \left(z + e_0\right) = \dim \mathcal{O} + \dim G - \dim M$. 
    \item[$\mathrm{(ii)}$] If $\mathcal{O}' \boldprec \mathcal{O}$, then $\J{G}\left(M;\mathcal{O}'\right) \boldprec \J{G}\left(M;\mathcal{O}\right)$. 
\end{itemize}
\begin{proof}
    Using $\Co{G}\left(z + e_0\right) = \Co{M}e_0$, we have that $\dim G - \dim G \cdot \left(z + e_0\right) = \dim M - \dim M \cdot e_0$.
    Since $M \cdot e_0 = \mathcal{O}$, (i) follows after rearranging.

    We then observe that $z + \overline{\mathcal{O}} \subseteq G \cdot \left( z + \overline{\mathcal{O}} \:\!\right) \subseteq \overline{G \cdot \left( z + \mathcal{O} \right)} \subseteq \overline{G \cdot \left(\z(\fm)^\reg + \mathcal{O} \right)} = \overline{\J{G}\left(M;\mathcal{O}\right)}$.
    Thus $\z(\fm)^\reg + \mathcal{O}' \subsetneq \z(\fm)^\reg + \overline{\mathcal{O}} = \bigcup_{z \in \z(\fm)^\reg}\left( z + \overline{\mathcal{O}} \:\!\right) \subseteq \overline{\J{G}\left(M;\mathcal{O}\right)}$.
    Therefore, we have that $\z(\fm)^\reg +\nobreak \mathcal{O}' \subsetneq \z(\fm)^\reg + \overline{\mathcal{O}} \subseteq \overline{\J{G}\left(M;\mathcal{O}\right)}$, which implies that $\J{G}\left(M;\mathcal{O}'\right) \boldpreceq \J{G}\left(M;\mathcal{O}\right)$.
    Since $\mathcal{O}' \boldprec \mathcal{O}$ entails that $\dim \mathcal{O}' < \dim \mathcal{O}$, it follows from (i) that $\dim \J{G}\left(M;\mathcal{O}'\right) < \dim \J{G}\left(M;\mathcal{O}\right)$.
    Hence $\J{G}\left(M;\mathcal{O}'\right) \neq \J{G}\left(M;\mathcal{O}\right)$, which proves (ii).
\end{proof}
\end{lemma}

Consequently, the closure order on $\smallslant{\cN_{\fm}}{M}$ determines the closure order between the relevant decomposition classes in $\fD[G]$.
We equip $\mathfrak{C}_G$ with the partial order $\subseteq$ defined by inclusion of sets, and let $\subsetneq$ and $\subsetneqdot$ denote the corresponding strict partial order and covering relation, respectively.
The following lemma demonstrates the analogous result to \textsc{Lemma}~\ref{Lemma Nilpotent Order Inducing Decomposition Class Order}(ii) for the partial order on $\mathfrak{C}_G$, under the assumption that we are working with stabiliser-type Levi subgroups.

\begin{lemma}\label{Lemma Semisimple Order Inducing Decomposition Class Order}
    Suppose $L,L' \in \mathfrak{C}_G$ are both $($stabiliser-type$)$ Levi subgroups with $L \subseteq L'$, and $\mathcal{O} \in \smallslant{\cN_{\fl}}{L}$.
    \begin{itemize}
        \item[$\mathrm{(i)}$] $\J{G}(L';\Ind_{\fl}^{\fl'} \mathcal{O}) \subseteq \overline{\J{G}(L;\mathcal{O})}^{\,\reg}$.
        \item[$\mathrm{(ii)}$] If $L \subsetneq L'$, then $\J{G}(L';\Ind_{\fl}^{\fl'} \mathcal{O}) \boldprec \J{G}(L;\mathcal{O})$.
    \end{itemize}
\begin{proof}
    Suppose $z \in \z(\fl')^\reg$, from which \textsc{Lemma}~\ref{Lemma Stabiliser-Type Levis in Decomposition Data} implies that $\fl' = \fc_\g z$, and note that $z \in \z(\fl)$.
    Then \textsc{Corollary}~\ref{Corollary Regular Closure of Levi-type Decomposition Class}(i) shows that $G \cdot ( z + \Ind_{\fl}^{\fl'} \mathcal{O} ) = \Ind_{\fl}^{\g} L \cdot (z + \mathcal{O}) \subseteq \overline{\J{G}(L;\mathcal{O})}^{\,\reg}$.
    It follows that $\J{G}(L';\Ind_{\fl}^{\fl'} \mathcal{O}) \subseteq \overline{\J{G}(L;\mathcal{O})}^{\,\reg}$, which proves (i).
    
    Moreover, (i) demonstrates that $\J{G}(L';\Ind_{\fl}^{\fl'} \mathcal{O}) \boldpreceq \J{G}(L;\mathcal{O})$.
    Using \textsc{Corollary}~\ref{Corollary Alternative Dimension Formula} twice, we know that $\dim \J{G}(L;\mathcal{O}) = \dim G - \dim L + \dim \mathcal{O} + \dim \z(\fl)$ and $\dim \J{G}(L';\Ind_{\fl}^{\fl'} \mathcal{O}) = \dim G - \dim L' + \dim \Ind_{\fl}^{\fl'} \mathcal{O} + \dim \z(\fl')$.
    Combining this with \textsc{Theorem}~\ref{Theorem Parabolic-Independence and Dimension Formula}(ii) shows that $\dim \J{G}(L';\Ind_{\fl}^{\fl'} \mathcal{O}) = \dim G - \dim L + \dim \mathcal{O} + \dim \z(\fl')$.
    If $L \subsetneq L'$, then $\z(\fl') \subsetneq \z(\fl)$, and thus $\dim \z(\fl') < \dim \z(\fl)$.
    Therefore, $\dim \J{G}(L';\Ind_{\fl}^{\fl'} \mathcal{O}) < \dim \J{G}(L;\mathcal{O})$, and so $\J{G}(L';\Ind_{\fl}^{\fl'} \mathcal{O}) \neq \J{G}(L;\mathcal{O})$, from which the result follows.
\end{proof}
\end{lemma}

\subsection{Types of Edges (assuming good characteristic)}

For the remainder of this paper, we assume that $p \geq 0$ is good for $G$.
It follows that every connected stabiliser of a semisimple element is a Levi subgroup, and therefore $\mathfrak{C}_G$ is precisely the set of stabiliser-type Levi subgroups.
We note that, if $G$ has no components of type $\mathrm{A}_{rp-1}$ (for all $r \in \Np$), then $p \geq 0$ is very good for $G$, and in this case $\mathfrak{C}_G$ is actually the set of all Levi subgroups.

Using the assumption of good characteristic, we can prove that each edge of the Hasse diagram $\bm{\Gamma}\fD[G]$ is either determined by the covering relation on $\mathfrak{C}_G$, or the covering relation on nilpotent orbits.

\begin{theorem}\label{Theorem Hasse Diagram Edges}
    Suppose $p \geq 0$ is good for $G$.
    Then every edge of $\bm{\Gamma}\fD[G]$ is of exactly one of the following types\emph{:}  

    \begin{itemize}
        \item[$\mathrm{(I)}$] An edge from $\J{G}(L';\Ind_{\fl}^{\fl'}\mathcal{O})$ to $\J{G}(L;\mathcal{O})$, for $L,L' \in \mathfrak{C}_G$ with $L \subsetneqdot L'$ and $\mathcal{O} \in \smallslant{\cN_{\fl}}{L}$.
        \item[$\mathrm{(II)}$] An edge from $\J{G}(L;\mathcal{O}')$ to $\J{G}(L;\mathcal{O})$, for $L \in \mathfrak{C}_G$ and $\mathcal{O},\mathcal{O}' \in \smallslant{\cN_{\fl}}{L}$ with $\mathcal{O}' \boldprecdot \mathcal{O}$.
    \end{itemize}

\begin{proof}

    Suppose $\fJ, \fJ' \in \fD[G]$ are such that $\fJ' \boldprecdot \fJ$, and $(L;\mathcal{O})$ is any saturated decomposition datum corresponding to $\fJ$.
    Since $p \geq 0$ is good, $\Co{G}z \subseteq G$ is a Levi subgroup for each $z \in \z(\fl)$.
    Applying \textsc{Corollary}~\ref{Corollary Closure of Induced Orbit} to \textsc{Theorem}~\ref{Theorem Decomposition Variety as Union of Orbit Closures}, we can see that $\overline{\fJ} = \bigcup_{z \in \z(\fl)} G \cdot \pigl(z + \overline{\Ind_{\fl}^{\mathfrak{c}_{\g}z}\mathcal{O}}\,\pigr)$.
    Therefore, there exists $z \in \z(\fl)$ and $e_0 \in \overline{\Ind_{\fl}^{\mathfrak{c}_{\g}z} \mathcal{O}}$ such that $z + e_0 \in \fJ'$.
    If we let $L' = \Co{G}z \in \mathfrak{C}_G$ and $\mathcal{O}' = L' \cdot e_0$, then $\fJ' = \J{G}(L';\mathcal{O}')$.

    First, suppose that $z \notin \z(\fl)^\reg$.
    Since $\fl \subsetneq \fc_{\g}z = \fl'$, we know that $L \subsetneq L'$, and we also have that $e_0 \in \overline{\Ind_{\fl}^{\fl'}\mathcal{O}}$.
    Suppose, for contradiction, that $e_0 \notin \Ind_{\fl}^{\fl'}\mathcal{O}$.
    Then $\mathcal{O}' \boldprec \Ind_{\fl}^{\fl'}\mathcal{O}$, and so \textsc{Lemma}~\ref{Lemma Nilpotent Order Inducing Decomposition Class Order}(ii) implies that $\J{G}(L';\mathcal{O}') \boldprec \J{G}(L';\Ind_{\fl}^{\fl'}\mathcal{O})$.
    Since $\J{G}(L';\Ind_{\fl}^{\fl'}\mathcal{O}) \boldprec \J{G}(L;\mathcal{O})$, by \textsc{Lemma}~\ref{Lemma Semisimple Order Inducing Decomposition Class Order}(ii), we have that $\fJ' \boldprec \J{G}(L';\Ind_{\fl}^{\fl'}\mathcal{O}) \boldprec \fJ$ which contradicts $\fJ' \boldprecdot \fJ$.
    Therefore, $e_0 \in \Ind_{\fl}^{\fl'}\mathcal{O}$, and thus $\fJ' = \J{G}(L';\Ind_{\fl}^{\fl'}\mathcal{O})$.

    Now suppose, for contradiction, that there exists $L_0 \in \mathfrak{C}_G$ such that $L \subsetneq L_0 \subsetneq L'$.
    Then \textsc{Theorem}~\ref{Theorem Transitivity of Arbitrary Induction} implies that $\Ind_{\fl}^{\fl'}\mathcal{O} = \Ind_{\fl_0}^{\fl'} \Ind_{\fl}^{\fl_0} \mathcal{O}$, and applying \textsc{Lemma}~\ref{Lemma Semisimple Order Inducing Decomposition Class Order}(ii) twice yields $\fJ' = \J{G}(L';\Ind_{\fl_0}^{\fl'} \Ind_{\fl}^{\fl_0} \mathcal{O}) \boldprec \J{G}(L_0; \Ind_{\fl}^{\fl_0} \mathcal{O}) \boldprec \J{G}(L;\mathcal{O}) = \fJ$.
    This again contradicts $\fJ' \boldprecdot \fJ$, and thus $L \subsetneqdot L'$, which means that $\mathrm{(I)}$ holds.

    Otherwise, $z \in \z(\fl)^\reg$.
    Thus $L' = L$, and $\mathcal{O}' = L \cdot e_0 \subseteq \overline{\mathcal{O}}$.
    Since $\J{G}(L;\mathcal{O}') = \fJ' \boldprec \fJ = \J{G}(L;\mathcal{O})$, we must have $\mathcal{O}' \neq \mathcal{O}$, and thus $\mathcal{O}' \boldprec \mathcal{O}$.
    If there exists $\mathcal{O}_0 \in \smallslant{\cN_{\fl}}{L}$ such that $\mathcal{O}' \boldprec \mathcal{O}_0 \boldprec \mathcal{O}$, then \textsc{Lemma}~\ref{Lemma Nilpotent Order Inducing Decomposition Class Order}(ii) implies that $\J{G}\left(L;\mathcal{O}'\right) \boldprec \J{G}\left(L;\mathcal{O}_0\right) \boldprec \J{G}\left(L;\mathcal{O}\right)$, which would contradict $\fJ' \boldprecdot \fJ$.
    Therefore, $\mathcal{O}' \boldprecdot \mathcal{O}$, which means that $\mathrm{(II)}$ holds. \qedhere
\end{proof}

\end{theorem}

Suppose $\fJ, \fJ' \in \fD[G]$ are such that $\fJ' \boldprecdot \fJ$.
If \textsc{Theorem}~\ref{Theorem Hasse Diagram Edges}$\mathrm{(I)}$ holds, then we say that $(\fJ',\fJ)$ is a \emph{semisimple edge} of $\bm{\Gamma}\fD[G]$.
Otherwise, \textsc{Theorem}~\ref{Theorem Hasse Diagram Edges}$\mathrm{(II)}$ holds, and we call $(\fJ',\fJ)$ a \emph{nilpotent edge} of $\bm{\Gamma}\fD[G]$.
The final result of this paper shows that we can determine the type of edge by comparing the stabiliser dimensions of its vertices.

\begin{proposition}\label{Proposition Edge Types and Stabiliser Dimension}
    Suppose $p \geq 0$ is good for $G$, and $x,y \in \g$ are such that $\J{G}y \boldprecdot \J{G}x$.
    \begin{itemize}
        \item[$\mathrm{(i)}$] The edge $\left(\J{G}y,\J{G}x\right)$ is semisimple if and only if $\dim \C_G y = \dim \C_G x$.
        \item[$\mathrm{(ii)}$] The edge $\left(\J{G}y,\J{G}x\right)$ is nilpotent if and only if $\dim \C_G y > \dim \C_G x$.
    \end{itemize} 
\begin{proof}
    Using \textsc{Lemma}~\ref{Lemma Stabiliser/Centraliser Dimension and Closure Order}(a), we know $\dim \C_G y \geq \dim \C_G x$.
    Since each edge of $\bm{\Gamma}\fD[G]$ is either semisimple or nilpotent, it suffices to prove that being a semisimple edge implies that $\dim \C_G y = \dim \C_G x$, and being a nilpotent edge implies that $\dim \C_G y > \dim \C_G x$.

    First, suppose that the edge from $\J{G} y$ to $\J{G} x$ is semisimple.
    Then \textsc{Theorem}~\ref{Theorem Hasse Diagram Edges}(I) implies that there exists $L,L' \in \mathfrak{C}_G$ with $L \subsetneqdot L'$ and $\mathcal{O} \in \smallslant{\cN_{\fl}}{L}$, such that $\J{G} x = \J{G}(L;\mathcal{O})$ and $\J{G} y = \J{G}(L';\Ind_{\fl}^{\fl'}\mathcal{O})$.
    Then \textsc{Lemma}~\ref{Lemma Semisimple Order Inducing Decomposition Class Order}(i) shows that $\J{G}y \subseteq \overline{\J{G}x}^{\,\reg}$, and thus $\dim \C_G y = \dim \C_G x$.

    Otherwise, the edge from $\J{G} y$ to $\J{G} x$ is nilpotent.
    Then \textsc{Theorem}~\ref{Theorem Hasse Diagram Edges}(II) implies that there exists $L \in \mathfrak{C}_G$ and $\mathcal{O},\mathcal{O}' \in \smallslant{\cN_{\fl}}{L}$ with $\mathcal{O}' \boldprecdot \mathcal{O}$, such that $\J{G} x = \J{G}(L;\mathcal{O})$ and $\J{G} y = \J{G}\left(L;\mathcal{O}'\right)$.
    It follows that $\dim \fd_{\g} x\s = \dim \z(\fl) = \dim \fd_{\g} y\s$.
    Since $\dim \J{G}y < \dim \J{G}x$, \textsc{Corollary}~\ref{Corollary Alternative Dimension Formula} implies that $\dim \C_G y > \dim \C_G x$.
\end{proof}
\end{proposition}

\newpage

\printbibliography

\newpage

\appendix

\section{\texorpdfstring{Example: $\mathrm{Sp}_4$-Decomposition Classes}{Example: Sp4-Decomposition Classes}}\label{Appendix: Example}

In this appendix, we shall describe the $\mathrm{Sp}_4$-decomposition classes, and observe the changes that occur in bad characteristic ($p = 2$).
As before, $\mathbb{K}$ is an algebraically closed field of characteristic $p \geq 0$, and all algebraic groups and Lie algebras will be defined over $\mathbb{K}$.

\subsection{Roots}
Let $I_2 \in \mathfrak{gl}_2$ denote the $2 \times 2$ identity matrix, and let $J \coloneqq \begin{psmallmatrix} 0 & I_2 \\ -I_2 & 0 \end{psmallmatrix} \in \mathfrak{gl}_4$.
We then define $G = \mathrm{Sp}_4 \coloneqq \left\{ g \in \mathrm{GL}_4 \midd g^\top J g = J\right\}$, which can be shown to equal \[ \left\{\vbigm \begin{psmallmatrix} A & B \\ C & D \end{psmallmatrix} \in \mathrm{GL}_4 \midd A,B,C,D \in \mathfrak{gl}_2, A^\top C = (A^\top C)^\top, B^\top D = (B^\top D)^\top, A^\top D - C^\top B = I_2 \right\}, \]
where $g^\top$ denotes the transpose of the matrix $g$.
Its corresponding Lie algebra is then $\mathfrak{g} = \mathfrak{sp}_4 \coloneqq \left\{ x \in \mathfrak{gl}_4 \midd x^\top J + Jx = 0 \right\} = \left\{\vbigm \begin{psmallmatrix} A & B \\ C & -A^\top \end{psmallmatrix} \in \mathfrak{gl}_4 \midd A,B,C \in \mathfrak{gl}_2, B = B^\top, C = C^\top \right\}$.
Let $\phi \colon \mathbb{K}^4 \rightarrow \mathfrak{gl}_4$ send a quadruple of scalars to the corresponding diagonal matrix.
Then $T \coloneqq \left\{ h_{s,t} \coloneqq \phi(s,t,s^{-1},t^{-1}) \midd s,t \in \mathbb{K}^{\times} \right\} \subseteq G$ is a maximal torus, with corresponding Lie algebra $\mathfrak{t} = \left\{ z_{s,t} \coloneqq \phi(s,t,-s,-t) \midd s,t \in \mathbb{K} \right\} \subseteq \mathfrak{g}$.
Define $\alpha \in \mathrm{X}(T)$ via $\alpha(h_{s,t}) \coloneqq st^{-1}$, and $\beta \in \mathrm{X}(T)$ via $\beta(h_{s,t}) \coloneqq t^2$.
Then $\Pi \coloneqq \{ \alpha, \beta \}$ is a choice of simple roots, with corresponding positive roots $\Phi^+ \coloneqq \{ \alpha, \beta, \alpha + \beta, 2\alpha + \beta \}$, inside the root system $\Phi \coloneqq \Phi(G,T) = \Phi^+ \sqcup -\Phi^+$.
For each $\gamma \in \Phi$, we define $e_{\gamma} \in \mathfrak{g}$ such that \[\sum_{\gamma \in \Phi} \lambda_{\gamma} e_{\gamma} = \scalemath{1}{\begin{psmallmatrix}
    0 & \lambda_{\alpha} & \lambda_{2\alpha + \beta} & \lambda_{\alpha+\beta} \\
    \lambda_{-\alpha} & 0 & \lambda_{\alpha+\beta} & \lambda_{\beta} \\
    \lambda_{-2\alpha-\beta} & \lambda_{-\alpha-\beta} & 0 & -\lambda_{-\alpha} \\
    \lambda_{-\alpha-\beta} & \lambda_{-\beta} & -\lambda_{\alpha} & 0
\end{psmallmatrix}}, \] for all $\lambda_{\gamma} \in \mathbb{K}$.
The corresponding root subspaces are then $\mathfrak{g}_{\gamma} = \mathbb{K}e_{\gamma}$, for each $\gamma \in \Phi$, with corresponding root subgroups $\mathrm{U}_{\gamma} = I_4 + \mathbb{K}e_{\gamma}$.
For each $\gamma \in \Phi$, we define the isomorphism of algebraic groups $u_{\gamma} \colon \mathbb{G}_{\mathrm{a}} \rightarrow \mathrm{U}_{\gamma}$ via $u_{\gamma}(\xi) \coloneqq I_4 + \xi e_{\gamma}$.
One can verify that the family $(u_{\gamma})_{\gamma \in \Phi}$ is a \emph{realisation} of $\Phi$ in $G$, as defined in \cite[\S 8.1]{S98}.

\subsection{Stabilisers of Semisimple Lie Algebra Elements}

Up to $G$-conjugation, the Levi subgroups of $G$ are $T$, $\left\langle T, \mathrm{U}_{\alpha}, \mathrm{U}_{-\alpha} \right\rangle$, $\left\langle T, \mathrm{U}_{\beta}, \mathrm{U}_{-\beta} \right\rangle$, and $G$ itself.
On the other hand, up to $G$-conjugation, \cite[Lemma 3.7]{S75} shows that each connected stabiliser of a semisimple Lie algebra element is of the form $\Co{G}(z_{s,t}) = \left\langle T, \mathrm{U}_{\gamma} \midd \gamma \in \Phi_{z_{s,t}} \right\rangle$, for some $s,t \in \mathbb{K}$.
Since each $\Phi_{s,t} \coloneqq \Phi_{z_{s,t}} = \left\{ \gamma \in \Phi \midd \mathrm{d}\gamma(z_{s,t}) = 0 \right\}$ is a symmetric subset of roots, it suffices to describe $\Phi^+_{s,t} \coloneqq \Phi^+ \cap \Phi_{s,t}$, for each $s, t \in \mathbb{K}$: we have that $\alpha \in \Phi^+_{s,t}$ if and only if $s-t = 0$; $\beta \in \Phi^+_{s,t}$ if and only if $2t = 0$; $\alpha + \beta \in \Phi^+_{s,t}$ if and only if $s+t = 0$; $2\alpha + \beta \in \Phi^+_{s,t}$ if and only if $2s = 0$.

\subsubsection{}

When $p \neq 2$, it follows that there are are six options for $\Phi^+_{s,t}\colon$
\begin{IEEEeqnarray*}{LLL}
    \bullet~ \text{If } s = t = 0, & \text{ then } & \Phi^+_{s,t} = \Phi^+; \\
    \bullet~ \text{If } s \neq t = 0, & \text{ then } & \Phi^+_{s,t} = \{ \beta \}; \\
    \bullet~ \text{If } t \neq s = 0, & \text{ then } & \Phi^+_{s,t} = \{ 2\alpha + \beta \}; \\
    \bullet~ \text{If } s = t \neq 0, & \text{ then } & \Phi^+_{s,t} = \{ \alpha \}; \\
    \bullet~ \text{If } s = -t \neq 0, & \text{ then } & \Phi^+_{s,t} = \{ \alpha+\beta \}; \\
    \bullet~ \text{Otherwise,} & & \Phi^+_{s,t} = \emptyset.
\end{IEEEeqnarray*}
We note that the second and third options (which correspond to the long roots) result in conjugate Levi subgroups of $G$, as do the fourth and fifth options (which correspond to the short roots).
One can then prove that $\varphi_1 \colon \bigl(a, \begin{smallmatrix} b_1 & b_2 \\ b_3 & b_4 \end{smallmatrix}\bigr) \mapsto \scalemath{0.85}{\begin{psmallmatrix}
    a & 0 & 0 & 0 \\
    0 & b_1 & 0 & b_2 \\
    0 & 0 & a^{-1} & 0 \\
    0 & b_3 & 0 & b_4
\end{psmallmatrix}}$ is an isomorphism of algebraic groups $\mathbb{G}_{\mathrm{m}} \times \mathrm{SL}_2 \rightarrow \Co{G} z_{1,0} = \left\langle T, \mathrm{U}_{\beta}, \mathrm{U}_{-\beta} \right\rangle$.
Hence there are two nilpotent $\Co{G} z_{1,0}$-orbits, of which the non-zero orbit is represented by $\mathrm{d}\varphi_1\bigl(0, \begin{smallmatrix} 0 & 1 \\ 0 & 0 \end{smallmatrix}\bigr) = e_{\beta}$.
Likewise, $\varphi_2 \colon A \mapsto \begin{psmallmatrix} A & 0 \\ 0 & \scalemath{0.65}{(A^\top)^{-1}} \end{psmallmatrix}$ is an isomorphism of algebraic groups $\mathrm{GL}_2 \rightarrow \Co{G} z_{1,1} = \left\langle T, \mathrm{U}_{\alpha}, \mathrm{U}_{-\alpha} \right\rangle$, and there are two nilpotent $\Co{G} z_{1,1}$-orbits, with the non-zero orbit represented by $\mathrm{d}\varphi_2\bigl(\begin{smallmatrix} 0 & 1 \\ 0 & 0 \end{smallmatrix}\bigr) = e_{\alpha}$.

\subsubsection{}

On the other hand, there are only two options for $\Phi^+_{s,t}$ when $p = 2$: if $s = t$, then $\Phi^+_{s,t} = \Phi^+$; otherwise, $\Phi^+_{s,t} = \{ \beta, 2\alpha + \beta \}$, the set of positive long roots.
Therefore, up to conjugation, there are only two connected stabilisers of semisimple Lie algebra elements: $\Co{G} z_{1,0} = \Co{G} z_{s,t} = \left\langle T, \mathrm{U}_{\beta}, \mathrm{U}_{-\beta}, \mathrm{U}_{2\alpha+\beta}, \mathrm{U}_{-2\alpha-\beta} \right\rangle$, for any $s \neq t$, and $G$ itself.
Moreover, any $x \in \mathfrak{g}$ with $x_{\mathrm{s}} \notin \mathfrak{z}(\mathfrak{g})$ is not Levi-type.
Similarly to before, one can prove that $\varphi_3 \colon \bigl(\begin{smallmatrix} a_1 & a_2 \\ a_3 & a_4 \end{smallmatrix}, \begin{smallmatrix} b_1 & b_2 \\ b_3 & b_4 \end{smallmatrix}\bigr) \mapsto \scalemath{0.8}{\begin{psmallmatrix}
    a_1 & 0 & a_2 & 0 \\
    0 & b_1 & 0 & b_2 \\
    a_3 & 0 & a_4 & 0 \\
    0 & b_3 & 0 & b_4
\end{psmallmatrix}}$ is an isomorphism of algebraic groups $\mathrm{SL}_2 \times \mathrm{SL}_2 \rightarrow \Co{G} z_{1,0}$.
It follows that there are four nilpotent $\Co{G} z_{1,0}$-orbits, and the non-zero orbits are represented by $\mathrm{d}\varphi_3(\begin{smallmatrix} 0 & 1 \\ 0 & 0 \end{smallmatrix},0) = e_{2\alpha+\beta}$, $\mathrm{d}\varphi_3(0,\begin{smallmatrix} 0 & 1 \\ 0 & 0 \end{smallmatrix}) = e_{\beta}$, and $\mathrm{d}\varphi_3(\begin{smallmatrix} 0 & 1 \\ 0 & 0 \end{smallmatrix},\begin{smallmatrix} 0 & 1 \\ 0 & 0 \end{smallmatrix}) = e_{\beta} + e_{2\alpha+\beta}$.

\subsection{Nilpotents} 

Let $V \coloneqq \mathbb{K}^4$, and have $\mathfrak{g} = \mathfrak{sp}_4$ act on $V$ via left-multiplication.
Fix $x \in \mathcal{N}_{\mathfrak{g}}$, and recall that a \emph{partition} of $n \in \mathbb{N}^+$ is a tuple of non-increasing positive integers $(\lambda_1,\ldots,\lambda_r)$ such that $\lambda_1 + \cdots + \lambda_r = n$.
Using \cite[\S1.3]{H79}, there exists a unique partition $\bm{\lambda}(x) = (\lambda_1,\ldots,\lambda_r)$ of 4 with the following property: there exists a family of vectors $v_1, \ldots, v_r \in V$ such that $\left\{ x^av_i \midd 1 \leq i \leq r, 0 \leq a \leq \lambda_i\right\}$ is a basis of $V$.
Moreover, $\bm{\lambda}(x)$ has exactly $2 \dim\ker(x^j) - \dim\ker(x^{j+1}) - \dim\ker(x^{j-1})$ parts of size $j$, for each $j \in \mathbb{N}^+$.
The corresponding \emph{index function} is $\chi(x)(m) \coloneqq \min\left\{ n \in \mathbb{N} \midd \forall v \in \ker(x^m) \colon (x^{n+1}v)^\top J x^n v = 0 \right\}$.
Hesselink proved, in \cite[Theorem 3.8]{H79}, that $x,x' \in \mathcal{N}_{\mathfrak{g}}$ are $G$-conjugate if and only if $\bm{\lambda}(x) = \bm{\lambda}(x')$ and $\chi(x) = \chi(x')$.

The \emph{exponent function} corresponding to the partition $\bm{\lambda}(x)$ is the map $n(x) \colon \mathbb{N}^+ \rightarrow \mathbb{N}$ which sends $m \in \mathbb{N}^+$ to the number of parts of $\bm{\lambda}(x)$ which are equal to $m$.
Suppose that $m_1 > \cdots > m_k$ are the elements of $\left\{ m \in \mathbb{N}^+ \midd n(x)(m) \geq 1 \right\}$ listed in decreasing order.
Then the \emph{symbol} corresponding to $x$ is defined as $\mathrm{S}(x) \coloneqq \pigl({m_1}_{\chi(m_1)}^{n(m_1)}, \ldots, {m_k}_{\chi(m_k)}^{n(m_k)}\pigr)$, where $n = n(x)$ and $\chi = \chi(x)$.
Hesselink also proved that $x,x' \in \mathcal{N}_{\mathfrak{g}}$ are $G$-conjugate if and only if $\mathrm{S}(x) = \mathrm{S}(x')$.

\begin{table}[H]
\centering
\begin{tabular}{c||c|rc|c|c|c}
{ $x$}                     & { $\bm{\lambda}(x)$}  & { } $\vpigm$              & { $\chi(x)(2)$} & { $\chi(x)(3)$} & { $\chi(x)(4)$} & { $\mathrm{S}(x)$}                    \\ \hline \hline
$e_{\alpha} +e_{\beta}$            & $(4)$                    &                      & $0$                & $1$                & $2$                & $\left(4_2\right)\vpigm$              \\ \hline
$e_{2\alpha+\beta} + e_{\beta}$ & $(2,2)$                  &                      & $1$                & $1$                & $1$                & $\left(2_1^2\right)\vpigm$      \\ \hline
\multirow{2}{*}{$e_{\alpha}$} & \multirow{2}{*}{$(2,2)$} & if $p \neq 2:$       & $1$                & $1$                & $1$                & $\left(2_1^2\right)\vpigm$      \\ \cline{3-7} 
                              &                          & if $p = 2:$          & $0$                & $0$                & $0$                & $\left(2_0^2\right)\vpigm$      \\ \hline
$e_{\beta}$                   & $(2,1,1)$                &                      & $1$                & $1$                & $1$                & $\left(2_1,1_0^2\right)\vpigm$ \\ \hline
$0$                           & $(1,1,1,1)$              & \multicolumn{1}{c}{} & $0$                & $0$                & $0$                & $\left(1_0^4\right)\vpigm$ \\    
\end{tabular}

\caption{Representatives of Nilpotent $\mathrm{Sp}_4$-Orbits}
\label{Table: The Classification of Nilpotent Sp4-Orbits}
\end{table}

\textsc{Table}~\ref{Table: The Classification of Nilpotent Sp4-Orbits} lists some nilpotent elements of $\mathfrak{g}$, along with $\bm{\lambda}(x)$, the values $\chi(x)(m)$ for $2 \leq m \leq 4$, and the symbol $\mathrm{S}(x)$.
We note that $\chi(x)(0) = \chi(x)(1) = 0$ and $\chi(x)(m) = \chi(x)(4)$ for all $x \in \mathcal{N}_{\mathfrak{g}}$ and $m \geq 4$.

Cross-referencing with \cite[Table 2]{H79}, we can see that \textsc{Table}~\ref{Table: The Classification of Nilpotent Sp4-Orbits} is a complete list of nilpotent $G$-orbit representatives, with the caveat that $G \cdot e_{\alpha} = G \cdot (e_{2\alpha+\beta} + e_{\beta})$ if $p \neq 2$.
Next we shall determine the closure order on $\smallslant{\mathcal{N}_{\mathfrak{g}}}{G}$.

\subsubsection{}

When $p \neq 2$, the closure order coincides with the dominance order on partitions, from which it follows that $\{0\} \boldprecdot G \cdot e_{\beta} \boldprecdot G \cdot e_{\alpha} = G \cdot (e_{2\alpha+\beta}+e_{\beta}) \boldprecdot G \cdot (e_{\alpha} + e_{\beta})$.

\subsubsection{}

For the $p=2$ case, we follow the method found in \cite[\S3]{S82}, and parametrise $\smallslant{\mathcal{N}_{\mathfrak{g}}}{G}$ using pairs of partitions whose total size is $2$.
We shall spare the reader the computational details, and state the outcome; namely, $\{0\} \boldprecdot G \cdot e_{\alpha} \boldprecdot G \cdot (e_{2\alpha+\beta}+e_{\beta}) \boldprecdot G \cdot (e_{\alpha} + e_{\beta})$ and $\{0\} \boldprecdot G \cdot e_{\beta} \boldprecdot G \cdot (e_{2\alpha+\beta}+e_{\beta})$.
This forms the following Hasse diagram for the closure order on $\smallslant{\mathcal{N}_{\mathfrak{sp}_4}}{\mathrm{Sp}_4}$ when $p = 2$, where the heights indicate the dimension of each nilpotent $\mathrm{Sp}_4$-orbit, as calculated using \cite[Table 2]{H79}.
All the edges are directed and oriented upwards, so the arrow heads have not been drawn.
The nilpotent orbits are labelled by their representatives, as given in \textsc{Table}~\ref{Table: The Classification of Nilpotent Sp4-Orbits}.

\begin{figure}[!ht]
    \centering
    \scalebox{1}{
    \begin{tikzpicture}[scale = 1]
            \node [style=white ellipse] (1) at (2, 0) {$e_{\alpha}+e_{\beta}$};
            \node [style=white ellipse] (2) at (2, -2) {$e_{2\alpha+\beta}+e_{\beta}$};
            \node [style=white ellipse] (3) at (0, -4) {$e_{\alpha}$};
            \node [style=white ellipse] (4) at (4, -4) {$e_{\beta}$};
            \node [style=white ellipse] (5) at (2, -7) {$0$};
            \node [style=white ellipse] at (8, 0) {8};
            \node [style=white ellipse] at (8, -2) {6};
            \node [style=white ellipse] at (8, -4) {4};
            \node [style=white ellipse] at (8, -7) {0};
            \draw (1) -- (2) ;
            \draw (2) -- (3) ;
            \draw (2) -- (4) ;
            \draw (3) -- (5) ;
            \draw (4) -- (5) ;
            \begin{scope}[on background layer]
                \foreach \y in {0,-2,-4,-7}
                \draw[gray!25, very thin] (-2,\y) -- (8,\y) ;
            \end{scope}
    \end{tikzpicture}
    }
    \caption{Hasse Diagram for $\smallslant{\mathcal{N}_{\mathfrak{sp}_4}}{\mathrm{Sp}_4}$, when $p = 2$}
    \label{Hasse Diagram Nilpotent Sp4 Orbits p=2}
\end{figure}

\subsection{Decomposition Classes} 

We can therefore determine the $G$-decomposition classes, and calculate their dimensions, stabiliser dimensions, and centraliser dimensions.

\subsubsection{}

For $p \neq 2$, there are nine decomposition classes, and we have chosen a representative for each in the table below.
The representatives are already written in their Jordan decompositions, with a semisimple part of the form $z_{s,t}$ and a nilpotent part $e_0$ on the list in \textsc{Table}~\ref{Table: The Classification of Nilpotent Sp4-Orbits}.
For each representative $z_{s,t} + e_0$, we then list $\dim \C_G(z_{s,t} + e_0)$ which equals $\dim \C_{\Co{G} z_{s,t}}$ by \textsc{Lemma}~\ref{Lemma Stabiliser/Centraliser of Jordan Decomposition}(i), and also equals $\dim \fc_{\mathfrak{g}}(z_{s,t} + e_0)$ by \cite[Lemma 2.6.2(ii)]{L05}; these were either calculated directly, or by using \cite[Table 2]{H79}.
We also list $\dim \fd_{\mathfrak{g}} z_{s,t}$ and $\dim \J{G}(z_{s,t}+e_0)$, as calculated using \textsc{Corollary}~\ref{Corollary Alternative Dimension Formula}.

\begin{table}[H]
\centering
\begin{tabular}{c||c|c|c}
$x \vbigm$                    & $\dim \C_G x$ & $\dim \fd_{\mathfrak{g}} x_{\mathrm{s}}$ & $\dim \J{G}x$  \\ \hline \hline
$z_{2,1} \vbigm$              & $2$           & $2$                         & $10$           \\ \hline 
$z_{1,1}+e_{\alpha} \vbigm$   & $2$           & $1$                         & $9$            \\ \hline
$z_{1,1} \vbigm$              & $4$             & $1$                           & $7$              \\ \hline
$z_{1,0}+e_{\beta} \vbigm$    & $2$             & $1$                           & $9$              \\ \hline
$z_{1,0} \vbigm$              & $4$             & $1$                           & $7$              \\ \hline
$e_{\alpha}+e_{\beta} \vbigm$ & $2$             & $0$                           & $8$              \\ \hline
$e_{\alpha} \vbigm$           & $4$             & $0$                           & $6$              \\ \hline
$e_{\beta} \vbigm$            & $6$             & $0$                           & $4$              \\ \hline
$0 \vbigm$                    & $10$            & $0$                           & $0$             
\end{tabular}
\caption{$\mathrm{Sp}_4$-Decomposition Class Representatives, when $p \neq 2$}
\label{Table: Sp4-Decomposition Class Representatives, when p neq 2}
\end{table}

We can therefore determine that the level sets of $\mathfrak{g}$ are as follows:
\begin{itemize}
    \item $\mathfrak{g}_{(2)} = \mathfrak{g}_{[2]} = \J{G}(z_{2,1}) \sqcup \J{G}(z_{1,1}+e_{\alpha}) \sqcup \J{G}(z_{1,0}+e_{\beta}) \sqcup \J{G}(e_{\alpha} + e_{\beta})$;
    \item $\mathfrak{g}_{(4)} = \mathfrak{g}_{[4]} = \J{G}(z_{1,1}) \sqcup \J{G}(z_{1,0}) \sqcup \J{G}(e_{\alpha})$;
    \item $\mathfrak{g}_{(6)} = \mathfrak{g}_{[6]} = \J{G}(e_{\beta})$;
    \item $\mathfrak{g}_{(10)} = \mathfrak{g}_{[10]} = \J{G}(0) = \{0\}$.
\end{itemize}

\subsubsection{}

One might initially suspect that there are still nine $G$-decomposition classes in the $p = 2$ case, since there are four nilpotent $\Co{G}z_{1,0}$-orbits and five nilpotent $G$-orbits.
However, $g \coloneqq \scalemath{0.7}{\begin{psmallmatrix}
        0 & 1 & 0 & 0 \\
    1 & 0 & 0 & 0 \\
    0 & 0 & 0 & 1 \\
    0 & 0 & 1 & 0
\end{psmallmatrix}} \in G$ is such that $g \cdot \mathrm{d}\varphi_3 \bigl(\begin{smallmatrix} a_1 & a_2 \\ a_3 & a_1 \end{smallmatrix}, \begin{smallmatrix} b_1 & b_2 \\ b_3 & b_1 \end{smallmatrix}\bigr) = \mathrm{d}\varphi_3 \bigl(\begin{smallmatrix} b_1 & b_2 \\ b_3 & b_1 \end{smallmatrix}, \begin{smallmatrix} a_1 & a_2 \\ a_3 & a_1 \end{smallmatrix}\bigr)$, hence $g \cdot (z_{s,t} + e_{2\alpha+\beta}) = z_{t,s} + e_{\beta}$, for all $(s,t) \in \mathbb{K}^2$.
It follows that $\J{G}(z_{1,0} + e_{2\alpha+\beta}) = \J{G}(z_{1,0} + e_{\beta})$, and so there are at most eight $G$-decomposition classes in the $p = 2$ case.

\begin{table}[!ht]
\centering
\begin{tabular}{c||c|c|c|c}
$x \vbigm$                    & $\dim \C_G x$ & $\dim \fc_{\mathfrak{g}} x$ & $\dim \fd_{\mathfrak{g}} x_{\mathrm{s}}$ & $\dim \J{G}x$  \\ \hline \hline
$z_{1,0} + e_{2\alpha+\beta} + e_{\beta} \vbigm$              & $2$ & $4$           & $2$                         & $10$           \\ \hline 
$z_{1,0} + e_{\beta} \vbigm$              & $4$ & $5$             & $2$                           & $8$              \\ \hline
$z_{1,0} \vbigm$    & $6$ & $6$            & $2$                           & $6$              \\ \hline
$e_{\alpha} + e_{\beta} \vbigm$              & $2$ & $4$             & $1$                           & $9$              \\ \hline
$e_{2\alpha+\beta} + e_{\beta} \vbigm$ & $4$ & $6$            & $1$                           & $7$              \\ \hline
$e_{\alpha} \vbigm$           & $6$ & $6$             & $1$                           & $5$              \\ \hline
$e_{\beta} \vbigm$            & $6$ & $7$            & $1$                           & $5$              \\ \hline
$0 \vbigm$                    & $10$ & $10$           & $1$                           & $1$             
\end{tabular}
\caption{$\mathrm{Sp}_4$-Decomposition Class Representatives, when $p = 2$}
\label{Table: Sp4-Decomposition Class Representatives, when p = 2}
\end{table}

Once again, we have chosen representatives of the form $z_{s,t} + e_0$ and calculated both $\dim \C_G(z_{s,t}+e_0)$ and $\dim \fc_{\mathfrak{g}}(z_{s,t}+e_0)$ for each, which we note do not always coincide when $p = 2$.
The values for the top three rows in \textsc{Table}~\ref{Table: Sp4-Decomposition Class Representatives, when p = 2} were determined using the isomorphism $\mathrm{d}\varphi_3 \colon \mathfrak{sl}_2 \oplus \mathfrak{sl}_2 \rightarrow \fc_{\mathfrak{g}}z_{1,0}$ alongside \textsc{Lemma}~\ref{Lemma Stabiliser/Centraliser of Jordan Decomposition}, and the others were determined using \cite[Table 2]{H79}.
We then observe that $\fd_{\mathfrak{g}}z_{1,0} = \mathrm{d}\varphi_3(\mathfrak{z}(\mathfrak{sl}_2), \mathfrak{z}(\mathfrak{sl}_2)) = \mathrm{d}\varphi_3(\mathbb{K}I_2, \mathbb{K}I_2)$ and $\fd_{\mathfrak{g}}z_{0,0} = \mathfrak{z}(\mathfrak{g}) = \mathbb{K}I_4$.
This allows us to calculate $\dim \fd_{\mathfrak{g}} z_{s,t}$, before finally using \textsc{Corollary}~\ref{Corollary Alternative Dimension Formula} to determine $\dim \J{G}(z_{s,t}+e_0)$.

Since no two rows in \textsc{Table}~\ref{Table: Sp4-Decomposition Class Representatives, when p = 2} contain the same values, this confirms that all eight of our representatives lie in distinct decomposition classes.
We observe that, unlike the $p \neq 2$ case, the maximal decomposition class $\J{G}^{\mathrm{max}} = \J{G}(z_{1,0} + e_{2\alpha+\beta} + e_{\beta})$ does not consist of regular semisimple elements.
This is equivalent to the observation that the maximal torus $T \subseteq G$ is not a stabiliser-type Levi subgroup, which happens precisely when $p \mkern 1.5mu \mathrm{X}(T) \cap \Phi$ is non-empty.

We can also determine that the stabiliser level sets of $\mathfrak{g}$ are as follows:
\begin{itemize}
    \item $\mathfrak{g}_{(2)} = \J{G}(z_{1,0} + e_{2\alpha+\beta}+e_{\beta}) \sqcup \J{G}(e_{\alpha}+e_{\beta})$;
    \item $\mathfrak{g}_{(4)} = \J{G}(z_{1,0}+e_{\beta}) \sqcup \J{G}(e_{2\alpha+\beta}+e_{\beta})$;
    \item $\mathfrak{g}_{(6)} = \J{G}(z_{1,0}) \sqcup \J{G}(e_{\alpha}) \sqcup \J{G}(e_{\beta})$;
    \item $\mathfrak{g}_{(10)} = \J{G}(0) = \mathfrak{z}(\mathfrak{g})$.
\end{itemize}
In contrast to the $p \neq 2$ case, the centraliser level sets of $\mathfrak{g}$ are (in general) distinct from the stabiliser level sets of $\mathfrak{g}$ when $p = 2$.
\begin{itemize}
    \item $\mathfrak{g}_{[4]} = \J{G}(z_{1,0} + e_{2\alpha+\beta}+e_{\beta}) \sqcup \J{G}(e_{\alpha}+e_{\beta})$;
    \item $\mathfrak{g}_{[5]} = \J{G}(z_{1,0}+e_{\beta})$;
    \item $\mathfrak{g}_{[6]} = \J{G}(z_{1,0}) \sqcup \J{G}(e_{2\alpha+\beta}+e_{\beta}) \sqcup \J{G}(e_{\alpha})$;
    \item $\mathfrak{g}_{[7]} = \J{G}(e_{\beta})$;
    \item $\mathfrak{g}_{[10]} = \J{G}(0) = \mathfrak{z}(\mathfrak{g})$.
\end{itemize}

\subsection{The Hasse Diagram}

We shall soon construct the Hasse diagram for the closure order on $\mathfrak{D}[\mathrm{Sp}_4]$, in both cases.
We first observe the following topological result, the proof of which is routine, and so shall be omitted.

\begin{lemma}\label{Lemma Limits of Morphisms}
    Suppose that $\psi \colon \mathbb{K} \rightarrow Y$ is a morphism of affine varieties, and that $X \subseteq Y$ is any subset.
    If $\psi(\xi) \in X$ for all but finitely many $\xi \in \mathbb{K}$, then $\psi(\mathbb{K}) \subseteq \overline{X}$.
\iffalse
\begin{proof}
    All morphisms of affine varieties are continuous with respect to the Zariski topology, hence $\psi(\overline{Z}) \subseteq \overline{\psi(Z)}$ for any subset $Z \subseteq \mathbb{K}$.
    Since $\mathbb{K}$ is an infinite field, the Zariski topology on $\mathbb{K}$ coincides with the cofinite topology.
    
    Suppose that $\psi(\xi) \in X$ for all but finitely many $\xi \in \mathbb{K}$, and let $Z \coloneqq \left\{ \xi \in \mathbb{K} \midd \psi(\xi) \in X \right\}$.
    Then $\overline{Z} = \mathbb{K}$, hence $\psi(\mathbb{K}) \subseteq \overline{\psi(Z)} \subseteq \overline{X}$.
\end{proof}
\fi
\end{lemma}

\begin{example}\label{Example Same Semisimple Part}
    For all characteristics $p \geq 0$, we have that $z_{1,0} + e_{\beta} \in \overline{\J{G}(z_{1,0} + e_{2\alpha+\beta}+e_{\beta})}$.
\begin{proof}
    Let $x \coloneqq z_{1,0} + e_{2\alpha+\beta}+e_{\beta} \in \mathfrak{g}$.
    First, suppose that $p \neq 2$, and let $g \coloneqq h_{1,-1} + \frac{1}{2}e_{2\alpha+\beta} \in \mathrm{GL}_4$.
    One can compute that $g^\top J g = J$ and $g \cdot x = z_{1,0} + e_{\beta}$, from which the result is clear.

    Otherwise, suppose that $p = 2$, and recall the isomorphism of algebraic groups \linebreak$\varphi_3 \colon \mathrm{SL}_2 \times \mathrm{SL}_2 \rightarrow \Co{G} z_{1,0}$ from earlier.
    Define a morphism $\psi \colon \mathbb{K} \rightarrow \mathfrak{g}$ via $\psi(\xi) \coloneqq z_{1,0} + \xi^2 e_{2\alpha+\beta} + e_{\beta} = \mathrm{d}\varphi_3\bigl(\begin{smallmatrix} 1 & \xi^2 \\ 0 & 1 \end{smallmatrix},\begin{smallmatrix} 0 & 1 \\ 0 & 0 \end{smallmatrix}\bigr)$.
    For each $\xi \neq 0$, $h_{\xi^{-1},1} \in G$ is such that $h_{\xi^{-1},1} \cdot \psi(\xi) = \varphi_3\bigl(\begin{smallmatrix} \xi^{-1} & 0 \\ 0 & \xi \end{smallmatrix},\begin{smallmatrix} 1 & 0 \\ 0 & 1 \end{smallmatrix}\bigr) \cdot \mathrm{d}\varphi_3\bigl(\begin{smallmatrix} 1 & \xi^2 \\ 0 & 1 \end{smallmatrix},\begin{smallmatrix} 0 & 1 \\ 0 & 0 \end{smallmatrix}\bigr) = \mathrm{d}\varphi_3(\begin{smallmatrix} 1 & 1 \\ 0 & 1 \end{smallmatrix},\begin{smallmatrix} 0 & 1 \\ 0 & 0 \end{smallmatrix}) = x$.
    Therefore, $\psi(\xi) \in \J{G}(x)$ for all $\xi \neq 0$, so \textsc{Lemma}~\ref{Lemma Limits of Morphisms} implies that $z_{1,0} + e_{\beta} = \psi(0) \in \overline{\J{G}(x)}$.
\end{proof}
\end{example}

We next record the following general property of \emph{nilpotent decomposition classes} -- defined as the decomposition classes of nilpotent elements -- which holds for all connected reductive algebraic groups.

\begin{lemma}\label{Lemma Nilpotent Decomposition Class Containment}
    Suppose that $H$ is a connected reductive algebraic group, over an algebraically closed field of any characteristic.
    Let $e_0 \in \mathcal{N}_{\mathfrak{h}}$ be any nilpotent element, and let $\mathfrak{J} \in \mathfrak{D}[H]$ be any $H$-decomposition class.
    Then $\J{H}(e_0) \boldpreceq \mathfrak{J}$ if and only if $e_0 \in \overline{\mathfrak{J}}$.

    Moreover, if $\mathfrak{J} = \J{H}(e_1)$ with $e_1 \in \mathcal{N}_{\mathfrak{h}}$, then $\overline{\mathfrak{J}}$ is precisely the union of $\J{H}(e_0)$ across all $e_0 \in \overline{H \cdot e_1}$.
\begin{proof}
    If $\J{H}(e_0) \boldpreceq \mathfrak{J}$ -- that is, $\J{H}(e_0) \subseteq \overline{\mathfrak{J}}$ -- then \textsc{Theorem}~\ref{Theorem Initial Decomposition Class Properties}(v) immediately shows that $e_0 \in \overline{\mathfrak{J}}$.
    Conversely, suppose that $e_0 \in \overline{\mathfrak{J}}$.
    Recall that $\overline{\mathfrak{J}}$ is $H$-stable and, using \textsc{Corollary}~\ref{Corollary Decomposition Varieties Properties}(i), also stable under addition by elements of $\mathfrak{z}(\mathfrak{h})$.
    Therefore, $H \cdot e_0 + \mathfrak{z}(\mathfrak{h}) \subseteq \overline{\mathfrak{J}}$, hence the main result follows from \textsc{Theorem}~\ref{Theorem Initial Decomposition Class Properties}(v).

    Suppose that $\mathfrak{J} = \J{H}(e_1)$, for some nilpotent $e_1 \in \mathcal{N}_{\mathfrak{h}}$.
    Then $\mathfrak{J} = \J{H}(H;e_1)$ is a Levi-type decomposition class, so it follows from \textsc{Proposition}~\ref{Proposition Levi-type Decomposition Variety as Union of Decomposition Classes} that $\overline{\mathfrak{J}}$ is the union of the decomposition classes it intersects.
    If $x \in \overline{\mathfrak{J}}$, then \textsc{Corollary}~\ref{Corollary Decomposition Varieties Properties}(ii) implies that $x_{\mathrm{s}} \in \mathfrak{z}(\mathfrak{h})$ and $x_{\mathrm{n}} \in \overline{H \cdot e_1}$, hence $\J{H}(x) = \J{H}(x_{\mathrm{n}})$.
    Therefore, the decomposition classes that $\overline{\mathfrak{J}}$ intersects are precisely those of the form $\J{H}(e_0)$ for $e_0 \in \overline{H \cdot e_1}$.
\end{proof}
\end{lemma}

The following result regarding decomposition classes of the form $\J{G}(z_{1,0}+e_1)$ will only be used in the $p = 2$ case, but it is easy to prove that it also holds when $p \neq 2$.

\begin{lemma}\label{Lemma Non-Nilpotent Sp4-Decomposition Classes}
    Let $e_0, e_1 \in \{0,e_{\beta},e_{2\alpha+\beta}+e_{\beta}\}$.
\begin{itemize}
    \item[$\mathrm{(i)}$] $\overline{\J{G}(z_{1,0}+e_1)}$ is a union of decomposition classes.
    \item[$\mathrm{(ii)}$] $\J{G}(z_{1,0} + e_0) \subseteq \overline{\J{G}(z_{1,0}+e_1)}$ if and only if $\dim \J{G}(z_{1,0} + e_0) \leq \dim \J{G}(z_{1,0} + e_1)$.
\end{itemize}
\begin{proof}
    First, suppose that $p \neq 2$.
    Then $\J{G}(z_{1,0}+e_1)$ is necessarily Levi-type, and so (i) follows from \textsc{Proposition}~\ref{Proposition Levi-type Decomposition Variety as Union of Decomposition Classes}.
    Using the argument from the proof of \textsc{Example}~\ref{Example Same Semisimple Part}, we know that $\J{G}(z_{1,0}+e_{\beta}) = \J{G}(z_{1,0}+e_{2\alpha+\beta} + e_{\beta})$.
    Therefore, (ii) follows from \textsc{Table}~\ref{Table: Sp4-Decomposition Class Representatives, when p neq 2}, alongside the observation from \textsc{Figure}~\ref{Hasse Diagram Sp4 Decomposition Classes p neq 2} that $\J{G}(z_{1,0}) \subseteq \overline{\J{G}(z_{1,0}+e_{\beta})}$.

    Now, suppose that $p = 2$.
    Let $\mathfrak{J} \coloneqq \J{G}(z_{1,0}+e_1)$, and suppose that $x \in \overline{\mathfrak{J}}$.
    Since $\overline{\mathfrak{J}}$ is $G$-stable, we can assume (without loss of generality) that $x_{\mathrm{s}} \in \mathfrak{t}$ and is thus of the form $z_{s,t}$ for some $s,t \in \mathbb{K}$.
    If $s = t$, then $x_{\mathrm{s}} \in \mathfrak{z}(\mathfrak{g})$, so $x_{\mathrm{n}} \in \overline{\mathfrak{J}}$, and the same argument as in \textsc{Lemma}~\ref{Lemma Nilpotent Decomposition Class Containment} shows that $\J{G}x = \J{G} x_{\mathrm{n}} \subseteq \overline{\mathfrak{J}}$.

    Otherwise, $(s,t) \in \mathbb{K}^2_{\neq}$.
    Fix any $(a,b) \in \mathbb{K}^2_{\neq}$, and let $c \coloneqq \frac{a+b}{s+t} + b \in \mathbb{K}$.
    Since $\overline{\mathfrak{J}}$ is $\mathbb{K}$-stable, and stable under addition by elements of $\mathfrak{z}(\mathfrak{g})$ -- see \textsc{Corollary}~\ref{Corollary Decomposition Varieties Properties}(i) -- it follows that $z_{a,b} + \frac{a+b}{s+t} x_{\mathrm{n}} = \frac{a+b}{s+t} x + z_{c,c} \in \overline{\mathfrak{J}}$.
    Using \cite[Lemma 2.10]{J04}, there exists $h \in \Co{G}z_{a,b}$ such that $h \cdot x_{\mathrm{n}} = \frac{s+t}{a+b} x_{\mathrm{n}}$.
    Since $\overline{\mathfrak{J}}$ is $G$-stable, it follows that $z_{a,b} + x_{\mathrm{n}} = h \cdot \left(z_{a,b} + \frac{a+b}{s+t} x_{\mathrm{n}}\right) \in \overline{\mathfrak{J}}$.
    Therefore, $\fd_{\mathfrak{g}}^{\mathrm{reg}}x_{\mathrm{s}} + x_{\mathrm{n}} = \left\{ z_{a,b} + x_{\mathrm{n}} \midd (a,b) \in \mathbb{K}^2_{\neq} \right\} \subseteq \overline{\mathfrak{J}}$, and so (i) follows from \textsc{Theorem}~\ref{Theorem Initial Decomposition Class Properties}(i).

    For (ii), the forwards implication is given by \textsc{Corollary}~\ref{Corollary Closure Order for Decomposition Classes}, so suppose that \linebreak$\dim \J{G}(z_{1,0} + e_0) \leq \dim \J{G}(z_{1,0} + e_1)$.
    Using \cite[Proposition 2.11]{J04}, we know that $z_{1,0} \in \overline{G \cdot (z_{1,0} + e_1)}$, from which it follows that $z_{1,0} \in \overline{\J{G}(z_{1,0} + e_1)}$.
    Therefore, (i) implies that $\J{G}(z_{1,0}) \subseteq \overline{\J{G}(z_{1,0} + e_1)}$.
    Using \textsc{Table}~\ref{Table: Sp4-Decomposition Class Representatives, when p = 2}, we are in (at least) one of the following cases: $e_0 = e_1$; $e_0 = 0$; or $e_0 = e_{\beta}$ and $e_1 = e_{2\alpha+\beta} + e_{\beta}$.
    If $e_0 = e_1$, then $\J{G}(z_{1,0} + e_0) \subseteq \overline{\J{G}(z_{1,0}+e_1)}$ is vacuously true, and if $e_0 = 0$, then we have already shown that $\J{G}(z_{1,0}+e_0) \subseteq \overline{\J{G}(z_{1,0} + e_1)}$.
    Otherwise, $e_0 = e_{\beta}$ and $e_1 = e_{2\alpha+\beta} + e_{\beta}$, in which case \textsc{Example}~\ref{Example Same Semisimple Part} shows that $z_{1,0} + e_0 \in \overline{\J{G}(z_{1,0} + e_1)}$, hence the result follows from (i).
\end{proof}
\end{lemma}

\subsubsection{}

Now, assume that $p \neq 2$.
Using the results in \cite[\S7.3]{CM93} (which we note, hold in all good characteristics), we can determine the following Lusztig--Spaltenstein inductions:
\begin{itemize}
    \item $\Ind_{\mathfrak{c}^{\vphantom{\circ}}_{\mathfrak{g}}z^{\vphantom{\circ}}_{1,1}}^{\mathfrak{g}} (\Co{G}z_{1,1} \cdot e_{\alpha}) = \Ind_{\mathfrak{c}^{\vphantom{\circ}}_{\mathfrak{g}}z^{\vphantom{\circ}}_{1,0}}^{\mathfrak{g}} (\Co{G}z_{1,0} \cdot e_{\beta}) = G \cdot (e_{\alpha} + e_{\beta})$;
    \item $\Ind_{\mathfrak{c}^{\vphantom{\circ}}_{\mathfrak{g}}z^{\vphantom{\circ}}_{1,1}}^{\mathfrak{g}} \{0\} = \Ind_{\mathfrak{c}^{\vphantom{\circ}}_{\mathfrak{g}}z^{\vphantom{\circ}}_{1,0}}^{\mathfrak{g}} \{0\} = G \cdot e_{\alpha}$.
\end{itemize}

We can then use \textsc{Theorem}~\ref{Theorem Hasse Diagram Edges} to determine the Hasse diagram of the $\mathrm{Sp}_4$-decomposition classes, which we present below.
As with \textsc{Figure}~\ref{Hasse Diagram Nilpotent Sp4 Orbits p=2}, all the edges are oriented upwards, and so their arrow heads have not been drawn.
We have recorded the dimension of each decomposition class from \textsc{Table}~\ref{Table: Sp4-Decomposition Class Representatives, when p neq 2}, as indicated by the labels to the right and the thin grey lines at the relevant heights.
Moreover, since $p \neq 2$ is good for $\mathrm{Sp}_4$, \textsc{Theorem}~\ref{Theorem Hasse Diagram Edges} shows that every edge is either semisimple or nilpotent.
We have drawn the semisimple edges with solid lines, and the nilpotent edges with dashed lines.
The decomposition classes themselves are labelled with their representatives, as given in \textsc{Table}~\ref{Table: Sp4-Decomposition Class Representatives, when p neq 2}.

\begin{figure}[!ht]
    \centering
    \scalebox{0.72}{
    \begin{tikzpicture}[scale = 0.82]
            \node [style=white ellipse] (A1) at (0, 0) {$z_{2,1}$};
            \node [style=white ellipse] (B1) at (4, -2) {$z_{1,1}+e_{\alpha}$};
            \node [style=white ellipse] (B2) at (4, -6) {$z_{1,1}$};
            \node [style=white ellipse] (C1) at (8, -2) {$z_{1,0}+e_{\beta}$};
            \node [style=white ellipse] (C2) at (8, -6) {$z_{1,0}$};
            \node [style=white ellipse] (D1) at (12, -4) {$e_{\alpha}+e_{\beta}$};
            \node [style=white ellipse] (D2) at (12, -8) {$e_{\alpha}$};
            \node [style=white ellipse] (D3) at (12, -11) {$e_{\beta}$};
            \node [style=white ellipse] (D4) at (12, -16) {$0$};
            \node [style=white ellipse] at (16, 0) {10};
            \node [style=white ellipse] at (16, -2) {9};
            \node [style=white ellipse] at (16, -4) {8};
            \node [style=white ellipse] at (16, -6) {7};
            \node [style=white ellipse] at (16, -8) {6};
            \node [style=white ellipse] at (16, -11) {4};
            \node [style=white ellipse] at (16, -16) {0};
            \draw (A1) -- (B1) ;
            \draw (A1) -- (C1) ;
            \draw (B1) -- (D1) ;
            \draw (C1) -- (D1) ;
            \draw (B2) -- (D2) ;
            \draw (C2) -- (D2) ;
            \draw[dashed] (B1) -- (B2) ;
            \draw[dashed] (C1) -- (C2) ;
            \draw[dashed] (D1) -- (D2) ;
            \draw[dashed] (D2) -- (D3) ;
            \draw[dashed] (D3) -- (D4) ;
            \begin{scope}[on background layer]
                \foreach \y in {0,-2,...,-8,-11,-16}
                \draw[gray!25, very thin] (-2,\y) -- (16,\y) ;
            \end{scope}
    \end{tikzpicture}
    }
    \caption{Hasse Diagram for $\mathfrak{D}[\mathrm{Sp}_4]$, when $p \neq 2$}
    \label{Hasse Diagram Sp4 Decomposition Classes p neq 2}
\end{figure}

\subsubsection{}

Determining the edges of the Hasse diagram is much harder when $p = 2$, mainly since we cannot use \textsc{Theorem}~\ref{Theorem Hasse Diagram Edges}.
We have to rely on ad-hoc methods to determine whether $\mathfrak{J}_1 \subseteq \overline{\mathfrak{J}_2}$, for decomposition classes $\mathfrak{J}_1,\mathfrak{J}_2 \in \mathfrak{D}[G]$.

Using \textsc{Lemma}~\ref{Lemma Nilpotent Decomposition Class Containment}, we know that the edges between nilpotent $G$-decomposition classes are completely determined by the closure order on nilpotent $G$-orbits.
It follows from \textsc{Figure}~\ref{Hasse Diagram Nilpotent Sp4 Orbits p=2} that such edges correspond precisely to the covering relations on this list:
\begin{itemize}
    \item $\J{G}(0) \boldprecdot \J{G} (e_{\alpha}) \boldprecdot \J{G} (e_{2\alpha+\beta}+e_{\beta}) \boldprecdot \J{G} (e_{\alpha} + e_{\beta})$;
    \item $\J{G}(0) \boldprecdot \J{G} (e_{\beta}) \boldprecdot \J{G} (e_{2\alpha+\beta}+e_{\beta})$.
\end{itemize}

Combining \textsc{Lemma}~\ref{Lemma Non-Nilpotent Sp4-Decomposition Classes} with \textsc{Table}~\ref{Table: Sp4-Decomposition Class Representatives, when p = 2}, we can conclude that the only edges between pairs of non-nilpotent decomposition classes correspond precisely to the covering relations $\J{G}(z_{1,0}) \boldprecdot \J{G} (z_{1,0} + e_{\beta}) \boldprecdot \J{G} (z_{1,0} + e_{2\alpha+\beta} + e_{\beta})$.
Therefore, it remains to determine all covering relations of the form $\J{G}(e_0) \boldprecdot \J{G}(z_{1,0}+e_1)$.

\begin{lemma}\label{Lemma Closure of Semisimple Elements}
    Suppose that $p = 2$.
    Then $\overline{\J{G}(z_{1,0})} = \overline{G \cdot \mathfrak{t}} = \mathfrak{t} \oplus \mathfrak{g}_{\alpha} \oplus \mathfrak{g}_{-\alpha} \oplus \mathfrak{g}_{\alpha+\beta} \oplus \mathfrak{g}_{-\alpha-\beta}$.
\begin{proof}
    Since $G \cdot \mathfrak{t}$ is precisely the set of all semisimple elements of $\mathfrak{g}$, and $\J{G}(0) = \mathfrak{z}(\mathfrak{g})$, it follows that $\J{G}(z_{1,0}) = G \cdot \mathfrak{t} \setminus \mathfrak{z}(\mathfrak{g})$.
    Using \textsc{Corollary}~\ref{Corollary Decomposition Varieties Properties}(i), we know that $\mathfrak{z}(\mathfrak{g}) \subseteq \overline{\J{G}(z_{1,0})}$.
    Therefore, $G \cdot \mathfrak{t} \subseteq \overline{\J{G}(z_{1,0})} \subseteq \overline{G \cdot \mathfrak{t}}$, from which the first equality in the result follows.

    Recall that the family of maps $u_{\gamma} \colon \xi \mapsto I_4 + \xi e_{\gamma}$, for each $\gamma \in \Phi$, defines a realisation of $\Phi$ in $G$.
    One can verify via case-by-case computation that, since $p = 2$, we have that $u_{\gamma}(\xi) \cdot e_{\delta} \in \mathfrak{h} \coloneqq \mathfrak{t} \oplus \mathfrak{g}_{\alpha} \oplus \mathfrak{g}_{-\alpha} \oplus \mathfrak{g}_{\alpha+\beta} \oplus \mathfrak{g}_{-\alpha-\beta}$, for all $\gamma \in \Phi$, $\xi \in \mathbb{K}$, and $\delta \in \{\alpha,-\alpha,\linebreak\alpha+\beta,-\alpha-\beta\}$.
    Using \cite[\S2.2.2]{L05}, one can also verify that $u_{\gamma}(\xi) \cdot z \in \mathfrak{h}$, for each $z \in \mathfrak{t}$.
    Since $\mathfrak{h}$ is clearly $T$-stable, it is therefore stable under the adjoint action of $G = \left\langle T, \mathrm{U}_{\gamma} \midd \gamma \in \Phi \right\rangle$.
    Since $\mathfrak{t} \subseteq \mathfrak{h}$, it follows that $G \cdot \mathfrak{t} \subseteq \mathfrak{h}$.
    Finally, we observe -- using \textsc{Table}~\ref{Table: Sp4-Decomposition Class Representatives, when p = 2} -- that $\overline{\J{G}(z_{1,0})}$ and $\mathfrak{t} \oplus \mathfrak{g}_{\alpha} \oplus \mathfrak{g}_{-\alpha} \oplus \mathfrak{g}_{\alpha+\beta} \oplus \mathfrak{g}_{-\alpha-\beta}$ are both closed irreducible sets of dimension $6$, from which the second equality follows.
\end{proof}
\end{lemma}

\begin{remark}
    \textsc{Lemma}~\ref{Lemma Closure of Semisimple Elements} demonstrates that the semisimple Lie algebra elements do not form a dense subset of $\mathfrak{g}$.
    As before, this is a consequence of the fact that $p \mkern 1.5mu \mathrm{X}(T) \cap \Phi \neq \emptyset$, and hence the maximal torus $T \subseteq G$ is not a stabiliser-type Levi subgroup.
\end{remark}

Using \textsc{Lemma}~\ref{Lemma Closure of Semisimple Elements}, the only Hasse diagram edge with terminal vertex $\J{G}(z_{1,0})$ must begin at $\J{G}(e_{\alpha})$.
Likewise, since the maximal decomposition class $\J{G}^{\mathrm{max}} = \J{G}(z_{1,0} + e_{2\alpha+\beta} + e_{\beta})$ contains all other decomposition classes in its closure, the only Hasse diagram edge with terminal vertex $\J{G}^{\mathrm{max}}$ must begin at the maximal nilpotent decomposition class $\J{G}(e_{\alpha} + e_{\beta})$.
Therefore, using \textsc{Lemma}~\ref{Lemma Nilpotent Decomposition Class Containment}, it remains to determine which nilpotents $e_0 \in \left\{ e_{\alpha},e_{\beta},e_{2\alpha+\beta}+e_{\beta},e_{\alpha}+e_{\beta} \right\}$ satisfy $e_0 \in \overline{\J{G}(z_{1,0} + e_{\beta})}$.

Consider $\psi \colon \mathbb{K} \rightarrow \mathfrak{g}$ defined via $\psi(\xi) \coloneqq \scalemath{0.75}{\begin{psmallmatrix}
        0 & \xi^2 & 1 & 0 \\
        0 & \xi^2 & 0 & 1 \\
        0 & 0 & 0 & 0 \\
        0 & 0 & \xi^2 & \xi^2
    \end{psmallmatrix}}$, and define $g_{\xi} \coloneqq \scalemath{0.75}{\begin{psmallmatrix}
        \xi^{-1} & 1 & 0 & 0 \\
        0 & 1 & \xi^{-1} & \xi^{-2} \\
        0 & 0 & \xi & 0 \\
        0 & 0 & \xi & 1
    \end{psmallmatrix}} \in \mathfrak{gl}_4$, for each $\xi \neq 0$.
Since each entry of $\psi(\xi)$ is a polynomial in $\mathbb{K}[\xi]$, it follows that $\psi \colon \mathbb{K} \rightarrow \mathfrak{g}$ is a well-defined morphism of affine varieties.
It is simple to compute that ${g_{\xi}}^\top J g_{\xi} = J$, so $g_{\xi} \in G$ for each $\xi \neq 0$.
Moreover, one can compute that $g_{\xi}(z_{0,\xi^2}+e_{\beta}){g_{\xi}}^{-1} = \psi(\xi)$, hence $\psi(\xi) \in G \cdot (z_{0,\xi^2}+e_{\beta}) \subseteq \J{G}(z_{0,\xi^2}+e_{\beta})$, for each $\xi \neq 0$.
Since $\Co{G}(z_{0,\xi^2}) = \Co{G}(z_{1,0})$, for each $\xi \neq 0$, we have that $\psi(\xi) \in \J{G}(z_{1,0} + e_{\beta})$.
It then follows from \textsc{Lemma}~\ref{Lemma Limits of Morphisms} that $e_{2\alpha+\beta} + e_{\beta} = \psi(0) \in \overline{\J{G}(z_{1,0}+e_{\beta})}$.

\begin{remark}
    Our definition of $\psi \colon \mathbb{K} \rightarrow \mathfrak{g}$ provides little to no insight into how one might arrive at such a morphism.
    As with much of the results in this subsection, this was the result of an ad-hoc process, involving a lot of trial and error and computation.
    We hope to later develop a more concrete method for determining Hasse diagram edges for $\mathfrak{D}[G]$ in bad characteristics, akin to the classification in \textsc{Theorem}~\ref{Theorem Hasse Diagram Edges}.
\end{remark}

By comparing dimensions using \textsc{Table}~\ref{Table: Sp4-Decomposition Class Representatives, when p = 2}, we know from \textsc{Corollary}~\ref{Corollary Closure Order for Decomposition Classes} that $e_{\alpha} + e_{\beta} \notin \overline{\J{G}(z_{1,0} + e_{\beta})}$.
Since \textsc{Figure}~\ref{Hasse Diagram Nilpotent Sp4 Orbits p=2} shows that $e_{\alpha},e_{\beta} \in \overline{G \cdot (e_{2\alpha+\beta} + e_{\beta})}$, it follows that the only Hasse diagram edge with terminal vertex $\J{G}(z_{1,0}+e_{\beta})$ must begin at $\J{G}(e_{2\alpha+\beta} + e_{\beta})$.

Therefore, we can draw the Hasse diagram for $\mathfrak{D}[G]$ using (mostly) the same conventions as \textsc{Figure}~\ref{Hasse Diagram Sp4 Decomposition Classes p neq 2} above.
The only difference is that we can no longer classify the Hasse diagram edges using \textsc{Theorem}~\ref{Theorem Hasse Diagram Edges}.
We have decided to draw the edges corresponding to minimal nilpotent coverings with dashed lines, and the other edges with solid lines.

\begin{figure}[!ht]
    \centering
    \scalebox{0.7}{
    \begin{tikzpicture}[scale = 0.9]
            \node [style=white ellipse] (A1) at (0, 0) {$z_{1,0}+e_{2\alpha+\beta}+e_{\beta}$};
            \node [style=white ellipse] (A2) at (0, -4) {$z_{1,0}+e_{\beta}$};
            \node [style=white ellipse] (A3) at (0, -8) {$z_{1,0}$};
            \node [style=white ellipse] (B1) at (6, -2) {$e_{\alpha}+e_{\beta}$};
            \node [style=white ellipse] (B2) at (6, -6) {$e_{2\alpha+\beta}+e_{\beta}$};
            \node [style=white ellipse] (B3) at (6, -10) {$e_{\alpha}$};
            \node [style=white ellipse] (B4) at (9, -10) {$e_{\beta}$};
            \node [style=white ellipse] (B5) at (9, -14) {$0$};
            \node [style=white ellipse] at (13, 0) {10};
            \node [style=white ellipse] at (13, -2) {9};
            \node [style=white ellipse] at (13, -4) {8};
            \node [style=white ellipse] at (13, -6) {7};
            \node [style=white ellipse] at (13, -8) {6};
            \node [style=white ellipse] at (13, -10) {5};
            \node [style=white ellipse] at (13, -14) {1};
            \draw (A1) -- (B1) ;
            \draw (A2) -- (B2) ;
            \draw (A3) -- (B3) ;
            \draw[dashed] (A1) -- (A2) ;
            \draw[dashed] (A2) -- (A3) ;
            \draw[dashed] (B1) -- (B2) ;
            \draw[dashed] (B2) -- (B3) ;
            \draw[dashed] (B2) -- (B4) ;
            \draw[dashed] (B3) -- (B5) ;
            \draw[dashed] (B4) -- (B5) ;
            \begin{scope}[on background layer]
                \foreach \y in {0,-2,...,-10,-14}
                \draw[gray!25, very thin] (-2,\y) -- (13,\y) ;
            \end{scope}
    \end{tikzpicture}
    }
    \caption{Hasse Diagram for $\mathfrak{D}[\mathrm{Sp}_4]$, when $p = 2$}
\end{figure}

Finally, we shall record certain properties of $\mathrm{Sp}_4$-decomposition classes that we can deduce are true in all characteristics.
In particular, this provides examples of bad characteristic situations where the assumptions in \textsc{Corollary}~\ref{Corollary Locally Closed Decomposition Classes} and \textsc{Proposition}~\ref{Proposition Coincide with Irreducible Component iff Isolated} still hold.

\begin{theorem}
    Suppose that $\mathfrak{J} \in \mathfrak{D}[\mathrm{Sp}_4]$, for any characteristic $p \geq 0$.
\begin{itemize}
    \item[$\mathrm{(i)}$] $\overline{\mathfrak{J}}$, $\overline{\mathfrak{J}}^{\,\mathrm{reg}}$, and $\overline{\mathfrak{J}}^{\,\mathfrak{g}\text{-}\mathrm{reg}}$ are each unions of decomposition classes.
    \item[$\mathrm{(ii)}$] $\mathfrak{J}$ is locally closed.
    \item[$\mathrm{(iii)}$] Suppose that $\mathfrak{g}_{\langle m \rangle}$ is a level set of $\mathfrak{g}$ such that $\mathfrak{J} \subseteq \mathfrak{g}_{\langle m \rangle}$.
    Then $\mathfrak{J}$ coincides with a sheet of $\mathfrak{g}_{\langle m \rangle}$ if and only if it is isolated in $\mathfrak{g}_{\langle m \rangle}$.
\end{itemize}
\begin{proof}
    First, assume that we are in good characteristic: $p \neq 2$.
    Then $\mathfrak{J}$ is necessarily Levi-type, so (i) follows from \textsc{Proposition}~\ref{Proposition Levi-type Decomposition Variety as Union of Decomposition Classes} and the subsequent paragraph.
    Similarly, (ii) is given by \textsc{Corollary}~\ref{Corollary Good Characteristic Decomposition Classes are Locally Closed}.
    Since every sheet of $\mathfrak{g}$ is also Levi-type, (iii) follows from \textsc{Proposition}~\ref{Proposition Coincide with Levi-Type Sheets iff Isolated}.

    Now, assume that we are in bad characteristic ($p = 2$), and suppose that $\mathfrak{J} = \J{G}(e_1)$ for some nilpotent $e_1 \in \mathcal{N}_{\mathfrak{g}}$.
    Then $\mathfrak{J} = \J{G}(G;e_1)$ is a Levi-type decomposition class, so (i) follows from \textsc{Proposition}~\ref{Proposition Levi-type Decomposition Variety as Union of Decomposition Classes}, as before.
    Otherwise, $\mathfrak{J} = \J{G}(z_{1,0} + e_1)$ for some $e_1 \in \left\{ 0, e_{\beta}, e_{2\alpha + \beta} + e_{\beta} \right\}$.
    Using \textsc{Lemma}~\ref{Lemma Non-Nilpotent Sp4-Decomposition Classes}, we know that $\overline{\mathfrak{J}}$ is a union of decomposition classes, and the rest of (i) follows from \textsc{Lemma}~\ref{Lemma Closure inside a Level Set} alongside \textsc{Lemma}~\ref{Lemma Constant Dimension}.
    
    Suppose that $x \in \overline{\mathfrak{J}}$, and observe that (i) implies that $\overline{\J{G}x}$ is also a union of decomposition classes.
    Therefore, (ii) follows from \textsc{Corollary}~\ref{Corollary Locally Closed Decomposition Classes}.

    Suppose that $\mathfrak{g}_{\langle m \rangle}$ is a level set of $\mathfrak{g}$.
    Using \textsc{Corollary}~\ref{Corollary Irreducible Components of Level Sets}(ii) and the proof of \textsc{Proposition}~\ref{Proposition Irreducible Components and Maximal Elements}(ii), the map $\mathfrak{J}' \mapsto \overline{\mathfrak{J}'} \cap \mathfrak{g}_{\langle m \rangle}$ is a bijection between the maximal decomposition classes in $\mathfrak{g}_{\langle m \rangle}$ and the sheets of $\mathfrak{g}_{\langle m \rangle}$.
    Therefore, (i) implies that every sheet of $\mathfrak{g}_{\langle m \rangle}$ is a union of decomposition classes, hence (iii) follows from \textsc{Proposition}~\ref{Proposition Coincide with Irreducible Component iff Isolated}.
\end{proof}
\end{theorem}

\end{document}